\documentclass {book}

\usepackage{verbatim}
\usepackage{graphicx}

\newtheorem{theorem}{Theorem}[chapter]
\newtheorem{lemma}{Lemma}[chapter]

\newtheorem{corollary}{Corollary}[chapter]

\font \m=msbm10
\newcommand{\R}{{\hbox {\m R}}}
\newcommand{\C}{{\hbox {\m C}}}
\newcommand{\Z}{{\hbox {\m Z}}}

\newcommand{\U}{{\hbox {\m U}}}
\newcommand {\HH}{{\hbox {\m H}}}
\newcommand {\eps} {{\varepsilon}}
\newcommand{\CL}{{\mathcal L}}
\newcommand {\WH}{{\mathcal W}}
\newcommand {\BL}{{\mathcal B}}

\newtheorem {conjecture}[lemma]{Conjecture}
\newtheorem {prop}[lemma]{Proposition}
\def \qed {{}}

\begin{document}
\author{Wendelin Werner}
\date{Universit\'e Paris-Sud and Ecole Normale Sup\'erieure\footnote {Supported by Agence Nationale pour la Recherche, ANR-BLAN-0058}}
\title{Lectures on two-dimensional critical percolation}


\maketitle

\tableofcontents

\mainmatter


%

\section*{Overview}
This school is designed for graduate students in mathematics
specializing in probability theory and statistical mechanics.  This
course consists of six lectures, plus four exercise sessions (two
exercise sheets) that are taken care of by Pierre Nolin and G\'abor Pete.

During the first half of this summer school, Greg Lawler has lectured
on Schramm-Loewner Evolutions (SLE), its definition and
properties. So, I expect the students/readers to have been exposed to
SLE before.

The main goal of my course is to explain the relation between lattice
models and these two-dimensional random curves. In fact, I choose to
focus on critical site percolation on the triangular lattice, which is
the model where the entire picture is now clear and rigorous: One can
start with the lattice-model, study its properties, prove its
conformal invariance and deduce the relation to SLE and finally use
the properties of SLE to derive the results about discrete percolation
that one is aiming for. This combines several steps and ideas, and I
will try to describe them one by one, hopefully with sufficient
details so that students will have a global view by the end of the
lecture series. For other lattice models and other lattices, conformal
invariance and the canvas that I will present do conjecturally hold
too, but in general, at least some mathematical steps are missing (for
instance, to my knowledge, scaling relations are proved to hold only
for percolation) even though a lot of important progress is being
made, such as Stas Smirnov's proof of conformal invariance of the
Ising model \cite {Sm2,Sm3,Sm4}.

Focusing on a specific model allows also to go into some technical
details.  Of course, the counterpart (and I apologize for this!) is
that I will not mention a lot of recent developments, such as the
relation with the Gaussian Free Field, Conformal Loop Ensembles,
questions related to noise etc.

We will try to make these lectures as self-contained as possible.
References on which these lectures are based are the following:
\begin {itemize}
\item For the proof of Cardy-Smirnov formula: Smirnov's self-contained original paper \cite {Sm}.
\item For the convergence to SLE(6), we follow the strategy outlined by Smirnov in \cite {Sm2} (see also our (we=Lawler-Schramm-Werner) paper on loop-erased random walks \cite {LSWlesl}).
\item We will follow some SLE(6) computations from \cite {LSW2,LSW5}.
\item For how to deduce the critical exponents for critical percolation from the above, we follow \cite {SmW, LSW5}.
\item For how to deduce the exponents associated to near-critical percolation, we use Kesten's hyperscaling relations proved in \cite {Ke}, following the presentation of the upcoming review paper \cite {No}.
\end {itemize}

The required background is:
\begin {itemize}
\item
Some probability theory knowledge: The content of the beginning of Grimmett's book \cite{G} on percolation, and some basic stochastic calculus (It\^o's formula, semi-martingales).
\item
Very basic complex analysis knowledge, Riemann's mapping theorem, see e.g. \cite {Ah}
\item
The definition of SLE: This can be found in Oded Schramm's original paper \cite {S1}, in Greg Lawler's notes \cite {Lln}, or in other introductory material such as \cite {Lbook, Wln}.
Note that we will use almost no general property of SLE here. We will for instance get for free the fact that SLE(6) is almost surely a curve by proving that it is the scaling limit of discrete percolation interfaces.
\end {itemize}

I hope at some point to extend these notes into a book on the subject.

\medbreak

The nice pictures in these lecture notes are not mine (I produced only the sketchy basic ones): Many thanks to Julien Dub\'edat, Pierre Nolin and Oded Schramm for allowing me to copy and use their postscript files!

I would also like to thank the organizers of the summer school for making it such a success, and Pierre Nolin and G\'abor Pete for proofreading the exercise sheets. I also thank Christophe Garban and Jeff Steif for their comments and all those who have pointed out misprints  in a former version of these notes (but I suspect that there are quite a few ones left...).

\chapter{Introduction and tightness}

\section {2D percolation}

Let me start with the definition of the percolation model.
We consider a regular planar periodic lattice -- soon we will be focusing on the triangular lattice where each face is an equilateral triangle.
For each site of the considered
infinite lattice, toss a coin independently: With probability $p$ the site is declared open, and with probability $1-p$, it is declared vacant. Then, one is interested in the connectivity properties of the set of open sites.
One can consider the ``connected components'' of open sites: An open connected component $C$ is a connected family of open sites, such that all the sites neighboring $C$ are closed.

For nice planar lattices,
one can prove (e.g. \cite {G}) that there exists a (lattice-dependent) value $p_c$ called the critical probability such that:
\begin {itemize}
\item If $p \le p_c$, then with probability one, there is no infinite connected component of open sites.
\item If $p > p_c$, then with probability one, there is a unique infinite connected component of open sites.
\end {itemize}

In fact, even though there is also no infinite cluster at $p=p_c$, the behavior is quite different than when $p < p_c$. Indeed, the probability $\pi(n)$  that a given fixed site in the lattice is in a cluster of diameter at least $n$ can be shown to decay exponentially fast as soon as $p< p_c$
(i.e. for each $p < p_c$, there exists $\psi (p) > 0 $ such that $\pi (n) \le \exp (-n \psi(p))$ for all large $n$), while it is believed that
when $p=p_c$,
\begin {equation}
\label {e00}
\pi (n) = n^{-5/48 + o (1)} \hbox { as } n \to \infty .
\end {equation}

When the lattice is transitive (i.e. all sites are ``equivalent'', the graph looks the same seen from any site) then one can define the probability $\theta (p)$ that a given site is in the infinite connected open component.
It is not very difficult to prove that $\theta (p)$ is continuous on $[p_c, 1]$ (and therefore on $[0,1]$ as soon as  $\theta (p_c) = 0$ since $\theta=0$ on $[0, p_c)$), see \cite {G} for details.

One way to interpret $\theta (p)$ is to say that it is the ``density'' of the infinite cluster $C$. Indeed, one can prove that for a given increasing nested sequence of sets $\Lambda_n$,
${\# (\Lambda_n \cap C) } / ({\# \Lambda_n}) \to  \theta (p)$
almost surely (see the first exercise sheet).

\medbreak

It is conjectured that for all transitive planar periodic lattices, the way this density vanishes as $p$ tends to the critical value $p_c$ from above is the following:
\begin {equation}
\label {e000}
\theta (p) = (p-p_c)^{5/36 + o(1)}.
\end {equation}

\medbreak

The main goal of these lectures is to explain to you  in an almost self-contained way the different steps that lead to the proof of the conjecture (\ref {e000})  in the special case where the lattice is the standard triangular lattice. We shall also derive (\ref {e00}) on the way.
As you will see, complex analysis and SLE will be instrumental.

\medbreak

One should emphasize that for other planar lattices, these conjectures are at present (i.e. in 2007) not proved.

\section {Notations and prerequisites}
We are from now on going to suppose that
the lattice is the regular triangular lattice $T$, where faces are equilateral triangles.
Until Lecture 6, we will also suppose that the probability $p$ is equal to $1/2$.

\begin{figure}
\centerline{\includegraphics*[height=1.2in]{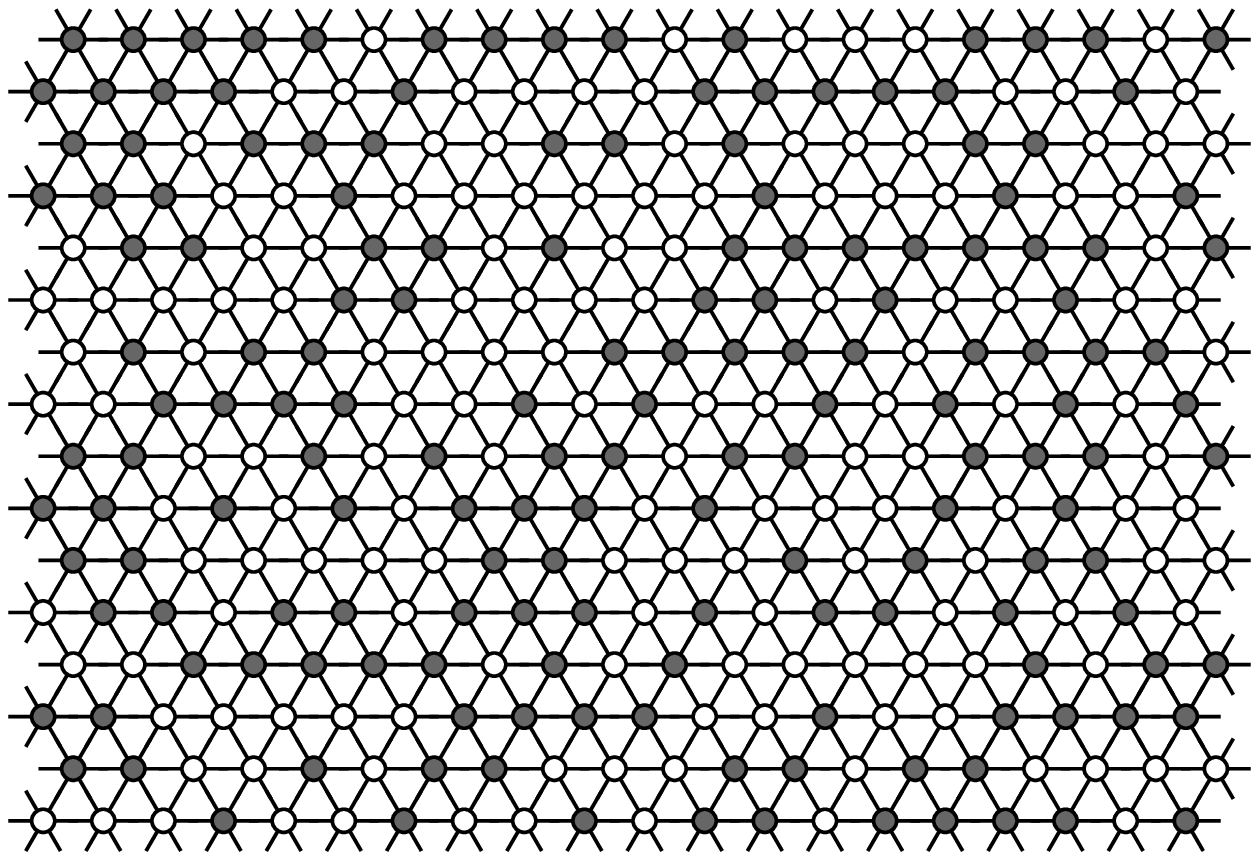}}
\centerline{\includegraphics*[height=1.2in]{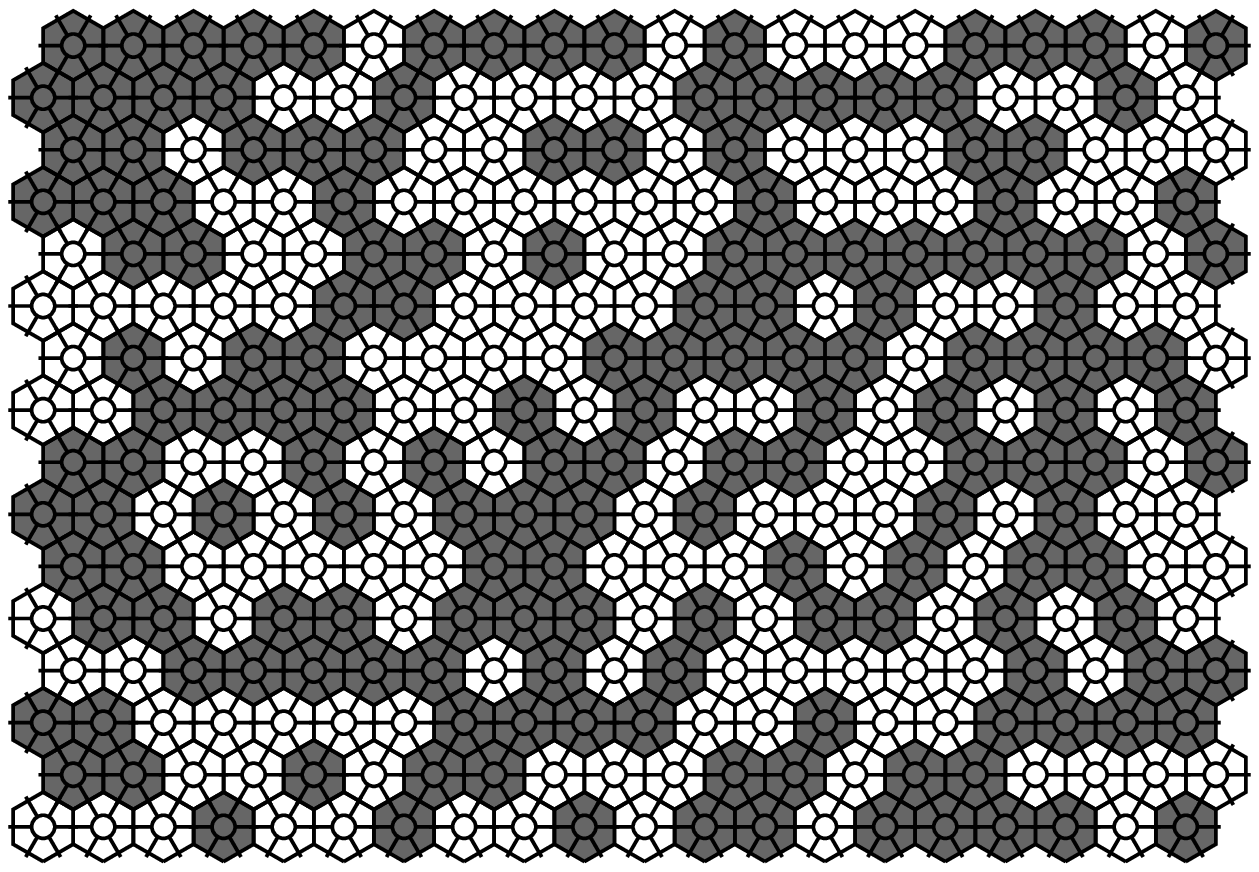}}
\centerline{\includegraphics*[height=1.2in]{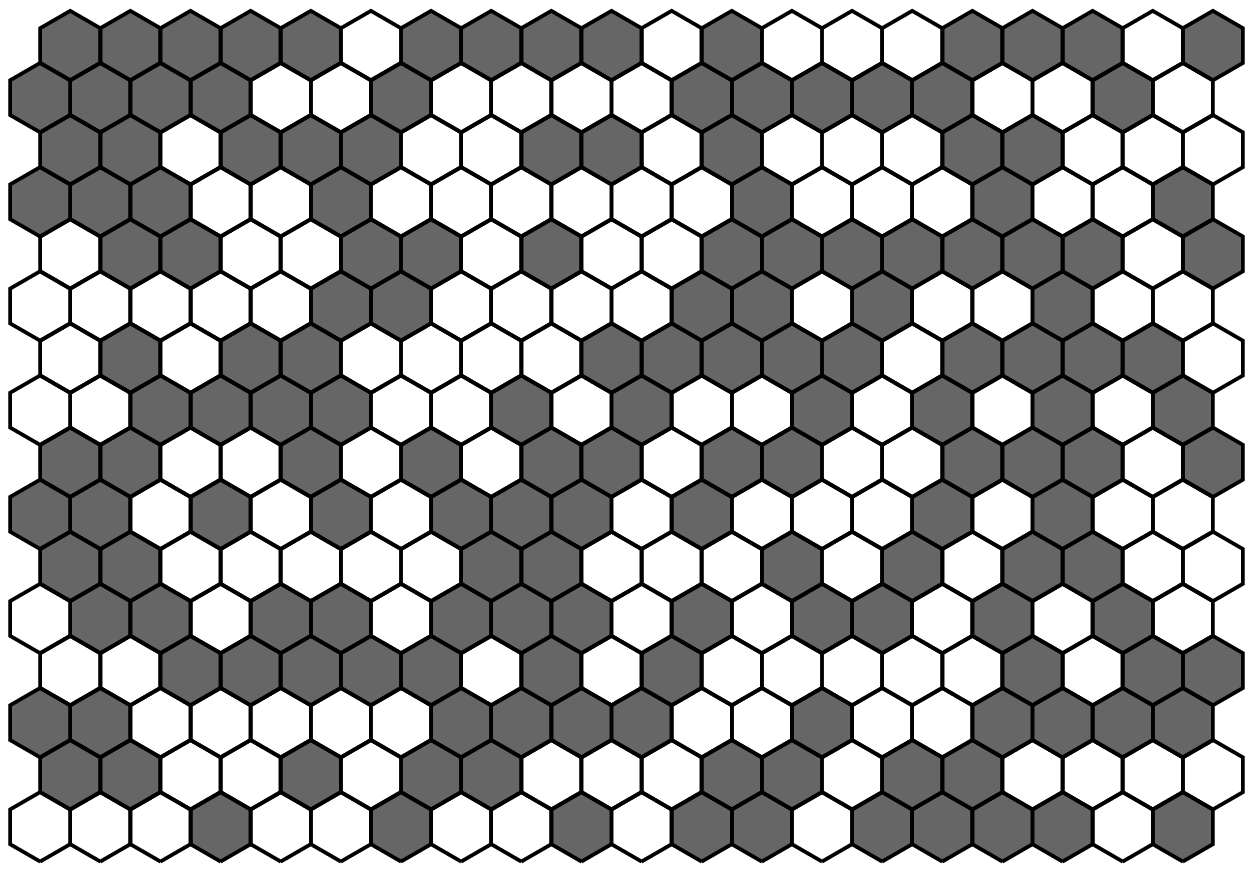}}
\caption {From sites to cells (picture by Oded Schramm)}
\label {fsc}
\end {figure}

Note that coloring the sites of the triangular lattice in black and white is equivalent to coloring the cells (or faces) of a honeycomb lattice as in Figure \ref {fsc}. Each site of the triangular lattice is the center of an hexagonal cell on the honeycomb lattice, and it is connected to each of the centers of the six neighboring hexagons. This representation is more convenient for our eyes to detect crossings and interfaces.

We will sometimes call sites/cells open or occupied. In these lectures, we say that black=open=occupied and white=closed=vacant.

 A point with cartesian coordinates $x$ and $y$ will be denoted by the complex number $x+i y$
(we will use the notation $(x,y)$ for other purposes).
The axis of our triangular lattice are in the (complex) directions $1$, $\exp(i\pi/3)$ and $\tau = \exp ( 2i \pi/3)$. We define $T$ in such a way that its sites are the points in $\Z + e^{i\pi/3} \Z$, so that neighboring sites are at distance one from each other in the complex plane.

We will also use the coordinates defined via the vectors $1$ and $e^{i\pi/3}$: The point $(u,v)$ will denote the point $u + v e^{i \pi /3}$.
Hence, $[0,n] \times [0, n]$ will denote a rhombus with side-length $n$.
More generally, when $I= [a,b]$ and $I'=[a', b']$ are two intervals on the real line, $I\times I'$ will denote a parallelogram, with two horizontal sides of length $b-a$ and two sides of length $b'-a'$ that are parallel to the direction $e^{i\pi/3}$.
We say that a nearest-neighbor path in $I \times I'$ is a top-to-bottom (resp. horizontal) crossing of $I \times I'$ if it joins the two horizontal boundary segments (resp. the two other boundary segments).

\medbreak

An important result due to Kesten and Wierman is that for this triangular lattice, $p_c= 1/2$. This builds on symmetry properties of the model but requires also additional arguments, see \cite {G}. In fact, we will not use the fact that $1/2$ is the critical probability until Lecture 5, but we will make extensive use of the symmetries when $p$ is equal to $1/2$.

\begin{figure}
\label {rhombus}
\centerline{\includegraphics*[height=2in]{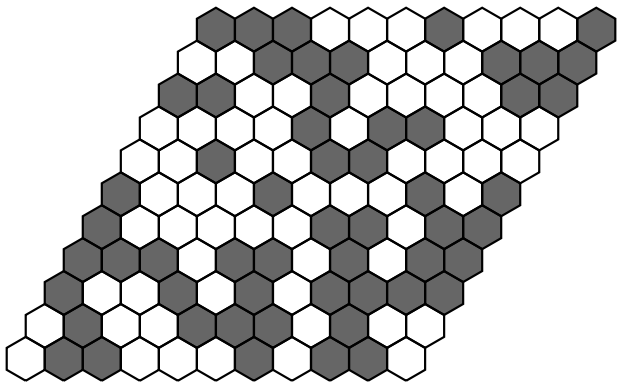}}
\caption{White top-to-bottom crossing vs. black horizontal crossing}
\end {figure}

One key-observation is for instance that when one cuts out a rhombus in our lattice as in Figure \ref {rhombus}, then the probability of an open horizontal crossing (i.e. that there exists a horizontal crossing consisting only of open sites) is identical to the probability of a top-to-bottom closed crossing. Indeed, consider the set $S$ of all white cells that are connected to the top boundary. Either it touches the bottom boundary of the rhombus (and then, a left-to-right crossing by black hexagons can not exist), or it does not touch the bottom part of the rhombus, and in the latter case, this means that the ``lower boundary'' of this set $S$ contains a black left-to-right crossing of the rhombus (in fact this crossing is the ``highest'' possible such black crossing).

Hence, exactly one of the two events ``there exists a white top-to-bottom crossing'' and ``there exists a black left-to-right crossing'' occurs (and this is due to the fact that the lattice is triangular).
But since $p=1/2$, it is clear that for symmetry reasons, these two events have the same probability. This probability is therefore equal to $1/2$.
Note that this is true regardless of the size of the rhombus, and that it therefore already gives information on the large-scale properties of percolation which is our subject of study.

Note that this argument uses only the symmetry of the rhombus, and that it works also for any other shape symmetric to a given axis.

\medbreak

We now recall a classical result called Harris' inequality: We say that an event $A$ is increasing if for any two realizations $w$ and $w'$ of percolation on a given lattice such that any site that is open for $w$ is also open for $w'$ (one often says $w \le w'$) and such that $w \in A$, then $w' \in A$.
In other words, as soon as $w \in A$, then $A$ contains all the realizations that are obtained from $w$ by changing some sites from closed to open.
An example of an increasing event is the existence of a left-to-right open crossing in a rhombus.

Harris' inequality for percolation that is often referred to as the FKG inequality states that if $A$ and $B$ are two increasing events, then
$ P ( A \cap B ) \ge P( A) \times P(B)$.
In other words, if $P(B) > 0$, then $P (A | B) \ge P(A)$; knowing that $B$ holds helps $A$ to hold too.
This inequality can be proved first for finite graphs inductively on the number of sites, and then one can pass to the infinite lattice limit (see e.g. \cite {G} or prove it yourself as a homework exercise).

For example, a consequence of Harris' inequality is that for any given rhombus, the probability that there exists an open left-to-right crossing {\em and} an open top-to-bottom crossing is at least $1/2 \times 1/2 = 1/4$.

\section {Russo-Seymour-Welsh}

I now present a short proof of the Russo-Seymour-Welsh (RSW) estimates
for this particular lattice that I learned from Stas Smirnov (as most
of what I will present in the first three lectures). Let us
define ``rectangles'' in the triangular lattice. For any positive
integer $a$ and any real positive $b$, $R(a,b)$ is the set of sites $x+iy$ in
the triangular lattice such that $0 \le x \le a$ and $0 \le y \le b$. Recall that if we use the
tilted integer coordinates $(u,v)$, then  $x = u + v/2$ and $y = v \sqrt {3}
/ 2$.

We denote by $H(a,b)$ the event that there exists a horizontal open
crossing of $R(a,b)$ (i.e. there exists a nearest-neighbor path of
black sites in the rectangle that joins a point with $x$-coordinate
equal to $0$ to a point with $x$-coordinate equal to $a$). Note that
$a \mapsto P (H(a,b))$ is clearly a non-increasing function of $a$.

\begin {lemma}[RSW]
 One always has $P (H(2a, b)) \ge P (H(a,b))^2 / 4$.
\end {lemma}

Before proving this fact, note that a horizontal crossing of the
rhombus $[-n, n] \times [0,2n]$ must contain a horizontal crossing of
the rectangle $R(n, n \sqrt {3} )$.  Hence, $P ( H (n, n \sqrt {3}))
\ge 1/2$, and if we combine this with the lemma, we get that\linebreak $P (
H(2n,n \sqrt {3})) \ge 1/16$, and $P ( H (4n, n \sqrt {3})) \ge
1/2^{10}$.  By induction, we therefore get readily that:
\begin {corollary}
For any $k$, there exists $a_k>0$ such that
$ P ( H (kn, n )) \ge a_k $ for any $n>2$.
\end {corollary}
The actual value of $a_k$ is not so important. What is crucial is that $a_k >0$ and that this inequality holds simultaneously for all $n$.
In fact, in the sequel, we will only use the fact that for any $n$,
$P (H (4n, n) ) \ge a_4$.
\medbreak

\noindent
{\bf Warm-up.}
Discrete percolation is just a way to represent a sequence of coin-tosses. The clue to understand the connectivity properties is always to ``organize'' this randomness, and to explore the state of the different sites in a suitable way.

We want to compare $H(2a, b)$ and $H(a, b)$.
Clearly, in order to cross the rectangle $R(2a,  b)$, one has first to cross the left half-rectangle $R(a,b)$.
 One way to discover if there is a left-right open crossing of $R(a,b)$, is ``from top to bottom''. Let us give some details here, since this will be a recurrent theme in the lectures:
Imagine that one pours water on the top side of the rectangle and that water can only flow along white sites/cells. If there is a black left-to-right crossing, it will prevent the water from reaching the real axis. Let $\CL$ denote the set of sites filled with water; this is just the ``cluster'' (for percolation in the rectangle) of white sites that are connected to the top side via white cells. For each possible such fixed set $L$, the event $\CL=L$ means that the sites of $L$ are white, and that all sites in the rectangle that are at distance $1$ from $L$ are black.
 Hence, the event $\CL=L$
is measurable with respect to the state of the sites in the union of $L$ and of its neighboring sites. When $H(a,b)$ holds, then the ``lower'' boundary of $\CL$ contains a simple (i.e. without loops) path $\gamma$  consisting of black cells only. It is easy to check that this is the ``highest'' such crossing (we hope that the precise definition of ``highest'' is clear enough from the context).
Hence, we see that for each such simple path, the event $\gamma = g$ is independent of the state of sites ``below'' $g$. Another way to see this, is to note that $\gamma =g$ means that the sites of $g$ are black, and that the sites ``above'' $g$ are colored in such a way that there is no ``higher'' black crossing.

Yet another way to construct $\gamma$ is to start an ``exploration process'' (as we shall describe later) from the top-left corner that leaves black cells to its right and white cells to its left, and to stop this process whenever it hits the right side of the rectangle. The lower boundary of this exploration path is $\gamma$.

\medbreak
\noindent
{\bf Proof of the lemma.}
For each deterministic simple left-to-right crossing $g$ of the rectangle $R(a,b)$, we denote by $\tilde g$ its symmetric path with respect to the line $I_a= \{z \ : \Re (z) = a \} = a + i \R$, and we call $O$ the connected component of $R(2a, b) \setminus ( g \cup \tilde g)$ that contains $a$ (when $a \in g$ and $O = \emptyset$ it is simple to see that our argument will be valid as well). Note that $O$ is a simply connected set of sites that is symmetric with respect to $I_a$ and that the point $a$ is the middle of the bottom side of $R(2a, b)$.
Let us first suppose that $g$ does not touch the real axis. Then, the boundary of $O$ (i.e. the set of points that are distance $1$ from $O$) can be divided into four parts: $g$, $\tilde g$, a part $J$ that is on the ``left-and-bottom'' boundary of $R(a,b)$ and its symmetric image with respect to $I_a$, that we call $\tilde J$.
The same symmetry argument as in the case of the rhombus (color $g$ and $\tilde J$ in black and $\tilde g \cup J$ in white, and sample percolation in $O$ and look which colors manages to create a crossing) shows that the probability of a black crossing in $O$ joining $g$ to $\tilde J$ is $1/2$, regardless of our choice of $g$. We call this event $A(g)$.
Note that the events $\{\gamma= g \}$ and $A(g)$ are independent because the former is independent of the state of the cells below $g$. Furthermore, if these two events hold then there exists a black path in $R(2a,b)$  that joins the left boundary of this rectangle to the union of its right boundary with the right-half of its bottom boundary. We call $A'$ this event.

When $g$ touches the real axis, then $\tilde J$ just consists of a part of the lower-right boundary, but the same symmetry argument still shows that
$P(A(g)) =1/2$ and that $P (A' | \gamma = g ) \ge 1/2$.

\begin{figure}
\centerline{\includegraphics*[height=1.3in]{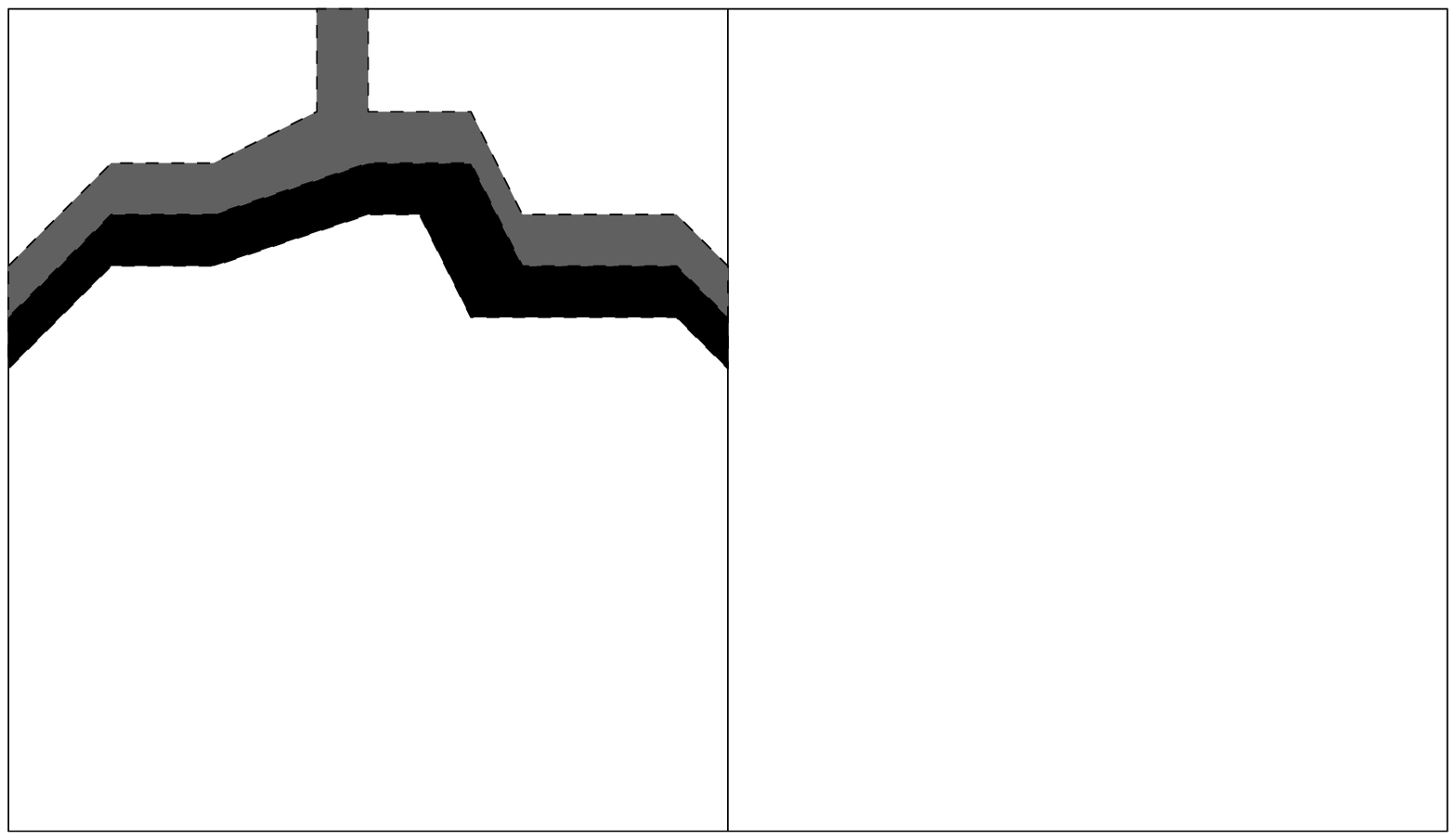}}
\centerline{\includegraphics*[height=1.3in]{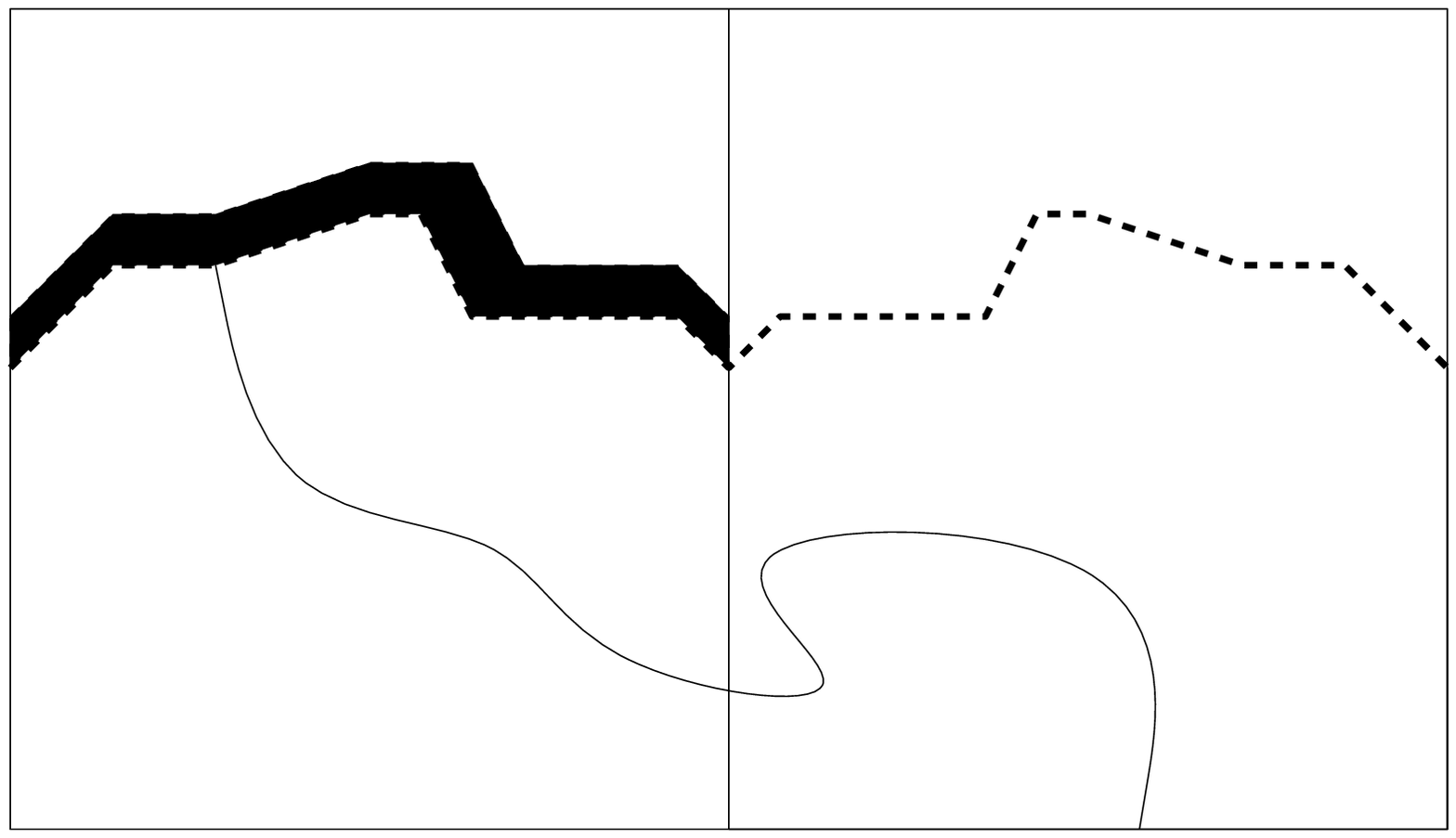}}
\centerline{\includegraphics*[height=1.3in]{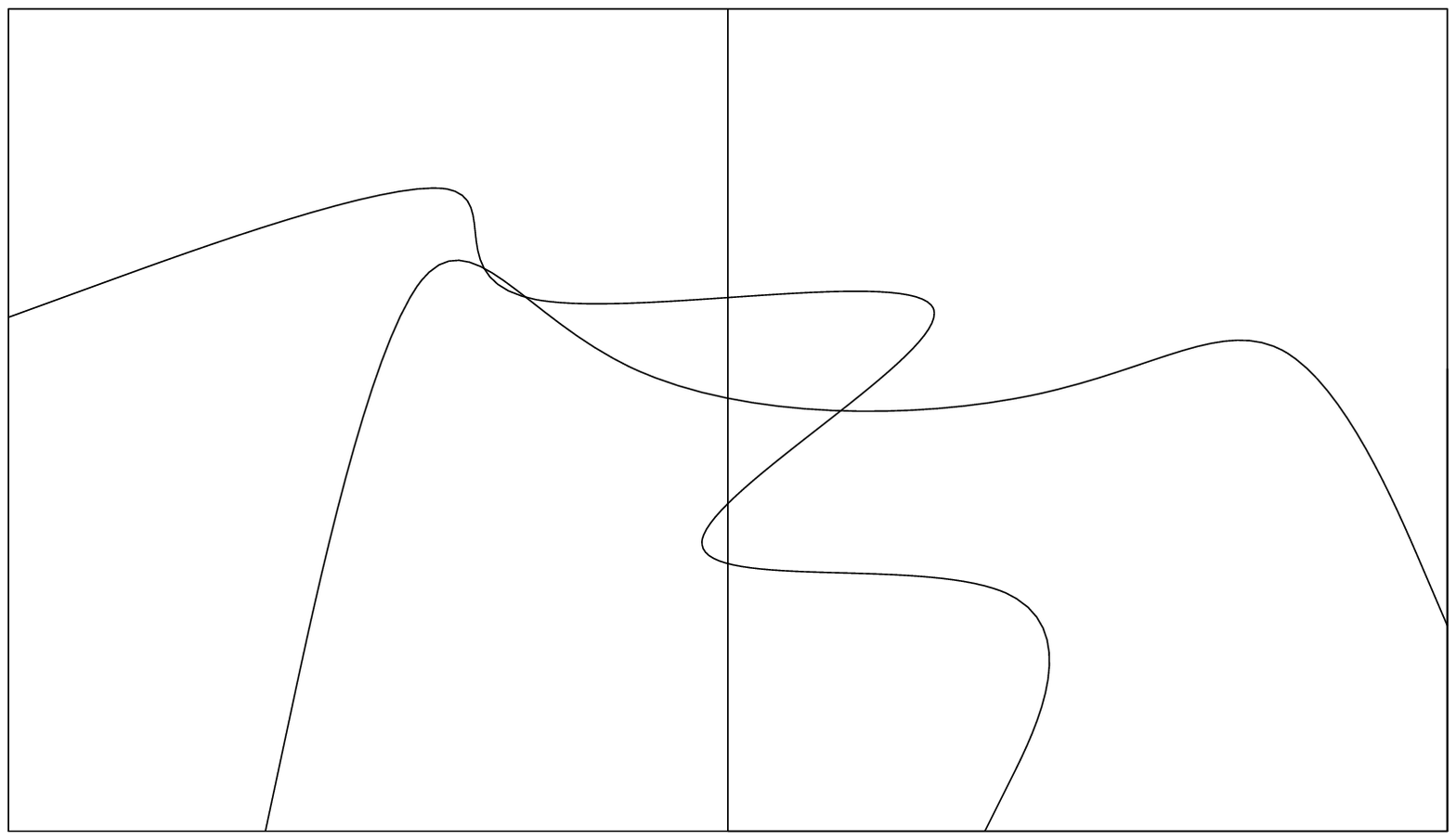}}
\caption {``Sketch'' of proof}
\end {figure}

Finally, note that $ A' \subset H(a,b) = \cup_g \{ \gamma = g \} $.
Putting all the pieces together, we get that
\begin  {eqnarray*}
P(A') &= &\sum_g P( A' \cap \{ \gamma = g \} )  \\
&\ge &\sum_g P ( A(g) \cap \{ \gamma = g \} )
= \sum_g P (A(g)) /2 = P( H(a,b))/ 2.
\end {eqnarray*}
Let now denote by $\tilde A$ the event that there exists a black crossing in $R(2a, b)$ joining the right boundary of the rectangle to the union of the left-boundary with the left-half of the bottom boundary. Symmetry shows that $P(\tilde A') = P (A')$. Furthermore, if the two increasing events $A'$ and $\tilde A'$ hold simultaneously, then $H(2a, b)$ holds too. Hence, we can conclude with the Harris inequality:
$$
 P( H(2a, b)) \ge P ( A' \cap \tilde A') \ge P (A') P(\tilde A') \ge P( H(a,b))^2 /4 .
$$
Note that we have not used at all the fact that the domain $R(2a, b)$ is a rectangle. The only feature that we used is that it is the union of a set
($R(a,b)$ in this case) and its symmetric image.

\medbreak
\noindent
{\bf Homework exercise:}
Try to prove similar results using (unions of) hexagons or rhombi instead of rectangles.

\section {First consequences}

Let us now define the hexagon centered at the origin with graph radius $n$, that consists of all points that can be joined to the origin by a nearest neighbor path on the lattice, with $n$ steps (or less).
In other words,
$$\Lambda_n = \{ ue^{ik \pi /3}  + v e^{i (k+1) \pi /3} \ : u \ge 0, v \ge 0 , u+ v \le n, k \in \{0,1,2, 3, 4, 5\}  \}.$$
We then define the concentric disjoint ``annuli'' $A_j= \Lambda_{2^{j+1}} \setminus \Lambda_{2^{j}}$.
We denote by $\partial \Lambda_n$ the set of points that are ``exactly'' at graph distance $n$ from the origin.

\medbreak

\begin{figure}
\centerline{\includegraphics*[height=1.9in]{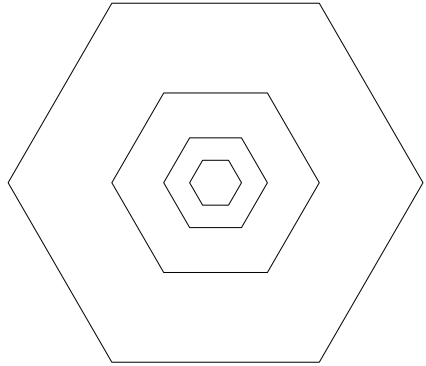}}
\caption {Concentric annuli}
\end {figure}

We say that there exists a closed circuit in $A_j$ if there exists a loop in $A_j$ that disconnects the origin from infinity, that consists only of closed (i.e. white) sites. The event is denoted by $C_j$.

\begin {lemma}
For any $j \ge 2$,
$P ( C_j) \ge (a_4)^6$.
\end {lemma}
We leave the proof of this lemma as a {homework exercise} (hint: One has to use the Harris inequality and note that if six well-chosen rotated rectangles are crossed by a closed path, then $C_j$ holds).

\medbreak

For any sets $D$ and $D'$ of sites, we denote by $D \leftrightarrow
D'$ the event that there exists a path of open sites joining a point
in $D$ to a point in $D'$.
\begin {corollary}
\label {co1.2}
For any $j,l \ge 2$,
$$P ( \partial \Lambda_{2^j} \leftrightarrow \partial
\Lambda_{2^{j+l}})   \le (1-a_4^6)^l.$$
\end {corollary}
This corollary follows immediately from the previous lemma because of
the independence of the realization of percolation in each of the
disjoint annuli $A_{j}$, $A_{j+1}$,$\ldots$, $ A_{j+l-1}$ and the fact that if
$\partial \Lambda_{2^j} \leftrightarrow \partial \Lambda_{2^{j+l}}$
then none of the independent events $C_{j}, \ldots , C_{j+l-1}$ holds.

\medbreak
 In particular, this shows that for any $n, m$ (such that $2^j \le n \le 2^{j+1}$, $2^l \le m < 2^{l+1}$),
$$
P ( \partial \Lambda_{n} \leftrightarrow \partial \Lambda_{nm} )
 \le P ( \partial \Lambda_{2^{j+1}}  \leftrightarrow \partial \Lambda_{2^{l+j}} )
 \le (1- a_4^6)^{l-1}
 \le  c m^{-\alpha}
$$
for two constants positive $c$ and $\alpha$ that are independent of $n$ and $m$.
In particular,
$$ P ( 0 \leftrightarrow \partial \Lambda_m ) \le c' m^{-\alpha}$$
for some constant $c'$.
\medbreak
\noindent
{\bf Homework exercise:} What sort of explicit lower bound does this argument give for $\alpha$?

 \medbreak
 \noindent
{\bf Easy homework exercise:} Prove that $P( 0  \leftrightarrow \partial \Lambda_m) \ge 1/(2(m+1))$
(hint: use the rhombus crossing probability).

\medbreak
Recall that we shall prove later in these lectures that in fact, $P ( 0  \leftrightarrow \partial \Lambda_m)
= m^{-5/48 + o(1)}$ as $m \to \infty$.

\medbreak
\begin {corollary}
\label {bigK}
There exist constants $K$ and $c$  such that for all $n  \ge 2$ and $m \ge 1$,
the probability that there exist $K$ disjoint open paths joining $\partial \Lambda_n$ to $\partial \Lambda_{mn}$
is smaller than $c \times m^{-3}$.
\end {corollary}

This follows readily from our estimates and the fact (called the Van den Berg-Kesten (BK) inequality, see \cite {G}) that
the
probability that there exist $K$ disjoint open paths joining $\partial \Lambda_n$ to $\partial \Lambda_{nm}$
is bounded from above by
$$ P ( \partial \Lambda_n \leftrightarrow \partial \Lambda_{nm})^K.$$
Then, we just have to take $K$ sufficiently large in such a way that $\alpha K > 3$.
We shall use this corollary a little later.

\medbreak
\noindent
{\bf A remark about tightness.} We want to understand the large-scale properties of our percolation model.
Let us for instance consider the rectangles $R(2n,n)$. At this point, we do not (yet) know that the crossing probabilities $P(H(2n, n))$ converge when $n \to \infty$. However, the Russo-Seymour-Welsh estimates give us a uniform lower bound on $P(H(2n, n))$. Similarly, they provide (using the lower bound for crossings by closed paths) a uniform upper bound. Hence, there exists $\epsilon > 0$ such that for any $n$, $\epsilon < P ( H(2n, n )) < 1 - \epsilon$. Consequently, one can find a subsequence $n_k \to \infty$ along which $P(H(2n_k,n_k))$ converges to a limit in $[\epsilon, 1 - \epsilon]$.

This can be iterated as follows: Since the sequence $P( H(4n_k, n_k))$ is also bounded away from $0$ and from $1$ for the same reasons, 
we can find a subsequence  $n_{\varphi (k)} \to \infty$ along which $P (H(4n_{\varphi(k)}, n_{\varphi(k)}))$ converges to some non-degenerate limit.
One can also then with the same argument find a sub-subsequence along which $P (H(3n, n ))$ also converges to some non-degenerate limit etc.
With the usual ``diagonal trick'', one can therefore find a subsequence $n_{\psi (k)}$ such that for all positive rational $\lambda$, the sequences $(P(H(\lambda n_{\psi (k)} , n_{\psi(k)} )), k \ge 1 )$ converge to some limit $f (\lambda) \in (0,1)$ (this limit could in principle depend on the choice of $\psi$).

Hence, Russo-Seymour-Welsh estimates show existence of subsequential limits with nice properties regarding crossing probabilities. We can also take other shapes than rectangles, and (as long as we have countably many of them) find subsequences for which all crossing probabilities converge.

As we shall see later, these estimates also show convergence of the law of renormalized exploration paths along subsequences.
This will turn out to be important in order to make the link to SLE in Lecture 3.

\chapter* {First exercise sheet}

\bigskip
\noindent
{\bf Color switching.}
Consider an $n \times m$ parallelogram and some integer $j \geq 2$. For any sequence of colors $\sigma=\sigma_1,\ldots,\sigma_j \in \{0,1\}$ ($0$ meaning ``white'' and $1$ ``black''), we consider the event
$C_{j,\sigma}$ that there exist $j$ disjoint left-to-right crossings, of respective colors given by $\sigma$ when ordered from bottom to top.

\smallskip
\noindent
1) Show that for any $\sigma$,
$$P(C_{j,\sigma}) = P(C_{j,\tilde{\sigma}}),$$
where $\tilde{\sigma}$ is defined by $\tilde{\sigma}_i = 1-\sigma_i$ for all $i=1,\ldots,j$.

\smallskip
\noindent
2) By conditioning on the lowest crossing, show that
$$P(C_{j,\sigma}) = P(C_{j,\sigma'}),$$
with $\sigma'_1 = \sigma_1$ and $\sigma'_i = 1-\sigma_i$ ($i \geq 2$).

\smallskip
\noindent
3) Generalize the previous result and deduce that the quantity $P(C_{j,\sigma})$ is independent of $\sigma$.

\smallskip
\noindent
4) Does this statement remain true if we replace the parallelogram by another simply connected subset on the lattice (like a rectangle $R(n,m)$) and consider  crossings from one fixed part of the boundary to another part?

\bigskip

\noindent
{\bf Two-arm exponent in the half-plane.}
We consider the ``half-hexagon'' $H_n = \{ z \in \Lambda_n \ : \ \Im (z)  \ge 0 \}$. Its boundary can be decomposed into two parts: A segment on the real axis and the ``angular semi-circle'' $h_n$.
We say that a point $x$ on the real line is $n$-good if there exist one open path originating from $x$ and one closed path originating from
$x+1$, that both stay in $x+H_n$ and join these points to $x+ h_n$ (these paths are called ``arms'').
Note that the probability $w_n$ that a point $x$ is $n$-good does not depend on $x$.

\smallskip
\noindent
1) We consider percolation in $H_{2n}$.

\noindent
a) Prove that with a probability that is bounded from below independently of $n$, there exists an open cluster $O$ and a closed cluster $C$, that both intersect the segment $[-n/2, n/2]$ and $h_{2n}$, such that $C$ is ``to the right'' of $O$.

\noindent
b) Prove that in the above case, the right-most point of the intersection of $O$ with the real line is $n$-good.

\noindent
c) Deduce that for some absolute constant $c$, $w_n \ge c/n$.

\smallskip
\noindent
2) We consider percolation in $H_n$.

\noindent
a) Show that the probability that there exists at least $k$ disjoint open paths joining $h_n$ to $[-n/2, n/2]$ in $H_n$ is bounded by $\lambda^k$ for some constant $\lambda$ that does not depend on $n$ (hint: use the BK inequality). Show then that the number $K$ of open clusters that join $h_n$ to $[-n/2, n/2]$ satisfies $P(K \ge k) \le \lambda^k$.

\noindent
b) Show that each $2n$-good point in $[-n/2, n/2]$ is the right-most point of the intersection of one of these $K$ clusters with the real line.

\noindent
c) Deduce from this that for some absolute constant $c'$,
$(n+1) w_{2n} \le  E(K) \le c'$.

\smallskip
\noindent
3) Conclude that for some positive absolute constants $c_1$ and $c_2$, $c_1 / n \le w_{n} \le c_2 / n$.

\bigskip

\noindent
{\bf Three-arm exponent in the half-plane.}
We say that a point $x$ is $n$-Good (mind the capital G) if it is the \emph{unique} lowest point in $x + H_n$ of an open cluster $C$ such that $C \not\subseteq x + H_n$.
Note that the probability $v_n$ that a point is $n$-Good does not depend on $x$.

\smallskip
\noindent
1) Show that this event corresponds to the existence of three arms originating from the neighborhood of $x$ in the half-hexagon $x + H_n$.

\smallskip
\noindent
2) Show that the expected number of clusters that join $h_{n/2}$ to $h_n$ is bounded. Compare this number of clusters with the number of $2n$-Good points in $H_{n/2}$ and deduce from this that for some constant $c_1$,
$v_n \le c_1 / n^2$.

\smallskip 
\noindent
3) Show that with probability bounded from below independently of $n$, there exists in $H_{n/2}$ an $n$-Good point (note that an argument is needed to show that with positive probability, there exists a cluster with a unique lowest point).
Deduce that for some positive absolute constant $c_2$, $v_n \ge c_2 / n^2$.

\bigskip

\noindent
{\bf Other half-plane estimates.}
Prove that the two-arm exponent in the half-plane, combined with color switching and the BK inequality implies that for some positive $\epsilon$,
$v_n \le c n^{-1- \epsilon}$.
What can one say about the four-arm exponent in the half-plane if one builds on the previous exercise?

\bigskip
\noindent
{\bf Five-arm exponent.}

\smallskip
\noindent
1) Consider the hexagon $\Lambda_{m}$ and, for even $m$, the annulus $A_m = \Lambda_{m} \setminus \Lambda_{m/2}$.

\noindent
a) Prove that the probability that there exist at least $k$ disjoint open paths joining the outer boundary of $A_m$ to its inner boundary
is bounded from above by $\lambda^k$, where $\lambda<1$ does not depend on $m$.

\noindent
b) Prove that the probability that there exist at least $k$ disjoint open clusters (if one considers percolation restricted to $\Lambda_m$) intersecting both the outer and the inner boundary of $A_m$ is bounded from above by $\lambda^k$. Deduce that the number $K$ of such clusters satisfies $E(K^2) < c$ for some absolute constant $c$.

\smallskip
\noindent
2) When $x$ is a site of the triangular lattice, we say that $U_m (x)$ is satisfied if
$x$ is closed and if out of its six neighbors originate five disjoint paths, three of them closed, two of them open, that are joining them to $x + \partial \Lambda_m$. We furthermore ask the colors of the five paths to be alternated: The two open paths are ``separated'' by closed paths in $x + \Lambda_m$. Note that the probability $u_m = P(U_m(x))$ does not depend on $x$.

\noindent
a) Suppose that $x \in \Lambda_{m/2}$ is such that $U_{2m} (x)$ holds. Show that it is on the boundary of two of the $K$ clusters.

\noindent
b) Conversely, suppose that two of the $K$ clusters are adjacent. Show that there are at most two points $x$ in $\Lambda_{m/2}$ on their joint boundary such that
 $U_{2m} (x)$ holds.

\noindent
c) Conclude that $(\# \Lambda_{m/2}) u_{2m} \le E( K^2) < c$.

\smallskip
\noindent
3) Prove using a Russo-Seymour-Welsh type argument that the probability that two of the $K$ clusters are adjacent, and only touch inside $\Lambda_{m/2}$ is bounded from below. Show that in this case, there exists at least one point in $\Lambda_{m/2}$ such that $U_{m/2} (x)$ holds.
Deduce that $(\# \Lambda_{m/2}) u_{m/2} \ge c'$ for some positive absolute constant $c'$.

\smallskip
\noindent
4) Conclude that there exist two constants $c_1$ and $c_2$ such that $c_1 m^{-2} \le  u_m \le c_2 m^{-2}$.

\bigskip
\noindent
{\bf Density of the infinite cluster.}
Our goal here is to show that the function $\theta(p) := P_p(0 \leftrightarrow \infty)$ can be viewed as the density of the infinite cluster 
${\mathcal C}$, in the sense that
$\# ({\mathcal C} \cap \Lambda_N) / \#\Lambda_N \to \theta (p)$ in $L^2$ (and in fact almost surely) as $N \to \infty$.

\smallskip
\noindent
1) Show that $E ( \# ({\mathcal C} \cap \Lambda_N )) = \theta (p) \times \# \Lambda_N $.

\smallskip
\noindent
Let us recall the following important property of percolation: for any $p < 1/2$, there exist constants $C,C'>0$ depending on $p$ such that for all $N$,
$$P_p(0 \leftrightarrow \partial \Lambda_N) \leq C' e^{-C N}.$$

\smallskip
\noindent
2) Prove that this implies the following for the supercritical regime: for any $p>1/2$, for all $N$,
$$P_p(0 \leftrightarrow \partial \Lambda_N\   | \ 0 \not \leftrightarrow \infty) \leq C' e^{-C N}.$$

\smallskip 
\noindent
3) Show that there exists a constant $C''>0$ (depending only on $p$) such that for all $N$,
$$\hbox {Var} ( \# ({\mathcal C} \cap \Lambda_N) )\leq C'' N^2.$$

\smallskip
\noindent
4) Show that the convergence takes place almost surely. Can one replace $(\Lambda_N)$ by any nested increasing family of domains?

\chapter{The Cardy-Smirnov Formula}

\section {Preliminaries}

We first collect some very basic facts:

\medbreak
\noindent
{\bf Analytic functions.}
We briefly make some elementary remarks  concerning analytic functions for those who never followed a complex analysis course or need to refresh their memories.

\begin {itemize}
\item
A $C^1$ complex-valued function $F$ defined in an open domain $O$ in the complex plane is called analytic if its differential at any point is just  multiplication by a complex number (recall that a differential is defined via a $2 \times 2$ real matrix i.e. four real numbers so that the fact that it is just multiplication by a complex number imposes two identities between these four real numbers). In other words, for each $z \in O$,  the limit of $(F(z+h) - F(z)) /h$ when $h$ tends to $0$ as a complex number exists. We call this limit $F'(z)$.

For each $\eta$ with $|\eta|=1$ (that should be thought of as a direction), we define
 $$\frac {\partial F}{\partial \eta} (z)= \lim_{h \to0, h \in \R} \frac {F(z+\eta h) - F(z)}{h}.$$
Clearly, if $F$ is analytic, then $(\eta)^{-1} \partial F / \partial \eta =  F'$ is independent of $\eta$.
In fact, the Cauchy-Riemann relation can be written as follows with this notation: $\partial F / \partial i = i \partial F / \partial 1$.
This identity between complex numbers corresponds to two equations between real numbers, and it suffices to check it in order to see that $F$ is analytic at $z$.

Recall that $\tau$ denotes the cubic root of unity $ \exp (2i \pi /3)$. Analyticity implies also that
$$ \partial F / \partial 1 = (\tau)^{-1} \partial F / \partial \tau = (\tau)^{-2} \partial F / \partial (\tau^2).$$
In particular, if we just take the real parts of these identities,
$$
 \partial ( \Re (F)) / \partial 1 = \partial ( \Re (F / \tau)) / \partial \tau = \partial ( \Re (F / \tau^2)) / \partial (\tau^2).
$$
These two equations also imply that the $C^1$ function $F$ is analytic ({homework exercise}).

\item
In fact, we will derive a discrete analog of the previous identities, but this is not the way in which we will prove that a certain function is analytic. We will rather use Morera's theorem involving contour integrals rather than derivatives. Let us briefly recall this result: Suppose that a function $F$ is continuous in $O$. Furthermore, assume that for any equilateral triangle in $O$ with one side parallel to the real axis, the contour integral of $F$ along this triangle vanishes. Then, $F$ is $C^1$ and analytic in $O$ (this is because one can approximate any contour integral of $F$ by sums of contour integrals along such equilateral triangles).

\end {itemize}

\medbreak
\noindent
{\bf Discrete exploration process.}
Suppose that a simply connected subgraph of the triangular lattice is given, as well as
two distinct points $A$ and $B$ on its ``boundary''. Then, the cells on $\partial D$
can be divided into two ``arcs'' $\BL$ and $\WH$ in
such a way that $A$, $\WH$, $B$ and $\BL$
are oriented clockwise ``around'' $D$.
Decide that all hexagons in $\BL$ are colored in black and
that all hexagons in $\WH$ are colored in white.
On the other hand, the cells in $D$ are chosen to
be black or white with probability $1/2$ independently of
each other. Consider now the (random) path $\gamma$
from $A$ to $B$ that separates the cluster of
black hexagons containing $\BL$ from the
cluster of white hexagons containing $\WH$.


\begin{figure}
\centerline{\includegraphics[height=2.3in]{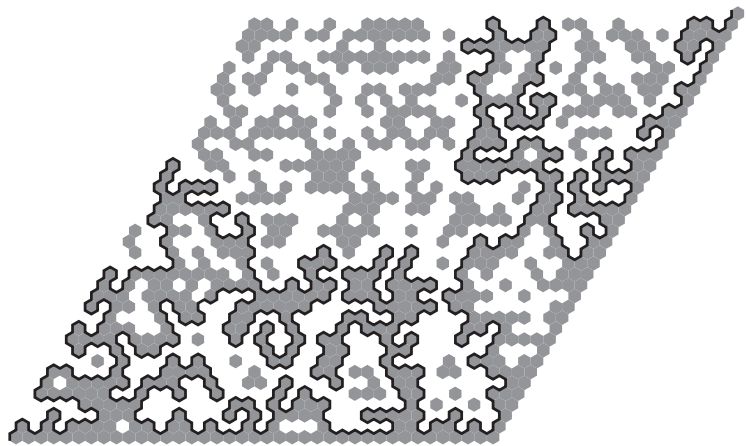}}
\caption{
An exploration process joining the two corners of a rhombus (picture by Julien Dub\'edat).}
\end {figure}

The curve $\gamma$ is therefore a simple (self-avoiding) path on the honeycomb lattice.
Another equivalent way ({homework exercise}) to define the interface $\gamma$ goes as follows:
It is a myopic self-avoiding walk. At each step $\gamma$ looks
at its three neighboring sites on the honeycomb lattice. It chooses at random
one of the sites that it has not visited yet (there are one or two such sites
since one site is anyway forbidden because it was the previous location of the
walk).

Hence, for each possible realization $g$ of $\gamma$, the event $\{\gamma =g \}$ means that all sites adjacent to the ``right'' of $g$ are black, and all sites to the left of $g$ are white. This event is therefore independent of the color of the cells that are not adjacent to $g$.
The fact that $\gamma$ can be explored as a myopic walk shows that the same holds for the beginning of $\gamma$. The event that the first $n$ steps of $\gamma$ follow the $n$ first edges of $g$ depends only on the color of the cells that are adjacent to the first $n$ steps of $g$.
This is the ``exploration property'' (that together with asymptotic conformal invariance should yield the conformal Markov property that characterizes SLE processes).

As we have already noticed in the proof of the Russo-Seymour-Welsh estimates,
exploration processes provide a way to discover ``left''-(resp. ``right'', ``top'', ``bottom'')-most
crossings.

\medbreak
\noindent
{\bf Conformal quadrilaterals.}
Suppose that $D$ is an open simply connected domain in the complex plane with $D \not= \C$. One way to state Riemann's mapping theorem is to say that there exist conformal (i.e. one-to-one analytic) maps $\Phi$ from $D$ onto the (inside of the) unit equilateral triangle $ABC$ in the complex plane ($A=0, B=1, C=e^{i\pi/3}$).
In fact there exist a three-parameter family of them and we can put further constraints on $\Phi$: if $a$, $b$, $c$ are three distinct points (or ``prime ends'') that are ordered anti-clockwise on $\partial D$, there exists exactly one conformal map $\Phi$ from $D$ onto the triangle, in such a way that $\Phi$ maps the three boundary points $a,b,c$ onto $A,B,C$ respectively.
If we now choose a fourth point $x$ on $\partial D$ (say, on the arc between $c$ and $a$), then the image $X = \Phi (x) \in [AC]$ of the fourth point is forced.

We call $(D,a,b,c,x)$ a conformal rectangle. It is possible to perform critical percolation of a fine-mesh (we call $\delta$ this mesh) lattice approximation of the domain $D$. Russo-Seymour-Welsh estimates yield readily that the probability of an open crossing from the (approximation) of the arc $(ab)$ to $(cx)$ is bounded away from $0$ and $1$, and one can wonder what the limit might be when $\delta \to 0$.

\section {Smirnov's theorem}

{\bf History of the problem and statement.}
Building on  the conformal field theory
ideas developed by Belavin, Polyakov and Zamolodchikov in  \cite {BPZ}, John Cardy \cite {Ca}
 gave an exact prediction
for the asymptotic value of crossing probabilities of conformal rectangles for critical percolation on
all canonical planar lattices, in terms of a
hypergeometric function. The meaning of conformal invariance for this model has been also investigated by Michael Aizenman \cite {Ai}.
Extensive numerical work (e.g., \cite {LPS})
 did comfort Cardy's
predictions. Lennart
Carleson  observed that Cardy's function is
closely related with the conformal maps
 onto equilateral triangles, and that Cardy's
prediction could be rephrased as follows:

\begin {conjecture}[Cardy's formula, Carleson's version]
If $D$ is conformally equivalent to the equilateral triangle $ABC$,
and if the four boundary points $a,b,c,x$
are respectively mapped onto $A, B, C, X \in [CA]$, then (in the
scaling limit when the mesh of the lattice goes to zero where one takes a lattice approximation of $D, a, b, c, x$),
the probability that there exists a
crossing in $D$ from the part $(ab)$ of $\partial D$ to $(cx)$
is equal to
$CX/CA$.
\end {conjecture}

We are now going to describe the proof of the following result:
\begin {theorem}[Smirnov \cite {Sm}]
Cardy's prediction is true in the case of critical site percolation
on the triangular lattice.
\end {theorem}

It is at present still a challenge to understand and prove why Cardy's formula holds for other lattices.

\medbreak

The striking feature in this statement is that crossing probabilities are conformally invariant (this is in a way more important than the actual precise formula) i.e. that any two conformally equivalent conformal quadrilaterals have the same asymptotic crossing probabilities.

\begin{figure}
\label {hit}
\centerline{\includegraphics*[height=3in, angle=-90]{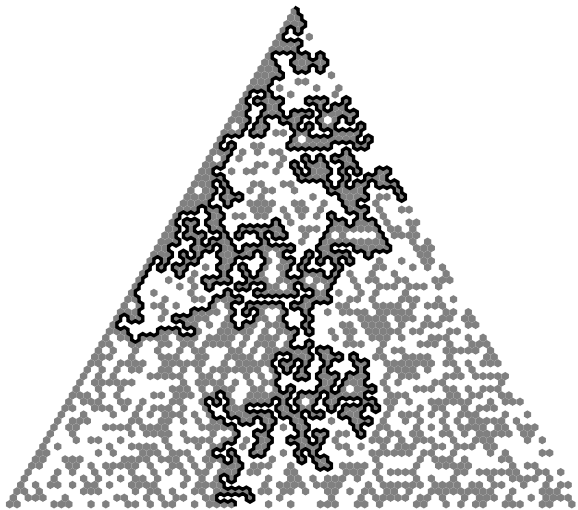}}
\caption{\label{f.exp2}
Exploration process in a triangle stopped at first hitting of the opposite side (pic. by Julien Dub\'edat).}
\end {figure}

\medbreak
\noindent
{\bf Homework exercise.}
Check that Smirnov's theorem implies that in the large-scale limit, the distribution of the hitting point of the exploration path with the bottom
boundary of the triangle (see Figure \ref{hit}) is uniform.

The rest of this lecture is devoted to the proof of Smirnov's theorem.
Suppose first for convenience that the domain $D$ is already the equilateral triangle and that $a=A$, $b=B$ and $c=C$.
For all $\delta = 1/n$, consider critical site percolation in
$ABC$ on the triangular grid with mesh-size $1/n$.
Write $A_{1}= A$, $A_\tau = A_2= B$ and $A_{\tau^2} = A_3= C$.
For each face $z$ of the triangular grid  (i.e. for each site of
the dual hexagonal lattice), let $E_1(z)$ denote the event
that there exists a simple open (i.e. black) path from $A_1A_\tau$
to $A_1A_{\tau^2}$ that separates $z$ from $A_{\tau}A_{\tau^2}$.
Similarly, define the event $E_\tau (z)$  corresponding
to the existence of a simple open path from $A_1 A_\tau$ to $A_\tau A_{\tau^2}$ separating $z$ from
$A_1A_{\tau^2}$,  and the event $E_{\tau^2}$ (corresponding to the existence of a path separating $z$ from $A_1A_\tau$).
Define finally for $j = 1, \tau , \tau^2$,
$$
 H_j^\delta (z) := P [ E_j^\delta (z) ].
$$
\medbreak
\noindent
{\bf Tightness.}
The arguments described in the previous lecture and in the
exercise session, ensure that the functions $H_j^\delta$
are uniformly  H\"older (actually, one
first has to define these functions in the continuous setting by interpolating them linearly in small triangles and keeping the original values of
$H_j^\delta$ at the center of the triangular faces).
Indeed, consider two points $z$ and $z'$ in the triangle with $|z-z'|< 1/100$, and suppose that for percolation with mesh-size $\delta$ (in the entire plane, not restricted to the domain $D_\delta$), there exists a closed circuit of diameter smaller than $1/4$ and an open circuit of diameter smaller than $1/4$ that surround both $z$ and $z'$. Then ({homework exercise}) this implies that $E_1^\delta (z)$ and $E_1^\delta (z')$ either both hold or both do not hold. In particular, it shows that
$|H_j^\delta (z) - H_j^\delta (z') |$ is bounded from above by the probability that there does not exist a closed circuit or an open circuit in some annuli (between radii $2 |z-z'|$ and $1/4$) around $z$, and this leads to a uniform H\"older bound
$$ |H_j^\delta (z) - H_j^\delta (z') | \le c |z - z'|^\beta$$
for suitable $c$ and $\beta$.

In particular, it shows that any for any sequence $\delta_n \to 0$,
the triplet of functions $(H_1^\delta, H_\tau^\delta, H_{\tau^2}^\delta)$
has a subsequential limit.
The  goal is now to identify the only possible
such subsequential limit.

Russo-Seymour-Welsh estimates also show that
when $d(z , A_{j \tau}A_{j \tau^2}) \to 0$, the function $H_j^\delta$
goes uniformly to zero, and that when $z \to A_j$, the function $H_j^\delta$
goes uniformly to one.
Hence, for any subsequential limit $(H_1, H_\tau, H_{\tau^2})$, one has
 $H_j (z) \to  0 $ when  $z$ approaches $A_{j\tau} A_{j \tau^2}$, and
 $H_j (z) \to 1$ when $z \to A_j$.

\medbreak
\noindent
{\bf Color switching.}
Now comes a key observation which is of a combinatorial nature:
Suppose that $z$
is the center of a triangular face (i.e. a site in the hexagonal lattice).
 Let $z_1, z_2,z_3$
denote the three (centers of the) neighboring faces
 (with the same orientation as the triangle $A_1A_2A_3$)
and $s_1,s_2, s_3$ the three corners of the face containing
$z$ chosen in such a way that $s_j$ is the corner ``opposite''
to $z_j$.
We focus on the event
$E_1^\delta (z_1) \setminus E_1^\delta (z)$.
This is the event that there exists three disjoint
paths $l_1$, $l_2$, $l_3$ such that
\begin {itemize}
\item
The two paths $l_2$ and $l_3$ are open and join
the two sites $s_2$ and $s_3$ to $A_1A_3$ and $A_1A_2$
respectively.
\item
The path $l_1$ is closed (i.e., it consists only of closed sites), and joins $s_1$ to $A_2 A_3$.
\end {itemize}
One way to check whether this event holds is to
start an exploration process from the corner $A_3$, say (leaving the
open sites on the side of $A_1$ and the closed sites on
the side of $A_2$).
If the event
$E_1^\delta (z_1) \setminus E_1^\delta (z)$ is true, then the exploration process
has to go through the face $z$, arriving into
$z$ through the
edge $z_3z$ dual to $s_1s_2$. In this way, one has
``discovered'' the simple paths $l_2$ and $l_1$ that
are ``closest'') to $A_3$.
Then, in the remaining (unexplored domain), there must
exist a simple open path from $s_3$
to $A_1 A_2$.
\begin{figure}
\centerline{\includegraphics[height=2in]{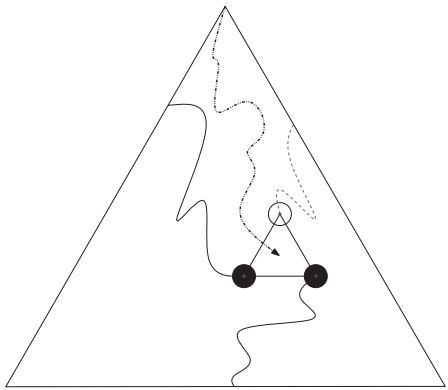}}
\caption{The three arms  (sketch)}
\end {figure}
But, the conditional probability of this event is the
same as that of the existence of a simple closed
path from $s_3$ to
$A_1 A_2$ (interchanging open and closed in the
unexplored domain does not change the probability
measure). Changing all colors once
again, shows finally  that $E_1^\delta  (z_1) \setminus E_1^\delta (z)$
has the same probability as the event that there exist
three disjoint
paths $l_1$, $l_2$, $l_3$ such that
\begin {itemize}
\item
The paths $l_1$ and $l_3$ are open and join
the two sites $s_1, s_3$ to $A_2A_3$ and $A_1A_2$
respectively.
\item
The path $l_2$ is closed, and joins $s_2$ to $A_1 A_3$.
\end {itemize}
This event is exactly
$E_\tau^\delta ( z_2) \setminus E_\tau^\delta (z)$.
Hence, we get that,
\begin {equation}
\label {der}
P [ E_1^\delta (z_1) \setminus E_1^\delta (z)]
=
P [ E_\tau^\delta (z_2) \setminus E_\tau^\delta (z)]
=
P [ E_{\tau^2}^\delta (z_3) \setminus E_{\tau^2}^\delta (z)]
.\end {equation}
These are reminiscent of the two equations relating derivatives of analytic function in three different directions that we recalled at the beginning of the lecture. So, it looks that we are on track to derive analyticity of a certain function.
Note however that
$$
H_1^\delta (z_1 ) - H_1^\delta (z) = P [ E_1^\delta (z_1) \setminus E_1^\delta (z)]
-P [ E_1^\delta (z) \setminus E_1^\delta (z_1)]
$$
so that the previous identities are not exactly identities between the discrete derivatives of the $H_j^\delta$'s in the three directions. Indeed, when one applies (\ref {der}) to the quantity $P [ E_1^\delta (z) \setminus E_1^\delta (z_1)]$ instead of 
$P [ E_1^\delta (z_1) \setminus E_1^\delta (z)]$, one ends up with an event involving $E_\tau^\delta(z')$ for some another neighbor $z'$ of $z_1$ (and not $z$).

\medbreak
\noindent
{\bf Contour integrals.}
(\ref {der})  can however be used to  show that
for any equilateral contour $\Gamma$ (inside the
equilateral triangle),
the contour integrals of $H_j^\delta$ for $j=1,\tau, \tau^2$
are very closely related. Loosely speaking,
$$
\oint_\Gamma dz H_1^\delta (z) =
\oint_\Gamma dz H_\tau^\delta (z) / \tau + O(\delta^\eps)
= \oint_\Gamma dz H_{\tau^2}^\delta (z) / \tau^2
+ O (\delta^\eps)
$$
when $\delta \to 0$ for some $\eps >0$.
More precisely:

\medbreak
\noindent
{\bf Homework exercise on telescoping sums.}
Define $$h_j^\delta (z, \eta ) = P [ E_j^\delta (z+\eta) \setminus E_j^\delta (z)]$$ where $\eta$ is chosen in such a way that
$z+\eta$ is one of the three neighbors of $z$.
Consider an approximation $\Gamma^\delta$ of $\Gamma$, such that $\Gamma^\delta$ is also an equilateral triangle consisting of
union of small triangular faces of the lattice with mesh $\delta$. Suppose for instance that $\Gamma$ is ``looking upwards'' (the horizontal side is the bottom side). Let ${\mathcal D}$ denote the set of centers of downwards looking small triangles in $\Gamma^\delta$.
Prove that:

\smallskip
\noindent
1) $h_j^\delta (z, \eta) = o ( \delta^\epsilon)$ for some $\epsilon >0$ (hint: use Russo-Seymour-Welsh if you do not wish to use directly the H\"older estimate).

\smallskip
\noindent
2) Show that for $\eta= i, i \tau, i\tau^2$,
$$ \sum_{z \in {\mathcal D}} ( H_1^\delta ( z + \eta) - H_1^\delta (z)) = \sum_{z \in {\mathcal D}} ( H_\tau^\delta (z + \eta \tau) - H_\tau^\delta (z) ) + o (\delta^{ \epsilon-1}).$$
(hint: Use the fact that $H_1^\delta (z+ \eta ) - H_1^\delta (z) = h_1^\delta (z, \eta) - h_1^\delta (z+ \eta , -\eta)$ and the identity (\ref {der})).

\smallskip
\noindent
3) If we call $I_\eta^\delta$ the previous identity, show that $i\delta ( I_i^\delta + \tau I_{i\tau}^\delta + \tau^2 I_{i\tau^2}^\delta)$ gives a meaning to the identity with integral contours. (hint: Note that almost all terms $H_j^\delta (z)$ disappear because they are multiplied by $1+ \tau + \tau^2$, and that the only remaining terms are those corresponding to triangles adjacent to the boundary of $\Gamma$).

\smallskip
\noindent
4) If the $H_j^{\delta_n}$'s converge to $H_j$'s, deduce from 3) that the contour integral of $H_1$ and of $H_\tau / \tau$ on $\Gamma$ are identical.

\medbreak

So, for any subsequential limit $(H_1,H_\tau, H_{\tau^2})$,
the contour integrals on equilateral triangles of
$H_1$ and $H_\tau/\tau$ -- and for analogous reasons of $H_{\tau^2}/\tau^2$ -- coincide.
 Hence, the contour integrals of the functions $H_1 + H_\tau + H_{\tau^2}$ and $G= H_1 + \tau H_\tau + \tau^2 H_{\tau^2}$ vanish too.
 By Morera's theorem, this ensures that these
functions are analytic. In particular $H_1 + H_2 + H_3$ is ananalytic real function -- it is therefore constant. Our estimates near to the corner $A$  show that
$H_1 + H_\tau+ H_{\tau^2} = 1$.
It follows that
$$\Re (G)= H_1 - 1/2 (H_2 + H_3) = H_1 - 1/2 (1 - H_1) = 3 (H_1 - 1/3) / 2$$
and similarly that $3H_\tau - 1 $ and $3 H_{\tau^2} -1$ are equal to  $2 \Re (G / \tau)$ and $2 \Re (G/ \tau^2)$.
We can indeed interpret the $(H_j - 1/3)$'s as the three conjugate functions of an analytic function and the relations (\ref {der}) as a discrete analog of the relation between their directional derivatives.

\medbreak
\noindent
{\bf Conclusion.}
So,  $H_{\tau^2}$ is harmonic because it is the real part of an analytic function.
We can note that $H_{\tau^2}$ extends continuously to the boundary of the triangle, that it takes the value $1$ at the top corner $C$, and that it vanishes on the bottom boundary $[AB]$.  Furthermore,  $G$
maps (for example) the segment $[AC]$ on a segment of the complex plane (because $H_{\tau}=0$ and
$H_1+ H_{\tau^2} =1$ on $[AC]$). A standard reflection argument (Schwarz reflection) shows that $G$ can be extended analytically to the neighborhood of this segment. In particular, it shows that the identities between directional derivatives still hold on these segments. As $H_{1} = 0$
on $BC$ and $H_\tau= 0$ on $AC$, it follows that the horizontal derivative of $H_{\tau^2}$ on $AC \cup BC$ vanishes.
Recall that $H_{\tau^2} (C) = 1$ and $H_{\tau^2} =0$ on $[AB]$.
The only harmonic function
in the equilateral triangle with these boundary conditions is
the normalized height
$$H_{\tau^2} (z) = 2 d(z, BC) / \sqrt {3}
$$
(it is unique for example because of the maximum principle, and this particular function satisfies all the conditions).
The particular case where $z \in [AC]$ is Cardy's formula when $D$ is the equilateral triangle.

\medbreak
\noindent
{\bf General domains.}
If $D$ is now any simply connected domain, and $a=a_1$, $b=a_{\tau}$, $c=a_{\tau^2}$ are boundary points,
the proof is almost identical. In its first part, the only difference is
that one replaces the straight boundaries $A_{j} A_{j\tau}$
by approximations of the boundary of $D$ on the triangular lattice that is between the
points $a_j a_{j \tau}$.
In exactly the same way, one obtains tightness and boundary estimates for the
discrete functions $H_j^\delta$. Also, the argument leading to the fact that
the contour integrals on equilateral triangles of
$H_1 + H_\tau + H_{\tau^2}$ and $G$ for any converging subsequential limit
vanish, remains unchanged. Hence, for any converging subsequential limit, one obtains a
triplet of functions $(H_1, H_\tau, H_{\tau^2})$ such that for $j=1, \tau, \tau^2$:
\begin {itemize}
\item
The function $H_1 + H_\tau + H_{\tau^2}$ is constant (because it is analytic and real) and equal to $1$ (because of RSW-estimates in the neighborhood of $a$).
\item
The function $G= H_1 + \tau H_\tau + \tau^2 H_{\tau^2}$ is analytic
\item
The function $H_j (x)$ tends to zero  when $x$ approaches the part of the
boundary between $a_{j\tau}$ and $a_{j\tau^2}$.
\item
The function $H_j (x)$ tends to one when $x \to a_j$.
\end {itemize}

The important feature is that this problem is conformally invariant:
If $\Phi$ denotes a conformal map from $D$ onto the equilateral triangle such
that $\Phi (a_j) = A_j$, and if $(H_1, H_\tau, H_{\tau^2})$ is such a triplet of
functions, then the triplet $(H_1 \circ \Phi^{-1}, H_\tau \circ \Phi^{-1},
H_{\tau^2} \circ \Phi^{-1})$ solves the same problem in the equilateral
triangle. In the latter case, we have just argued that the unique solution is given by
the functions $ 2 d( x, A_{j \tau} A_{j \tau^2}) / \sqrt {3}$.
This concludes the proof.

\medbreak

One should stress that this proves more than the
asymptotic behavior of the crossing probabilities:
It gives the asymptotic value of the probability of the events $E_j^\delta (z)$ also for $z$ inside the domain $D$ (and not only on its boundary).

In fact, the convergence to SLE(6) that we will derive in next lecture will imply the convergence of the probabilites of a much larger class of events.

\chapter{Convergence to SLE(6)}

\section {Our goal}

Let us fix a bounded simply connected domain $D$ with two distinct boundary points $x$ and $c$. We assume for convenience that $\partial D$ is in fact a continuous curve (and we use the term boundary points to designate prime ends i.e. the same point in the plane can for instance correspond to two boundary points of $D$, for instance in a slit domain).
As before, for each small $\delta$, we choose a lattice approximation $D_\delta, x_\delta, c_\delta$ of the triplet $D, x, c$ on the triangular lattice with mesh-size $\delta$. The boundary of $D_\delta$ is divided by $x_\delta$ and $c_\delta$ into two parts. Let us call $I_\delta$ one of them. We now perform critical percolation in $D_\delta$. We call $\gamma^\delta$ the outer boundary (i.e. viewed from $\partial D \setminus I_\delta$) of the union of all white cells that touch $I_\delta$. This is also the outer boundary (viewed from $I_\delta$) of the union of all black cells that are attached to
$\partial D_\delta \setminus I_\delta$.
This discrete interface $\gamma^\delta$  can be explored from $x_\delta$ to $c_\delta$
(or the other way around). We now also view these interfaces as continuous curves, and we choose to reparametrize $\gamma^\delta$ in such a way that they are defined for $u\in [0,1]$ with $\gamma^\delta (0) = x_\delta$ and $\gamma^\delta (1) = c_\delta$.


\begin{figure}[ht]
\centerline{\includegraphics[height=2.3in]{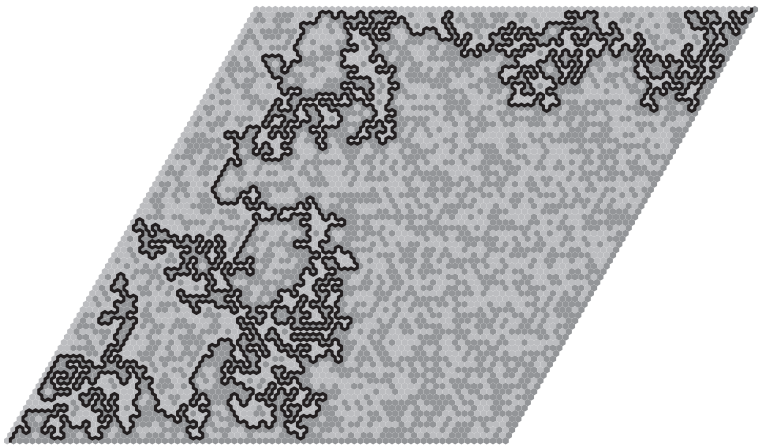}}
\caption{
An exploration process in a rhombus with small mesh-size (pic. by Julien Dub\'edat).}
\end {figure}

\begin {theorem}
\label {conv}
The law of $\gamma^\delta$ converges to the law of chordal SLE(6) from $x$ to $c$ in the domain $D$.
\end {theorem}

When one states a convergence in law, one has to be precise about the topology. Here, we use the distance between two continuous paths $\gamma$ and $\gamma'$ (with endpoints) defined by
$$ \inf_\varphi \sup_{u \in [0,1]} | \gamma (u) - \gamma' ( \varphi (u) )| $$
where the infimum is taken over all continuous increasing bijections $\varphi$ of $[0,1]$ into itself (see the second exercise sheet).

\medbreak

\section {Hand-waving argument}

Intuitively, the theorem is quite natural.
Recall from Greg Lawler's course that the SLE processes can be characterized in the following way: For each simply connected domain $D$ with two different boundary points,
we define a probability measure $P_{D,x,c}$ on parametrized (i.e. modulo linear time-reparametrization\footnote {This is due to the fact that in this set-up we want our curves to be invariant under the conformal maps from $D$ onto itself that preserve $0$ and $\infty$. In the case of the upper half-plane and the two boundary points, these are the maps $z \to \lambda z$ and they correspond to a linear time-change.}) non-self-traversing curves from $x$ to $c$ in $D$.
\begin {enumerate}
\item
Suppose first that this family is conformally invariant: For any conformal map $\Phi: D \to D'$, then the law of the image of $P_{D,x,c}$ under $\Phi$ is $P_{D',x', c'}$.
\item
Suppose furthermore that for any given $t$, the conditional law of $\gamma [t, \infty)$ given $\gamma [0,t]$ is $P_{D_t, \gamma_t, c}$, where
$D_t$ is the connected component of $D\setminus \gamma [0,t]$ that contains a neighborhood of $c$ (this is the domain Markov property).
\item
Suppose finally that for some $D,x,c$ that is symmetric with respect to the line $(xc)$, the law $P_{D,x,c}$ is symmetric with respect to this line.
\end {enumerate}
Then, the laws $P_{D,x,c}$ are that of SLE curves (parametrized by half-plane capacity) for some parameter $\kappa \ge 0$, respectively from $x$ to $c$ in the domain $D$.

Note that conditions 1. and 3. are easy to fulfill. Just take a symmetric scale-invariant measure on paths in $\HH, 0,\infty$, and map it conformally onto any other region $(D,a,b)$.

In our problem, we might argue that these three conditions should be satisfied for the laws of subsequential limits of exploration paths.
Indeed, symmetry and the domain Markov property hold in the discrete case, and we have just proved a conformal invariance statement
with the Cardy-Smirnov formula.
But putting this straight requires effort, for instance  to show that the discrete domain Markov property translates in the limit to the domain Markov property.
This is the strategy followed by Federico Camia and Chuck Newman in \cite {CN2}. We will follow here the somewhat different route outlined by Stas Smirnov in \cite {Sm2}.

In order to identify a candidate for the value of $\kappa$, one might note that in the scaling limit, we expect that if we choose two additional points $a$ and $b$ on the boundary of $D$, and if we map $D_t$ onto our favorite equilateral triangle $ABC$ by a conformal map $G_t$ in such a way that the three given marked points $a,b,c$ on $\partial D$  are mapped onto $A,B,C$, then the point $G_t (\gamma_t)$ is a martingale. Indeed, (if $x$ is on the arc between $a$ and $b$), it is the conditional probability given $\gamma [0,t]$ of the event that the curve $\gamma$ hits the arc $ca$ on $\partial D$
 before the arc $bc$. The computation (this is ``Cardy's formula computation for SLE'' in Greg Lawler's course) shows that SLE(6) is the only SLE with this property.

\section {Tightness}

Again, the first step will be to exploit the Russo-Seymour-Welsh estimates in order to get subsequential limits.
Basically, they ensure that the probability that the discrete curves $\gamma^\delta$ wiggle back and forth many times somewhere is very small, because this would imply multiple arm crossings of some small annuli with large modulus (and we know that such multiple crossings do not exist because we can show that their expected number vanishes uniformly in $\delta$ as the modulus goes to infinity -- this is where Corollary \ref {bigK} comes into play). Also, we have some a priori bounds (the one-arm estimate due to RSW)  that ensure that the curves $\gamma^\delta$ are far from plane-filling.
Hence, in our space of continuous functions, the laws of $\gamma^\delta$'s are tight, and existence of subsequential limits follow
(this result is  due to Aizenman and Burchard \cite {AB}).
We leave the details of the argument for the second exercise sheet.

A consequence of this result is that for any sequence of $\delta_n \to 0$, one can find an increasing sequence $n_k \to \infty$ such that
the subsequence $\gamma^{\delta_{n_k}}$ converges in law. Suppose that $\gamma$ follows the law of this subsequential limit. 
In fact, by Skorokhod's representation theorem, we can find a subsequence $n_k$ and a realization of all
 random curves $\gamma$ and $\gamma^{\delta_{n_k}}$ on the same probability space, such that
$\gamma^{\delta_{n_k}} \to \gamma$ almost surely.
In the following, we will therefore assume that $\delta_n \to 0$ is already chosen in such a way that the coupling with $\gamma^{\delta_n} \to \gamma$ a.s. exists (and we will work with this coupling).

It suffices to prove that the law of $\gamma$ is necessarily that of SLE(6) in order to complete the proof of Theorem \ref {conv}.

\section {Loewner chains}

Our goal is to prove that $\gamma$ is an SLE(6).
The first thing to check is that $\gamma$ can indeed be almost surely constructed via a Loewner chain.
Note that a priori, the curve $\gamma$ is fractal-like and has many double points, so that this is not a trivial statement.

Let us recall a few facts:

 \begin {itemize}
 \item
 Suppose that a continuous real-valued function $t \mapsto w_t$ is given.
 For each $z \not= 0$ in the closed upper half-plane $\overline \HH$, we define the function $t \mapsto g_t(z)$ that satisfies $g_0 (z) =z$ and
 \begin {equation}
 \label {loew}
  \partial_t g_t (z) = \frac {2}{g_t (z) - w_t}.
  \end {equation}
 This is well-defined at least for small $t$, and up to the first time
 $T(z)$ at which the distance between $g_t (z)$ and $w_t$ vanishes.
 If we define for each $t$ the set $K_t = \{ z \ : \ T(z) \le t \}$, then it turns out that $g_t$ is a conformal map from $\HH\setminus K_t$ onto $\HH$ such that
 $ g_t (z) = z + o(1)$ as $z \to \infty$.
 This increasing family $(K_t)$ is called a Loewner chain in the upper half-plane, generated by the driving function $w$. Note that with this construction, $K_t$ is an increasing family of compact sets parametrized by half-plane capacity (as defined in Greg Lawler's course)\footnote {The half-plane capacity of $K$ is the unique non-negative integer $a$ such that there exists a conformal map $g$ from $\HH \setminus K$ to $\HH$ with an expansion $g(z)= z + 2a/z + o (1/z)$ as $z \to \infty$. It is a nice way to measure the ``size'' of $K$ seen from infinity, especially when one iterates such conformal maps $g$.}.
 \item
 It is not difficult to prove that any (properly parametrized) {\em simple} continuous curve $\gamma$ in $\HH$ with $\gamma(0)=0$ is a time-changed Loewner chain. One just first time-changes $\gamma$ in such a way that the half-plane capacity of $\gamma [0,t]$ is $t$, then defines $g_t$ as the normalized conformal map from $\HH\setminus K_t$ onto $\HH$, and finally checks that $w_t= g_t (\gamma(t))$ indeed defines the set $K_t = \gamma [0,t]$ via the differential equation (\ref {loew}).  But, one should keep in mind that there exists many Loewner chains that are not simple curves.
\item
If we now fix a domain $D$ with two boundary points $x$ and $c$ as before, we choose a fixed conformal map $\Phi$ from $\overline D$ onto $\overline \HH$ such that $\Phi (x) = 0$ and $\Phi (c) = \infty$. We say that $(K_t)$ is a Loewner chain from $x$ to $c$ in $D$ if its image under $\Phi$ is a Loewner chain in $\HH$ as defined above (note that since $\Phi$ is not unique, the time-parametrization is defined up to linear time-change only).
 \end {itemize}
A general characterization for a family $(K_u)$ in $\HH$ to be a time-changed Loewner chain requires two things: First, the half-plane capacity has to be continuous and increasing; this makes it possible to reparametrize the Loewner chain by its half-plane capacity (that we call $t$). Then, one requires $(K(\cdot))$ to grow  ``locally seen from infinity'' i.e. that the diameter of the image of $K_{t+ s}\setminus K_t$ under $g_t$ is small for small $s$ (i.e. $K_t$ grows only at one point). More precisely,
 one has to check that for any $t_0 \ge 0$, for all $\epsilon$, there exists $s_0$ such that for any $t \le t_0$,
 the diameter of $g_t ( K_{t+s_0} \setminus K_{t} ) $ is smaller than $\epsilon$.
 Indeed, if these two criteria are satisfied, one can check that if one defines for each $t$, $w_t = \cap_{s >0} (\overline { g_t (K_{t+s} \setminus K_t)})$, the family $K_t$ (if parametrized by its half-plane capacity) is generated by the driving function $w_t$ (see \cite {LSW1} for details).

In fact (see \cite {LSW1} again), the local ``growth condition'' can be replaced by the following condition: for any $t_0 \ge 0$, for all $s_0$, there exists $\epsilon'$  such that for any $t \le t_0$, the set $K_{t+s} \setminus K_t$ is separated from infinity in $\HH \setminus K_t$ by a set of diameter not greater than $\epsilon'$.

To check that a family $(K_u)$ is a Loewner chain to $c$ in the domain $D$, one just has to check that $\Phi(K_u)$ satisfies the previous conditions. In the following, we will fix a given conformal map $\Phi$ from $(D,x,c)$ to $(\HH,0, \infty)$ and we will refer to the capacity of $K_u$ (in $D$, seen from $c$) as the half-plane capacity of $\Phi (K_u)$.

\medbreak

In our setup, the random curve $(\gamma_u, u \in [0,1])$ is given (and a priori, it is not parametrized by its capacity).  For each $u$, we define $D_u$ the connected component of $D \setminus \gamma [0,u]$ that has $c$ on its boundary, and $K_u = D \setminus D_u$. We can furthermore assume that $\gamma$ is non-constant on any interval.
Note that because $\Phi(\gamma)$ is a continuous curve, it follows immediately that the half-plane capacity of $\Phi(K_u)$ is continuous with respect to $u$.

If we could prove that the half-plane capacity of $\Phi(K_u)$ is (strictly) increasing, then  it would follow that the (time-changed) family $(K_t)$ is a Loewner chain, because the ``local growth'' condition follows immediately from the uniform continuity of $\Phi (\gamma)$ on any compact interval $[0,t_0+s]$.  Indeed, the set $\Phi (K_{t+s} \setminus K_t)$ is disconnected from infinity in $\HH \setminus \Phi (K_t)$ by $\Phi ( \gamma [t, t+s])$.

\section {Capacity increases}

In this section, we are going to argue that almost surely, the capacity of $(K_u)$ (seen from $c$ in $D$) is a strictly increasing function. This amounts to prove that $(K_u)$ itself is a strictly increasing family.

For all $z \in D$, let us define the time $\sigma_z$ (respectively $\sigma_z^\delta$) to be the first time $u$ at which $\gamma [0,u]$
(resp. $\gamma^\delta [0,u]$) does disconnect $z$ from $c$ in $D$ (resp. $z_\delta$ from $c_\delta$ in $D^\delta$). Recall that (because of RSW estimates), $\gamma$ does almost surely not hit the point $z$.

\begin {lemma}
\label {hitt}
For each fixed $z \in D$ with rational coordinates, $\sigma_z^{\delta_n} \to \sigma_z$ almost surely as $n \to \infty$.
\end {lemma}

Note that the convergence of $\gamma^{\delta_n}$ to $\gamma$ implies immediately that almost surely,
$$\sigma_z \ge \limsup_{_n \to \infty} \sigma_z^{\delta_n}.$$
What one has to prove here is that when $\gamma^{\delta_n}$ comes close to disconnect $z$ (this has to be the case at time $\sigma_z$ since $\gamma^{\delta_n}$ is close to $\gamma$), then it does disconnect $z$ quickly and with high probability.
This can be done using the Russo-Seymour-Welsh estimates in the
intersection of annuli with the domain that the exploration path has
not yet visited at this close-to-disconnection time
(see Figure~2).
We leave the proof of this lemma to a {homework exercise}.

\begin{figure}[ht]
\label{hint}
\centerline{\includegraphics*[height=1.3in]{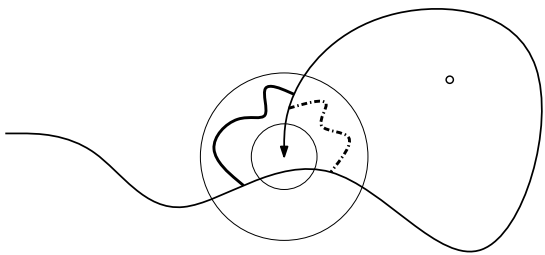}}
\caption{Getting close to disconnect implies likely to disconnect.}
\end {figure}

\medbreak
\noindent
{\bf Remark for those who have looked at exercise sheet \#1:}
Note that this lemma is closely related to the a priori estimates on 5 arms in the plane and 2 or 3 arms in the half-plane. Indeed, if the lemma would not be true, then this means that with positive probability, there exists a point $\gamma (\sigma_z)$ in the neighborhood of which
the approximating percolation models have $6$ arms (if this point is in $D$) or 3 arms (if the point is on $\partial D$), and these estimates indicate that this is not the case (one has first to note that the $6$-arm exponent is strictly larger than $2$ and to generalize the estimates to annuli). But proving the lemma this way would be longer than the suggested proof above.

\begin {lemma}
For all $u<u'$, there almost surely exists $v \in (u,u')$ such that $\gamma_v \notin \gamma[0,u] \cap \partial D$.
\end {lemma}

\noindent {\bf Proof.} Define a countable dense set of points $(x_j)$ on the union of the boundaries of the connected components of the complement of $\gamma [0,u]$. Because of the Russo-Seymour-Welsh estimates, almost surely, none of the points $x_j$ are visited by $\gamma$, and this implies the lemma.

\medbreak

These two  lemmas imply that
the map $u \mapsto K_u$ is strictly increasing. Indeed, almost surely, for any rational $0<u<u'$, there exists $v \in (u,u')$ with $\gamma (v) \notin \gamma [0,u]$. The point $\gamma_v$ is therefore in (the inside) of one of the connected components of $D \setminus \gamma [0,u]$. Suppose that it is not the connected component with $c$ on its boundary. $\gamma_v$  is the same connected component as some point $z$ with rational coordinates.  But at the time (or very shortly after) at which this connected component has been disconnected, the curves $\gamma^{\delta_n}$'s did disconnect $z$ as well (with high probability), and are therefore prevented from entering this component (the curves $\gamma^{\delta_n}$ do not have double points and have to end at $c$); this leads to a contradiction.
Hence, $\gamma_v$ is in the connected component with $c$ on its boundary, and we get that for all rational times $u<u'$, $K_u \not= K_{u'}$
almost surely.

\medbreak

Hence,  we see that if we time-change $\gamma$ in such a way that the half-plane capacity satisfies $\hbox {hcap}(\Phi (K_t))= t$, then $\Phi (K_t)$ is a Loewner chain in the upper half-plane generated by a random continuous function $w_t$.

Note that since $\gamma^{\delta_n} \to \gamma$, if we parametrize also each of the $\gamma^{\delta_n}$'s by its capacity, we get that
$ \sup_t | \gamma^{\delta_n} (t) - \gamma (t) | \to 0$ almost surely. This follows from the fact that the curve $\gamma$ is a.s. uniformly continuous.
We will from on now work with these curves parametrized by their half-plane capacity.

\section {Side-remark concerning the proof of Cardy-Smirnov formula}

In the proof of the Cardy-Smirnov formula (and in fact in its statement) and in the statement of Theorem \ref {conv}, we have not been precise on how the discrete domains $D_\delta$ were allowed to  approximate the continuous domain $D$ in order for the statement to hold. Let us temporarily come back to the notations used in the proof of Cardy's formula.  $D_\delta$ is a simply connected domain on the triangular grid of mesh-size $\delta$. Its boundary can be decomposed into three paths $\partial_1$, $\partial_\tau$ and $\partial_{\tau^2}$. Existence of subsequential limits for the function $H_j^\delta$ builds on Russo-Seymour-Welsh estimates and require just that any $z$ in the inside of $D$ is in the inside of $D_\delta$ when $\delta$ is sufficiently small. Proving the statement about contour integrals in the inside of the triangle requires no other condition on the $D_\delta$'s. Hence, analyticity of the functions $H_1 + H_\tau + H_{\tau^2}$ and $H_1 + \tau H_\tau+ \tau^2 H_{\tau^2}$ still holds.
When handling the boundary conditions for these functions, we use  the Russo-Seymour-Welsh estimates again. We leave it as a homework exercise to check that the following additional conditions are sufficient to conclude:
\begin {itemize}
\item
$(a_\delta, b_\delta, c_\delta) \to (a,b,c)$.
\item
For any $z \in \partial_j$, the distance of $z$ to $\partial_j^\delta$ converges to $0$ when $\delta \to 0$ uniformly with respect to $z$.
\end {itemize}

\section {Identifying continuous martingales}

Let us now come back to our setup (the domain $D$, boundary points $c$ and $x$, the almost sure limit $\gamma$).
Let us fix any two points $a$ and $b$ on the boundary of $D$ in such a way that $a,x,b,c$ are ordered anticlockwise on $\partial D$.
Let us define the event ${\mathcal A}^{\delta_n}$ that $\gamma^{\delta_n}$ hits the arc $a_{\delta_n} c_{\delta_n}$ before the arc $c_{\delta_n} b_{\delta_n}$. Then, since the connected component of $D \setminus \gamma^{\delta_n}[0,t]$ with $c_{\delta_n}$ on its boundary converges (and the four boundary points $(a_{\delta_n}, \gamma^{\delta_n}_t, b_{\delta_n}, c_{\delta_n})$) to that of $D \setminus  \gamma [0,t]$ (with the four points $(a, \gamma_t, b, c)$) in the above sense, we get that almost surely,
if $\gamma [0,t] \cap (ca \cup bc) = \emptyset$, then
$$P ( {\mathcal A}^{\delta_n} | \gamma^{\delta_n} ([0,t])) \to X_t$$
where $X_t$ is the image of $\gamma_t$ under the conformal map that maps
$D_t, a, b, c$ onto the equilateral triangle $ABC$.

Let ${\mathcal A}$ denote the event that $\gamma$ hits the arc $ac$ before the arc $bc$. Then, by Lemma \ref {hitt},
 the probability of the symmetric difference between ${\mathcal A}$ and ${\mathcal A}^{\delta_n}$ goes to $0$ as $n \to \infty$.
Indeed, the hitting time of $ac$ is equal to the disconnection time of some point $z$ with rational coordinates close to the corner $a$.

\medbreak

Let $T$ denote the hitting time of $bc \cup ca$ by $\gamma$. We define $X_t = 1_{\mathcal A}$ if $t \ge T$. Note that $X_t$ is  a deterministic
function of $\gamma [0,t]$.

\begin {lemma}
For all $t$,
$X_t= P ( {\mathcal A} | \gamma [0,t] )$.
\end {lemma}

\noindent
{\bf Proof.}
Let us consider any continuous bounded function $f$ on the space of curves.
Using three times dominated convergence, and our previously collected results,  we obtain that
\begin {eqnarray*}
  E ( 1_{\mathcal A} f(\gamma [0,t]))& =& \lim_{n \to \infty} E ( 1_{{\mathcal A}^{\delta_n}}  f ( \gamma^{\delta_n} [0,t] )) \\
  &=& \lim_{n \to \infty} E ( P ( {\mathcal A}^{\delta_n} | \gamma^{\delta_n} [0,t] ) f ( \gamma^{\delta_n} [0,t] )  )\\
&=& \lim_{n \to \infty} E ( X_t f( \gamma^{\delta_n} [0,t]) ) \\
&=& E ( X_t f ( \gamma [0,t] ))
\end {eqnarray*}
and this proves the lemma.

\medbreak

Note that its definition shows that the function $X_t$ is continuous with respect to $t$. Hence, it is a continuous martingale i.e. a time-changed one-dimensional Brownian motion.

\section {Recognizing SLE(6)}

We are now ready to conclude. Let us consider the Loewner chain $(\Phi (K_t))$ defined in the upper half-plane by $\Phi (\gamma (t))$.
Define the conformal maps $g_t$ associated to it, and its continuous and random driving function $w_t$.
Define also the points $a'= \Phi (a)$ and $b'= \Phi (b)$.

We shall use the Schwarz-Christoffel conformal map $\Psi$ from the upper half-plane onto the triangle such that $\Psi (0) = A$, $\Psi (1) = B$ and
$\Psi (\infty) = C$. Recall (e.g. \cite {Ah}) that the restriction of $\Psi$ to $[0,1]$ is given  by
$$
\Psi (x) = \hbox {cst} \times  \int_0^x \frac {dy} {y^{2/3} (1-y)^{2/3} }
$$
and that $\Psi$ satisfies the equation
$$
3\Psi''(x) + 2\left( \frac {1}{x} +\frac {1}{x-1} \right) \Psi'(x) = 0 .
$$

Note that  for all $t \le T$
$$ X_t = \Psi \left( \frac { w_t - g_t (a') } {g_t(b') - g_t (a') } \right)$$
(because $z \mapsto (z - g_t(a')) / (g_t(b') - g_t (a'))$ maps $\HH$ onto itself and the three boundary points $g_t(a')$, $g_t(b')$ and $\infty$ onto $0$, $1$ and $\infty$).
Hence,
\begin {equation}
\label {wt}
w_t = g_t (a') + (g_t (b') - g_t (a')) \Psi^{-1} (X_t)
.\end {equation}
Note that because
$\partial_t g_t (x) = 2/ (g_t(x) - w_t)$ for any fixed $x$, the functions $g_t(a')$ and $g_t (b')$  are  determined from $X_t$ as long as $t \le T$ (they are solution to a simple ODE, just replace $w_t$ by its value (\ref {wt}) in the differential equation for $g_t(a')$ and $g_t(b')$), and the function $w$ therefore too.  In other words, the random processes $ w_t, g_t (a')$ and $g_t (b')$ are measurable with respect to the filtration of $X$ (as long as $t \le T$).

Since $t \mapsto g_t (a')$, $t \mapsto g_t (b')$ are $C^1$ functions, and $\Psi^{-1}$ is $C^2$, it follows that
$w_t$ is a semi-martingale (see e.g. \cite {RY}) i.e.
one can write
it on $[0,T]$ as
$ w_t = M_t + V_t $
where $M$ is a local martingale and $V$ is a process with bounded variation measurable with respect to the same
filtration, and one can use It\^o's formula to study functions of this local martingale.

Recall that the choice of points $a$ and $b$ was arbitrary. The previous argument shows that for all $a' < w_0 < b'$ on the real line, the processes
$$
\Psi \left( \frac { w_t - g_t (a') } {g_t(b') - g_t (a') } \right)
$$
are local martingales i.e. that there is no drift term in the semi-martingale decomposition of these processes.
Applying the rules of stochastic calculus (homework exercise, this is almost identical to the computation for Cardy's formula in terms of SLE in Greg Lawler's course) together with the ODE for $\Psi$, we readily get that
if $Z_t = (w_t - g_t (a')) /(g_t(b') - g_t (a'))$,
$$
\frac {dV_t}{g_t (b') - g_t (a')} \Psi' (Z_t)
+
\frac {\Psi'' (Z_t)}{2(g_t (b') - g_t(a'))^2} ( d \langle M \rangle_t - 6 dt ) =0,
$$
(where $\langle M \rangle_t$ denotes the quadratic variation of $M$, see e.g. \cite {RY}).
Since this has to be true for all $a' < w_0 < b'$ up to $T=T(a',b')$ and $T(a',b') \to \infty$ as $a' \to -\infty$ and $b' \to +\infty$, we conclude that
$V_t =0$ and $\langle M \rangle_t = 6t$ for all $t \ge 0$.
 In other words, $w(t/6)$ is a standard one-dimensional Brownian motion, and $\gamma$ is an SLE(6) process. This concludes the proof of Theorem \ref {conv}.

\medbreak

A by-product of this theorem
is the locality property of SLE(6) (i.e. Theorem 4.34 in Greg Lawler's course, where it is derived directly in the SLE framework without any reference to percolation), because the corresponding property is clear in the discrete setting:

\begin {prop}
\label{locality}
Consider a chordal SLE(6) process $\gamma$ from $a$ to $b$ in the domain $D$, and a chordal SLE(6) process $\gamma'$ from $a$ to $b'$ in $D$.
Define $T$ (respectively $T'$) the first time at which $\gamma$ (resp. $\gamma'$) disconnects $b$ from $b'$ in $D$. 
Then, the two paths $\gamma [0,T]$ and $\gamma' [0,T']$ (defined modulo time-reparametrization) have the same law.
\end {prop}

Two other results (and these ones are not so easy to derive directly in the SLE framework, see Rohde-Schramm \cite {RS} for the first one -- and I do not know if there is another proof of the second one) that Theorem \ref {conv} implies immediately are the following:
\begin {itemize}
\item
SLE(6) is almost surely generated by a random continuous curve $\gamma$.
\item
SLE(6) is reversible: The law of the time-reversal of an SLE(6) from $c$ to $x$ in $D$ is that on SLE(6) from $c$ to $x$ in $D$ (after reparametrization).
\end {itemize}

\section {Take-home message}

At the end of the first three lectures, we can now look back at what we have proved so far: We started with general features of critical percolation on the triangular lattice, and we have identified a quantity that converges in the scaling limit (when the domain $D$ is fixed and the mesh of the lattice goes to zero) to a conformally invariant limit. Then, we have proved that in fact (for rather weak assumptions on the way the discrete lattice has to approximate the domain), the entire discrete exploration process converge in the scaling limit to a conformally invariant random curve called SLE(6) that we can perform computations with. Furthermore, and this was and will be essential, we did see that discrete hitting times by the discrete exploration process do converge to the hitting times by the SLE(6) curve. Hence, this basically shows that the probability of every event that can be described in terms of hitting times by an exploration process (or by properties of several exploration processes) converge in the scaling limit to the
probability of the corresponding event for SLE(6) (which is a conformally invariant quantity). This allows for instance to compute explicitly the limit of the probabilities of certain simple connectivity events, see e.g. \cite {S2, Dub}.

\chapter*{Second exercise sheet}

We still study critical site percolation on the triangular lattice with mesh $\delta>0$. We consider a discrete domain $\Omega^\delta$, approximating a simply connected compact domain $\Omega$. For each $\delta$, the law of the discrete  exploration interface $\gamma^\delta$ (corresponding to given boundary conditions) is a distribution on continuous curves. We would like to show that this family of probability measures is regular enough to possess subsequential limits (that are also laws of  random curves) when $\delta \to 0$. The goal of this exercise sheet is to present a proof of this result. In the original proof of Aizenman and Burchard \cite {AB}, stronger results are derived, namely more precise estimates of the H\"older regularity of the curve under the limiting measures as well as dimension bounds, but this is not needed for our purposes.

\medbreak

\medskip

Recall that a family $(\gamma^\delta, \delta >  0)$ of random variables with values in a metric space ${\mathcal S}$ is called {\em tight} if for all positive $\eps$, one can find a compact set $K$ in ${\mathcal S}$,   such that $\inf_\delta P(\gamma^\delta \in K) > 1-\eps$. We will need the following part of Prokhorov's Theorem:

\medbreak
\noindent{\bf Prokhorov's Lemma:}  If a sequence $(\gamma^{\delta(n)}, n \ge 0)$ is tight, then there exists $n(k) \to \infty$ such that $\gamma^{\delta(n(k))}$ converges in distribution as $k \to \infty$.

\medskip

Our {space of curves} ${\mathcal S}$ can be defined as the set of equivalence classes of continuous functions from $[0,1]$ to the plane, where two functions $f_1$ and $f_2$ represent the same curve if and only if there exists a continuous monotone bijection $\phi$ from $[0,1]$ onto itself such that $f_1  = f_2 \circ \phi$. We endow this space with the metric
$$d(f_1,f_2) = \inf_{\phi} \sup_{[0,1]} \|f_1(u) -f_2 \circ \phi (u)\|,$$
where the infimum is over the set of continuous monotone bijections $\phi$ (note that this distance does not depend on the choice of the functions $f_1$ and $f_2$ in their equivalence classes).
Note finally that the discrete interface $\gamma^\delta$ in $\Omega^{\delta}$ is also an element of $\mathcal{S}$.

\medbreak
\noindent
{\bf Question 1.} Check that $d$ is indeed a metric.

\medbreak
\noindent
For each $n$ and each curve $\gamma$, we define by induction $T_0^n = 0$, and for each $j \ge 1$, 
$$  T_j^n (\gamma) = \inf \{ t > T_{j-1}^n (\gamma) \ : \ | \gamma (t) - \gamma  (T_{j-1}^n (\gamma)) | > 16 \times 2^{-n} \}$$
and $M(n)=M (n, \gamma)$ the maximum number of such steps before the path reaches its end-point. 

\medbreak
\noindent
{\bf Question 2.} Suppose that a sequence $u(n)$ is fixed.
Prove that for each $n_0$, the set of curves $\gamma$ such that for all $n \ge n_0$, $M(n, \gamma) \le u(n)$ is compact in ${\mathcal S}$. 
(hint: Consider a sequence of such curves, and find a subsequence along which all $M(n,\gamma_k)$'s converge, as well as all 
$\gamma_k (T_j^n ( \gamma_k))$'s).

\medbreak
\noindent

We now would like to derive some properties of our random functions $\gamma^\delta$.
We start with the following very useful consequence of Corollary \ref {co1.2}:

\medbreak

\noindent 
{\bf Question 3.}
Show that there exist two positive constants $\alpha$ and $C$, such that for all $k$, $\delta$, $r$ and $x$, 
the probability that there exist $k$ disjoint open crossings of the annulus $A(x,r)= \{z \ : \ r < |z - x | < 4r \}$  is bounded by $C2^{-k\alpha}$.

\medbreak
\noindent
{\bf Question 4.}
Deduce from the previous question that the probability that $\gamma^\delta$ crosses a given annulus $A(x,r)$ more than $2k$ times is bounded by $C 2^{-k \alpha}$.

\medbreak
\noindent
{\bf Question 5.}
Prove that for some well-chosen fixed points $x_1, \ldots, x_{N}$ for $N = N(n) \le  C' 4^n$,  for each $j \le M(n, \gamma)$, at least one of the $N(n)$ annuli $A(x_i, 2^{-n})$ 
is crossed by $\gamma^\delta$ between times $T_{j-1}$ and $T_j$. 

\medbreak
\noindent 
{\bf Question 6.}
Prove that for a well-chosen large $K$, the probability that there exists an $i \le N$ such that $A(x_i, 2^{-n})$ is crossed at least $K n$ times is bounded by $C'' 2^{-n}$. 

\medbreak
\noindent
{\bf Question 7.}
Deduce from the previous questions that the probability that $M(n) \ge n K C' 4^n$ is bounded by $C'' 2^{-n}$ (independently of $\delta$). 

\medbreak
\noindent
{\bf Question 8.} Conclude that the family $(\gamma^\delta, \delta > 0)$ is tight.

\chapter{SLE(6) computations}

In this lecture, we will estimate probabilities of certain events involving SLE processes, following the papers \cite {LSW2, LSW5}. We will then use these results in the next lecture in order to get some information about discrete percolation.

\section {Radial SLE}

\noindent
{\bf Radial Loewner chains.}
We now want to find a nice way to encode
``continuously increasing'' families of compact subsets $(K_t, t \ge 0)$
of the closed unit disk $\overline \U$ that are growing from the
boundary point $1$ towards $0$ 
(note that $0$ is an inner point of the unit disk, so that the situation is different from chordal Loewner chains).
As in the chordal case, we are going to
focus on the conformal geometry of the complement $U_t$ of
$K_t$ in $\U$.
Let us first look for a suitable parametrization:
It turns out to be convenient to
define the conformal map $g_t$ from
$H_t$ onto $\U$ that is normalized
by
$$
g_t (0) = 0 \hbox { and } g_t'(0) > 0 .
$$
Note that $g_t'(0) \ge 1$ and that $t \mapsto g_t'(0)$ is increasing (we leave this as a {homework exercise}).
It is therefore natural to measure the ``size''
of $K_t$ by $a(K_t) = \log g_t'(0)$. 
This suggest to consider growing families of
compact sets that have been reparapetrized in such a way that $a(K_t) = t$ i.e. that $g_t'(0)= e^t$.

Suppose now that $(\zeta_t, t \ge 0)$
is a continuous
function on the unit circle $\partial \U$.
Define for all $z \in \overline \U$,
the solution $g_t(z)$ to the ODE
\begin {equation}
\label {rODE}
\partial_t g_t (z) = - g_t (z) \frac { g_t (z) + \zeta_t}{g_t (z) -
\zeta_t}
\end {equation}
with
$g_0 (z) = z$.
This solution is well-defined up to the
(possibly infinite) time
$$
T(z) = \sup \{ t> 0 \ : \
\min_{s \in [0,t)} |g_s (z)- \zeta_s|
>0 \}.
$$
We then define
$$
K_t = \{ z \in \overline \U \ : \ T(z) \le t \}
\hbox {
and }
U_t = \U \setminus K_t .$$
The family $(K_t, t\ge 0)$ is called the radial Loewner chain
associated to the driving function $\zeta$. Note that in this case, $a(K_t) =t$.

The general statements derived in Greg Lawler's course for chordal chains
have analogues in this radial case.
There is one nice additional feature:
It is possible to estimate the Euclidean distance $d_t$ from $0$
to $K_t$ in terms of $a(K_t)=t$.
Indeed, since $U_t$ contains the disc $d_t  \U$, it follows readily that $g_t' (0) \le 1/ d_t$.
On the other hand, Koebe's $1/4$ Theorem states that
$1/d_t \le 4 g_t' (0)$.
This is loosely speaking due to the fact that the best $K_t$
can do to get as close to $0$
in ``time $t$'' is to go straight i.e. to choose $\zeta =1 $.
Hence, for all $t \ge 0$,
\begin {equation}
\label {e.koebe}
 e^{-t}/4 \le d(0, K_t) \le e^{-t}.
\end {equation}

\medbreak
\noindent
{\bf Radial SLE$(\kappa)$.}
Radial SLE($\kappa$) is then the random family of
sets $(K_t, t \ge 0)$ that is obtained when
$$
\zeta_t = \exp ( i \sqrt {\kappa} B_t )
$$
where $\kappa>0$ is fixed and $(B_t, t \ge 0)$ is
standard one-dimensional Brownian motion.

As in the chordal case, one can then define radial $SLE$
from $a \in \partial D$
to $b \in D$ in any open simply connected domain $D$
by taking the image of radial $SLE$ in $\U$ under
the conformal map $\Phi$ from $\U$ onto $D$ such
that $\Phi(1)=a$ and $\Phi (0)=b$.
Note that in this radial case, the
time-parametrization is also
well-defined since there exists
only one such conformal map (recall that in the
chordal case, one had to invoke the
scaling property to make
sure that chordal SLE in other domains than the
half-space was properly defined).

Just as for chordal SLEs, radial SLEs are the only random Loewner chains that combine 
domain Markov property and symmetry (in symmetric domains).

\section {Relation between radial and chordal SLE(6)}

Radial SLE(6) and chordal SLE(6) are in fact very closely related, and this property is specific to the case $\kappa=6$:
radial SLE(6) and chordal SLE(6) are basically the same, up to the time at which the curve disconnect the target points
of the two processes from one another. 
This is similar to the locality property of SLE(6) (Proposition \ref {locality}).
  Here is the precise statement:

\begin {prop}
\label{radchor}
Consider a chordal SLE(6) process $\gamma^1$ from $1$ to $-1$ in the unit disc $\U$, and a radial SLE(6) process $\gamma^2$ from $1$ to $0$ in $\U$.
Define
$$T^l = \inf \{t > 0 \ : \ \gamma^l [0,t] \hbox { disconnect } 0 \hbox { from } -1 \hbox { in } \U \}$$
for $l=1,2$. Then, the two paths $\gamma^1 [0,T^1]$ and $\gamma^2 [0,T^2]$ defined modulo time-reparametrization have the same law.
\end {prop}

It is possible to prove this result directly via stochastic calculus methods \cite {LSW2, Wln, Lbook}. The computations are a little messy, but in spirit, this is just a straightforward consequence of It\^o's formula (applied several times).
Note that we have derived the locality property of chordal SLE(6) in the last lecture using its relation with critical percolation. In fact, it is also possible to view this proposition as a consequence of locality of chordal SLE(6). 
Here is an outline of how to proceed:

\medbreak
\noindent
{\bf Outline of the proof of Proposition \ref {radchor}.}
Consider the chordal SLE(6) $\gamma^1$ and define, for each $t \le T^1$, the conformal map $f_t$ from $U_t$ onto $\U$ such that $f_t (0)= 0$ and $f_t (\gamma^1_t) = 1$.
Chordal SLE(6) satisfies the conformal Markov property, and $f_t$ is a conformal map. Hence,  if $\gamma^1 [0,t]$ ($t$ can be a stopping time such as the first time at which capacity seen from the origin reaches a certain value) is given and $t <T^1$, the conditional law of $f_t (\gamma^1[t,\infty))$ is that of a chordal Loewner chain from $1$ to $f_t(-1)$ in $\U$. By the locality property of chordal SLE(6), this is (at least up to some stopping time) the same as the law of chordal SLE from $1$ to $-1$.
It follows readily that if we view $\gamma^1$ as a radial Loewner chain (and time-change it accordingly), it is driven by a continuous function $\zeta_t= \exp (i \theta_t)$ on the unit circle such that $\theta_t$ has independent increments. 

But chordal SLE(6) from $1$ to $-1$ in $\U$ is symmetric with respect to the real axis. It follows that $\theta_t$ and $- \theta_t$ have the same distribution. Hence, $\theta_t$ is a Brownian motion with no drift, running at a certain constant speed $\kappa$. In other words, up to $T^1$, $\gamma^1$ coincides with some radial SLE($\kappa$). An inspection of the behavior of chordal and radial SLEs at very small times yields that this speed $\kappa$ has to be equal to $6$.

\section {Relation to discrete radial exploration}

We now define the radial exploration process for percolation in the discrete setting.
We start with a fine-mesh (we call $\delta$ this mesh) approximation of the unit disc, the unit hexagon, or more generally any other simply connected domain $D$ with $1$ as boundary point and $O \in D$. Our goal is to define a path from the boundary point $1$ to the origin. We are going to define this path dynamically.
We start with the same rule as the exploration process from $1$ in the chordal case, except that we do not fix a priori the colors of the sites on the $\partial D$. Note that as long as the discrete exploration path does not disconnect the origin from infinity, there is some arc $I$ of points on  $\partial D$ that are connected to the origin without intersecting the exploration path. When the exploration path is inside $D$, then we continue it as if  we would do the exploration process from $1$ to one of the points in $I$. When the exploration process hits $\partial D$, we turn into the direction of $0$. The rules that we just described determine the exploration path up to the first time at which $I = \emptyset$, i.e. at which it disconnects the origin from $\partial D$. At that moment, note that the connected component of the complement of the path that contains the origin is simply connected, and that it has a boundary point at distance $\delta$ of the
tip of the exploration process. We then force the exploration process to move to this point. The exploration process is then at a boundary point of the connected component that contains the origin. Now, we start again, as if the colors of the boundary of this domain would not have been known, and we start exploring interfaces in this domain using the same algorithm.

Note that when the exploration process disconnects the origin from the unit disc, it means that it has discovered a monochromatic loop around the origin. We can in fact recognize the color of the loop by looking at the winding of the exploration path. If it disconnects the origin clockwise, then the inside boundary corresponds to an open loop, and if it disconnects it anti-clockwise, then the inside boundary corresponds to a closed loop.

\begin {theorem}
When the mesh of the lattice goes to zero, then the law of the radial discrete exploration process converges to that of radial SLE(6).
\end {theorem}

\noindent
{\bf Outline of the proof.}
 Basically, one can note that the fact that up to the first disconnection time $T_1$, the law of the discrete exploration process converges to that of the radial SLE(6) up to its first disconnection time, is a combination of the fact that chordal exploration does converge to chordal SLE(6), that the discrete disconnection time converge to the continuous disconnection time, and that radial SLE(6) is identical to (some) chordal SLE(6) up to the disconnection time.
Then, one notes that the discrete new domain that one has to explore converges almost surely to the continuous domain, and one iterates the procedure.
We leave the details of the proof of this theorem as a {homework exercise}.

\medbreak

Other consequences of the same arguments, together with Koebe's $1/4$ theorem go as follows. 
\begin {itemize}
\item
Let us define the probability $J_t$ of the event ${\mathcal J}_t$ that radial SLE(6) in the unit disc, parametrized by capacity has not closed any clockwise loop before time $t$. Define for each annulus $\{ z \ : \ r  <|z|< 1 \}$ (with $r<1/2$), the probability $j^\delta (r)$ that percolation on a $\delta$-grid approximation of the annulus contains an open path joining its inner and outer boundary.
Then for any small $\delta$,
$$ {J}_{4 \log (1/r)} \le  j^\delta (r) \le {J}_{\log (1/r) / 4 }.$$

\item
Define for $z= \exp (ix)$ on the unit circle, the event ${\mathcal H}_t(x) \subset {\mathcal J}_t$ that one radial SLE(6) started from $1$ did not disconnect the point $z$ from the origin before time $t$. Define for each annulus as before the probability $h^\delta (r,x)$ that there exist two arm of opposite colors joining the inner boundary to the outer boundary, such that the open one lands in the arc $\exp (\theta), \theta \in (0,x)$ and the closed one on the other arc between $z$ and $1$.
Then, for any sufficiently small $\delta$,
$$ P({\mathcal H}_{4 \log (1/r)}(x)) \le h^\delta (r,x) \le P ({\mathcal H}_{\log (1/r) / 4 } (x)) .
$$
\end {itemize}
We again leave the detailed proofs of these facts as {homework exercises}.

\section {Exponent computations}

Our goal is to estimate  $P( {\mathcal H}_t (x))$, ${J}_t$ and some related probabilities when $t \to \infty$.
To avoid to many subscripts, we now write ${\mathcal H} (x,t) = {\mathcal H}_t (x)$.

\bigbreak
\noindent {\bf Disconnection exponent.}
 Let us start with estimating $P({\mathcal H} (x,t))$:

\begin {prop}
\label {p.disco}
There exists a constant $c$ such that for all $t \ge 1$
and for all $x \in (0, 2 \pi)$,
$$
e^{- t/4} (\sin (x/2))^{1/3}
\le P ( {\mathcal H}(x,t) ) \le c e^{-t/4} (\sin (x/2))^{1/3}.
$$
\end {prop}

\noindent
{\bf Proof.}
We will use the notation
$$
f(x,t) = P ( {\mathcal H} (x,t) )
.$$
Let
$\zeta_t = \exp ( i \sqrt {6} B_t )$
be the driving process of the radial SLE(6), with $B_0=0$.
For all $x \in (0, 2 \pi)$,
we define $Y_t^x$ the continuous function (with respect to $t$) such that
$$ g_t (e^{ix} ) = \zeta_t \exp ( i Y_t^x ) $$
and $Y_0^x=x$.
The function $Y_t^x$ is defined as long as ${\mathcal H}(x,t)$ holds.
Since $g_t$ satisfies Loewner's differential equation, we get that
\begin{equation}\label{Ysde}
d Y_t^x = \sqrt {6} \  dB_t
+ \cot(Y_t^x/2)\ dt.
\end{equation}
Let
$$
\tau^x
:=
\inf\{t \ge 0 \ : \  Y_t^x\in\{0,2\pi\}\}
$$
denote the time at which $\exp (ix)$ is absorbed
by ${K_t}$, so that
$$
f (x,t) = P ( \tau_x > t ).
$$
We therefore want to estimate the
probability that the
diffusion $Y^x$ (started from $x$)
has not hit $\{0, 2 \pi \}$ before time $t$
as $t \to \infty$.
This is a rather standard problem. The general theory of diffusion processes can be used to argue that
	$f(x,t)$ is smooth on $(0,2\pi) \times \R_+$, and It\^o's
	formula immediately shows that
	\begin {equation}
	\label {e.od}
	3 \partial_{x}^2 f + \cot (x/2)
	\partial_x f = \partial_t f.
	\end {equation}
	Moreover, comparing $Y$ with a Bessel processes
	when $Y$ is small,
	one can easily see that
	for all $t >0$,
	\begin {equation}
	\label {bv}
	\lim_{x \to 0+} f(x,t)
	= \lim_{x \to 2 \pi - } f(x,t) = 0.
	\end {equation}
	Hence, $f$ is solution to (\ref {e.od})
	with boundary values (\ref {bv}) and $f(x,0) = 1$.
	This characterizes $f$, and its long-time
	behavior is described in terms of the
	first eigenvalue of the operator
	$3 \partial_x^2 + \cot (x/2) \partial_x$.
	It turns out that the first eigenfunction is $\sin (x/2)^{1/3}$. We therefore define
	$$
	F(x,t) = E(   1_{{\mathcal H}(x,t)} \sin (Y_t^x/2)^{1/3} ). $$
	Then, it is easy to see that $F$ also solves (\ref {e.od})
	with boundary values (\ref {bv}) but this time with initial
	data $F(x,0) = \sin (x/2)^{1/3}$.
	One can for instance invoke the maximum principle
	to construct a handcraft proof (as in \cite {LSW1})
	of the fact that this characterizes $F$.
	Since $e^{- t /4} \sin (x/2)^{1/3}$
	also satisfies these conditions, it follows that
	$$
	F(x,t) = e^{-t /4 } \sin (x/2)^{1/3}.
	$$
	Hence,
	$$
	f(x,t) = P ( {\mathcal H} (x,t))
	\ge E ( 1_{ {\mathcal H} (x,t)} \sin (Y_t^x /2)^{1/3})
	=  e^{- t / 4} \sin (x/2)^{1/3}.
	$$
	To prove the other inequality, one can for instance use
	an argument based on Harnack-type considerations:
	For instance, one can see
	that (uniformly in $x$) a positive fraction of
	the paths $(Y_t^x , t \in [0, 1])$
	such that $\tau_x > 1$ satisfy
	$Y_1^x \in [\pi/2, 3 \pi /2]$.
	This then implies readily (using the Markov
	property at time $t-1$) that
	for all $t \ge 1$,
	$$
	f(x,t)
	\le c_0 P ( \tau_x > t
	\hbox { and } Y_t^x \in [ \pi/2, 3\pi /2] )
	\le c_1
	 F(x,t) =
	c_1 e^{- t/4} \sin (x/2)^{1/3}.$$

\medbreak
\noindent
{\bf Disconnection exponent II.}
Let us now focus on the probability $J_t$ that $\gamma$ has not closed any clockwise loop before time $t$:

\begin {prop}
\label {disco2}
There exist two constants $c_1$ and $c_2$ such that for all $t \ge 2$,
$$ c_1 e^{-5t/48} \le J_t \le c_2 e^{-5t/48}.$$
\end {prop}

Note that at time $t$, the boundary of $U_t$ can be decomposed into
three parts. The arc on $\partial \U$ that the path has not yet disconnected, an arc $\partial^1_t$ that corresponds to ``the left'' of the curve $\gamma$ and the part $\partial^2_t$ that corresponds to ``its right''. The point $\gamma_t$ separates the two arcs $\partial_t^1$ and $\partial_t^2$.
Note that at the times at which the path $\gamma$ completes a clockwise loop, $\partial^2_t= \emptyset$ and $\partial^1_t= \partial U_t$, whereas when it completes an anticlockwise loop, $\partial_t^1 = \emptyset$.
These arcs correspond to the scaling limit of the colors on the boundaries of the approximations of $\gamma$.

Let us define $Y_t$ the arclength of $g_t ( \partial_t^1)$. We are interested in the probability that the process $Y_t$ does not hit $2 \pi$ before time $t$ (this is exactly $J_t$). We note that $Y$ is a Markov process started from $Y_0=0$.
One can note that when $Y_t \notin \{0, 2 \pi \}$, then $Y_t$ evolves exactly as the Markov process $Y_t$ described in the previous section:
$$
d Y_t = \sqrt {6} \  dB_t
+ \cot(Y_t/2)\ dt.
$$
We define the function
$$ f(x,t)= P ( 2 \pi \notin Y [0, t] | Y_0= x ).$$
As before, it satisfies the same PDE (\ref {e.od}) for all $x \in (0, 2\pi)$ and $t \ge 0$. We know the boundary condition $f(2 \pi-, t) = 0$ and
it remains to understand what happens to the function $f(x,t)$ in the neighborhood of $x=0$. It will turn out to be convenient to define
$h(x,t) = \int_0^1 f(x,t+s) ds$. Clearly, $h(x,t) \le f(x,t) \le h(x,t-1)$, $h$ satisfies the same PDE (\ref {e.od}) and $h(2 \pi-, t) = 0$.
We will prove in a moment that $\partial_x h (0,t)=0$ for all $t >0$.
Hence, just as in the previous case we are looking for the first eigenfunction of our operator but with Dirichlet boundary condition at $2\pi$ and Neumann at $0$.
It turns out that the positive eigenfunction is $\cos (x/4)^{1/3}$ and that the corresponding eigenvalue is $5/48$. The rest of the proof of the
proposition is then identical to that of Proposition \ref {p.disco}.

\medbreak
It now remains to show that for each $t$,
$ h(x, t) - h(0,t) = o(x)$
as $x \to 0$. We are going to use the same percolation configuration in the unit disc to approximate $h(0,t)$ and $h(x,t)$. These quantities  correspond asymptotically to the
probability that the exploration process does not discover a ``white'' loop before reaching capacity $t$. In the first case, one starts with the configuration where all points on $\partial D$ are black, and in the second one, all points on $\partial D$ are black, except an arc of length $x$ near the starting point. The two exploration processes are therefore the same, but the time at which they discover the white loop might be different.

The reason for working with $h$ rather than with $f$ is the following: The probability that the second exploration process discovers a white loop strictly before the first one (i.e. using the small white arc) is of order $x$, but when it does so at some time $t$, then with high probability,
the first one will discover a white loop before time $t + o(1)$ with probability $1- o(1)$ (as $x \to 0$). So, instead of proving that the law of this discovery time is diffuse, we just randomize the time $t$ and replace it by a uniformly chosen time in the interval $[t,t+1]$ as this gets rid of this issue directly.

We note that in order for the exploration process to discover a white loop using the initial boundary arc, the white loop has to go through this boundary arc, and therefore, there exist two disjoint white paths starting on this arc, that reach distance $1/4$ (from the point $1$). If this big loop exists, then indeed, the exploration process will ``surround it'' and hit the arc at some random time $\sigma$. This two-arm estimate tells us that $P ( \sigma < t+1) = O(x)$ as $x \to 0$ and $t$ is fixed.  This is not quite enough since we want to prove that $h(x,t) - h(0,t)=o(x)$.   Suppose that $\sigma < t+1$. At $\sigma$, the exploration process is on the arc.
Let us consider two cases:
\begin {itemize}
\item
In between the two discovered white arcs starting from the arc, there is a black path that reaches distance $x^{1/2}$. Note that this black path has to be in the region that has not yet been explored. Hence, we get that for some $\alpha$ and $c$ independent of everything (and in particular of the past of the path before $\sigma$, the probability of the existence of this third arm is bounded by $c x^{\alpha/2}$.
\item
If such a black arm does not exist, then it means that the exploration process will hit the first arc at some time $\sigma'$ before reaching distance $x^{1/2}$ of the point $1$. When it does so, it does discover a white loop around the origin that does not use the initial white arc.
Furthermore, it is clear that $\sigma' - \sigma$ in this case is not larger than the capacity seen from $0$ of $\U \setminus (1-x^{1/2}) \U$ in $\U$,   So, we get that in this case, $\sigma' - \sigma \le cst \times (x^{1/2})$.
\end {itemize}
Putting the pieces together, we get that indeed
$$|h(x,t) - h(0,t) | \le \hbox {cst} \times x \times ( x^{\alpha/2} + x^{1/2}) = o(x).
$$

\medbreak
\noindent
{\bf Derivative exponents.}
	The previous arguments can be generalized as follows in order to derive the
	value of other exponents:
	 We
	will focus on the moments of the  derivative of $g_t$ at
	$\exp (ix)$ when ${\mathcal H}(x,t)$ holds.
	Note that on a heuristic level, $|g_t'(e^{ix})|$ measures
	how ``far'' $e^{ix}$ is from the origin in $H_t$.

	More precisely, we fix $b \ge 0$, and we
	define
	$$
	f(x,t)
	:=
	E \Bigl(\left|g_t'(\exp(i x))\right|^b\,
	1_{{\mathcal H}(x,t)}\Bigr) .
	$$
	We also define the numbers
	\begin {eqnarray*}
	\nonumber
	q=q(b)
	&
	:= &  \frac{1+\sqrt{1+ 24b  }}{6}
	\\
	\lambda  = \lambda (b)
	&:=&
	\frac{4b +1 +\sqrt{1 + 24 b}}{8}
	.\end {eqnarray*}
	The main result of this paragraph is the following
	generalization of Proposition~\ref {p.disco}:
	\begin {prop}
	\label {p.exp}
	There is a constant $c>0$ such that
	for all $t \ge 1$, for all $x \in (0, 2\pi)$,
	$$
	 e^{-\lambda t} \bigl(\sin(x/2)\bigr)^q
	\le f(x,t)\le c  e^{-\lambda t} \bigl(\sin(x/2)\bigr)^q
	$$
	\end {prop}

	\noindent
	{\bf Proof.}
	We can assume that $b>0$ since the case $b=0$ was treated in
	the previous section.
	Let $Y^x_t$ be as before
	and define  for all $t < \tau^x$
$$
\Phi_t^x :=
\left| g_t'(\exp(ix))\right|
\,
.$$
On $t\ge \tau^x$ set $\Phi_t^x := 0$.
	Note that on $t<\tau^x$,
	$\Phi_t^x =\partial_x Y_t^x$ and
	$$Y_t^x = \sqrt {6} B_t
+ \int_0^t \cot(Y_s^x/2)\ ds.$$
Hence, we
	that for $t<\tau^x$
	\begin{equation}\label{e.logder}
	\partial_t \log \Phi_t^x = -
	\frac 1{2\sin^2(Y_t^x/2)}
	\end{equation}
	so that
	\begin{equation}\label{Phiis}
	(\Phi_t^x)^b =
	\exp
	\left( - \frac {b}{2}
	\int_0^t \frac {ds} {\sin^{2}(Y_s^x/2)} \right)
	\ ,
	\end{equation}
	for $t<\tau^x$.
	Hence,
	$$
	f(x,t) = E \left( 1_{{\mathcal H}(x,t)}
	\exp
	\left( - \frac {b}{2}
	\int_0^t \frac {ds} {\sin^{2}(Y_s^x/2)} \right)
	\right)
	.$$
	Again, it is not difficult to see
	that
	 the right hand side of~(\ref{Phiis}) is $0$
	when $t=\tau^x$
	and that
	\begin{equation}\label{sidebd}
	\lim_{x\to 0} f(x,t)=\lim_{x\to 2\pi} f(x,t)=0
	\end{equation}
	holds for all fixed $t>0$.

	Let $F:[0,2\pi]\to\R$ be a continuous function with
$F(0)=F(2\pi)=0$, which is smooth in $(0,2\pi)$, and set
$$
h(x,t)=h_F(x,t):=
E \Bigl((\Phi_t^x)^b\,F(Y_t^x)\Bigr).
$$
By~(\ref{Phiis}) and the general
theory of diffusion Markov processes,
we know that $h$ is smooth in $(0,2\pi)\times\R_+$.
The Markov property for $Y_t^x$ and~(\ref{Phiis}) show
that $h(Y_t^x,t'-t) (\Phi_t^x)^b$ is a local
martingale on $t<\min\{\tau^x,t'\}$.  Hence,
the drift term of the stochastic
differential $d\bigl(h(Y_t^x,t'-t) (\Phi_t^x)^b\bigr)$
is zero at $t=0$.  By It\^o's formula, this means
that
\begin{equation}\label{ekol}
\partial_t h =
\frac \kappa 2\, \partial_x^2 h
+\cot(x/2)\, \partial_x h
-\frac b {2 \sin^2 (x/2)}\, h
\,.
\end {equation}
The corresponding positive eigenfunction is
$ \bigl(\sin(x/2)\bigr)^q$. We therefore define $F$ to be this function, so that
$F(x) e^{-\lambda t} = h_F$ because both satisfy~(\ref{ekol}) on
$(0,2\pi)\times [0,\infty)$ and have the same boundary values.

Finally, one can conclude the proof using the same type of argument
as for Proposition \ref {p.disco}.

\chapter{Arm exponents}

\section {Some notations}

We will now state and prove results concerning critical percolation on the triangular lattice with lattice-mesh one.
Let us first tune our definitions:
For each $r$, $\Sigma_r$ denote a discrete approximation on this fixed lattice of the closed disc $r \partial \U$ of radius $r$ (we will for instance just take the set of sites on the triangular lattice that are contained in $r \overline \U$).
Recall that $\Lambda_n$ was used to denote  consisting of the sites that are at graph distance $n$ or less from the origin.

The set $\partial_r$ will be used to denote the set of sites that are on the boundary of $\Sigma_{r}$ (i.e. the set of sites in $\Sigma_r$ that have a neighbor outside of $\Sigma_r$.
We also use the set ${\mathcal A} ( r_1, r_2) = \Sigma_{r_2}\setminus \Sigma_{r_1}$ for $r_1 < r_2$ and its inner boundary $\partial_{r_1}'$ (i.e. the set of points outside of $\Sigma_{r_1}$ that have a neighbor in $\Sigma_{r_1}$).

We will often refer to $\partial_{r_1}$ and $\partial_{r_2}$ as the inner and outer circles (or boundaries) of the annulus ${\mathcal A}(r_1, r_2)$ (even if they are in fact hexagons).

We now define for all $r_1 < r_2$:
\begin {itemize}
\item The probability $\pi (r_1, r_2)$ of the event $\Pi(r_1,r_2)$ that there exists an open path joining $\partial_{r_1}'$ to $\partial_{r_2}$.
\item For each even $j$, the probability $\pi_j (r_1, r_2)$ of the event $\Pi_j (r_1, r_2)$ that there exist $j$ disjoint paths of alternating type (one open, then one closed, then one open etc.)  in ${\mathcal A}(r_1, r_2)$  joining the two boundary circles $\partial_{r_1}'$ to $\partial_{r_2}$, when ordered clockwise, say (but since $j$ is even, this does not need to be specified...).
\end {itemize}

The goal will be to derive the asymptotics of the quantities $\pi (0,n)$ and $\pi_j (0, n)$ as $n\to \infty$. The constraint to consider only the case where $j$ is even\footnote {when $j$ is odd, the colors of the crossings cannot all alternate because of parity -- anyway, it can be shown that the order of the colors does not matter} is not necessary, but it makes life a little simpler. Since our prime objective is to study $\pi$ and $\pi_4$ because these are the two quantities that will enable to make contact with near-critical percolation, we will restrict ourselves to these two cases.
Note that you have already seen in the first exercise session that $c_1 / n^2 \le \pi_5 (0,n) \le c_2 / n^2$ by arguments based on Russo-Seymour-Welsh theory only.

As we will see in the next lecture, the four-arm exponent plays a special role, as it is the one that controls the fluctuations of the macroscopic connectivity properties when one lets the percolation probability vary near to its critical value.

\section {One-arm exponent}

This case is the simplest. One reason is that the existence of one crossing is an increasing event, as opposed to the multiple crossing of different colors, so that the Harris inequality can be directly applied. Establishing the following theorem was one of the goals of these lecture series:

\begin {theorem}
For critical percolation on the triangular lattice, $P( 0 \leftrightarrow \partial \Lambda_n) = n^{-5/48 + o(1)}$ as $n \to \infty$.
\end {theorem}

Recall that we have proved in the previous lecture that for all $R > 1$,
$ \pi ( n, nR )$ converges as $n \to \infty$ to a non-disconnection probability $f(R)$ for the SLE(6) process (the probability that an SLE(6) in the unit disc does not close a clockwise loop before reaching the disc $\U / R$), and that for two constants $c_1$ and $c_2$, and all $R$ sufficiently large,
$$ c_1 R^{-5/48} \le f (R) \le c_2 R^{-5/48}.$$

We can note that on the one hand, for all $r_1 < r_2 < r_3$, one obviously has
$$ \pi (r_1, r_3 ) \le \pi (r_1, r_2) \pi (r_2, r_3)$$
because of the independence of $\Pi(r_1, r_2)$ and $\Pi(r_2, r_3)$.
On the other hand, Russo-Seymour-Welsh theory and the Harris inequality imply readily that for some absolute constant $c$ and
for all $r_1, r_2, r_3$ with $2r_1 < r_2 < r_3 / 2$,
$$ \pi (r_1, r_3) \ge c \times \pi (r_1, r_2) \pi (r_2, r_3)$$
Just note that  if $\Pi (r_1, r_2)$ and $\Pi (r_2, r_3)$ hold, and if there exist two open loops surrounding the origin in the annuli ${\mathcal A} (r_2 / 2 , r_2) $ and
${\mathcal A}(r_2,2r_2)$
 (the probabilities of these two events are bounded from below by RSW), then $\Pi (r_1, r_3)$ holds as well.

We now fix any small positive $\epsilon$. We then choose $R$ large enough so that
$$R^{-5/48 - \epsilon}  \le f(R) \le R^{-5/48 + \epsilon} \hbox { and  }R^{-\epsilon} \le c$$ for the previous constant $c$.

Then we get that for this fixed $R$,
$\pi ( R^j, R^{j+1} ) \to f(R)$ as $j \to \infty$, and hence, there exists $j_0$ so that for $j \ge j_0$,
$$
R^{-5/48 - 2 \epsilon} \le \pi (R^j, R^{j+1}) \le R^{-5/48 + 2 \epsilon}.$$

For each large $n$, we choose $j_1$ so that $R^{j_1} \le  n  <  R^{j_1  + 1}$. Then, we have on the one hand that
\begin {eqnarray*}
\pi (0, n) &\le&  \pi (0, R^{j_1})
\\ &\le& \pi (0, R^{j_0}) \prod_{j=j_0}^{j_1 -1} \pi (R^j, R^{j+1}) \\
&\le & \pi (0, R^{j_0}) \times  (R^{-5/48 + 2 \epsilon })^{j_1 - j_0} \\
& \le & \left( \pi (0, R^{j_0}) \times (R^{-5/48 + 2 \epsilon})^{-j_0 -1} \right) \times (R^{-5/48 + 2 \epsilon})^{j_1+1} \\
& \le & \hbox {cst} \times n^{-5/48 + 2 \epsilon }
\end {eqnarray*}
where the constant depends on $\epsilon$ but not on $n$.

On the other hand,
\begin {eqnarray*}
 \pi (0, n) &\ge& \pi (0,R^{j_1+1}) \\
  &\ge& \pi (0, R^{j_0+1}) \prod_{j=j_0+1}^{j_1} (c \times \pi (R^j, R^{j+1})) \\
&\ge & \pi (0, R^{j_0 +1}) (R^{-5/48 - 3 \epsilon })^{j_1 - j_0} \\
& \ge & \left( \pi (0, R^{j_0 +1}) \times (R^{-5/48 -3 \epsilon})^{-j_0} \right) \times \left(R^{-5/48 -3 \epsilon}\right)^{j_1} \\
& \ge & \hbox {cst} \times n^{-5/48 -3 \epsilon }
\end {eqnarray*}
and this concludes the proof of the theorem.

\section {Four-arms exponent}

We now want to exploit our SLE computations in order to derive the corresponding result for $\pi_4$:

\begin {theorem}
For critical site-percolation on the triangular lattice,
$\pi_4 (0, n) = n^{-5/4 + o(1)}$ as $n \to \infty$.
\end {theorem}
As we shall see, the proof is more complicated.

\medbreak
\noindent
{\bf Rephrasing the SLE computation.}
Recall from last lecture that we did compute the asymptotic behavior of the quantity
$ E (1_{{\mathcal H}(x,t)} | g_t' (e^{ix}) |)$ for radial SLE(6) started from $1$, where ${\mathcal H}(x,t)$ denotes the event that the curve up to time $t$ does not disconnect $e^{ix}$ from the origin. We showed that this quantity is (up to multiplicative constants) of the order of $e^{-5t/4}$. If we choose the point $e^{ix}$ randomly and uniformly on the unit circle, we get that
$$ c_1 e^{-5t/4} \le E \left( \int_0^{2\pi} dx | g_t' (e^{ix})| 1_{{\mathcal H}(x,t)} dt \right)
\le c_2 e^{-5t/4}$$
for large $t$ and absolute constants $c_1$, $c_2$.
Note that the quantity
$$l_t= \int_0^{2\pi} dx | g_t' (e^{ix})| 1_{{\mathcal H}(x,t)}$$
 is just the arclength of the image under $g_t$ of the set of points on $\partial \U$ that are not disconnected from the origin at time $t$. In other words, this is (up to a factor $2 \pi$) the harmonic measure (seen from the origin) of this arc i.e.  the probability that  a Brownian motion started from the origin hits $\partial \U$ before $\gamma[0,t]$.

\begin{figure}
\centerline{\includegraphics*[height=2in]{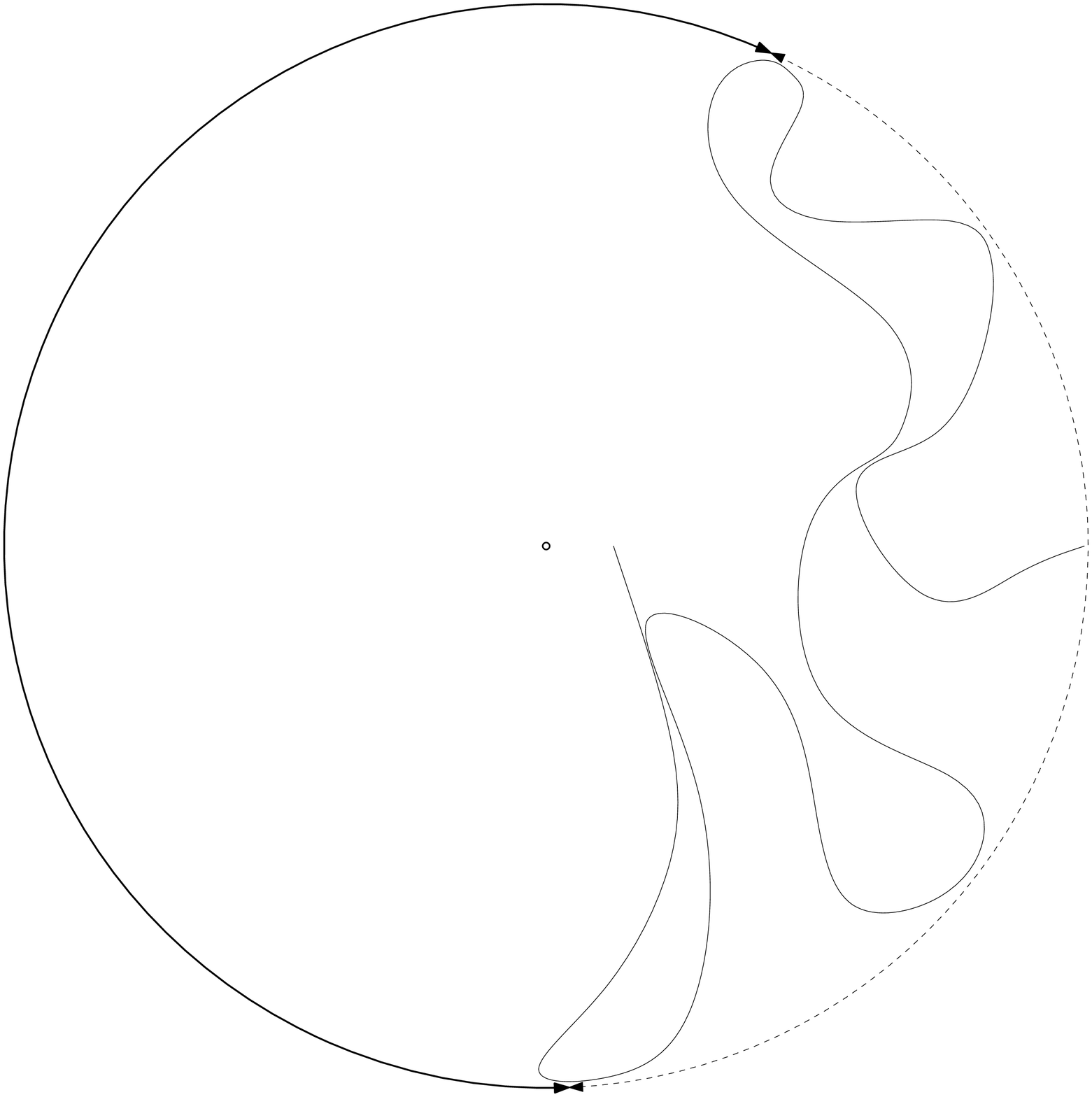}}
\caption {$l_t$ is the harmonic measure at the origin of $\partial \U$ in $\U \setminus \gamma [0,t]$.}
\end {figure}

We would like to argue that $l_t$ is comparable to the (asymptotic) probability of the existence of two further arms joining a small circle to $\partial \U$ in the domain left not yet explored by the SLE at time $t$. The following argument is rather easy and based on classical complex analysis only:

We first quickly recall the definition of extremal distance between two subsets $A$ and $A'$ in a (non-necessarily simply connected) domain $D$. This is the supremum over all smooth metrics (i.e. corresponding to the Lebesgue measure with piecewise smooth density) that give unit area to $D$, of the squared distance (in $D$) between $A$ and $A'$ (i.e. the square of the length of the shortest path in $D$ for this metric that joins these two sets). It is clear from this definition that this quantity is conformally invariant, and it is easy to check (homework exercise) that the extremal distance between the two vertical sides of an $a \times b$ rectangle is just $a/b$ (an optimal metric is a multiple of the Lebesgue measure). Let us call $\pi \ell (A,A', D)$ this extremal distance (we use this normalization in order not to have $\pi$'s later on). Mind the different notation $l$ and $\ell$! Note that clearly, $\ell$ is a decreasing function of $D$ when $A$ and $A'$ are fixed (there are more possible metrics to choose from) and a decreasing function of $A$ and $A'$ as well.  Furthermore, if $A''$ disconnects $A$ from $A'$ in $D$, then $\ell (A, A'', D ) \le \ell (A, A', D)$.

Suppose that we choose a small $r$. Define the path $\gamma$ up to its first hitting time $\tau$ of the circle of radius $r$. By Koebe's $1/4$ theorem, we know that $  \log (1/4r) \le   \tau \le \log (1/r)$.
We then define the extremal distance $\ell_r$ between $r\partial \U$ and $\partial \U$ in $\U \setminus \gamma [0, \tau_r] $.
The previous monotonicity remarks imply that
$$
\ell \left( (r/8) \partial \U , \partial \U, \U \setminus \gamma [0,\log (1/r)] \right) \le \ell_r
\le
\ell \left( r \partial  \U, \partial \U, \U \setminus \gamma [0, \log (1/8r) ]\right).
$$
We therefore define for all $t$,
$$ \tilde \ell_t = \ell \left( \frac {e^{-t}}{8} \partial \U, \partial \U, \U \setminus \gamma [0, t]\right).$$

\begin{figure}
\centerline{\includegraphics*[height=2in]{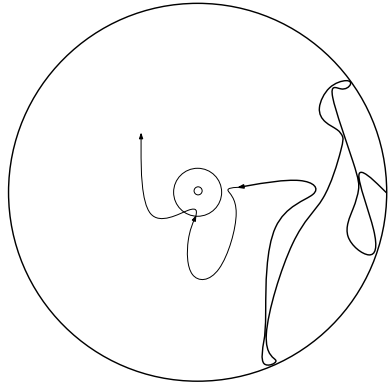}}
\caption {$\gamma$ up to the times $\log (1/8r)$, $\tau_r$, $\log (1/r)$ and the two circles $r\partial \U$ and $r \partial \U / 8$.}
 \end {figure}

As extremal distance is conformally invariant,
$$ \tilde \ell_t = \ell \left( g_t \left( \frac {e^{-t}}8 \partial \U \right), L_t, \U \right)$$
where $L_t$ is the arc $g_t (\partial \U)$ of length $l_t$. Furthermore, and we leave this as a homework exercise, it is easy to check that
$$ \frac 1 {32} \U \subset g_t \left(\frac {e^{-t}}8 \U \right) \subset \frac 1 2 \U$$
because the distance $d_t$ between the origin and $\gamma [0,t]$ satisfies  $e^{-t} /4  \le d_t \le  e^{-t}$
and that if $L(l)$ denotes an arc of length $l$ on $\partial \U$, then the three quantities
$l$, $\exp ( - \ell ( \U /32, L(l), \U) ) ) $ and $\exp ( - \ell ( \U/2, L(l), \U))$  are comparable (up to multiplicative constants)
when $l \to 0$.

We can then put the pieces together and conclude that there exist two constants $c_1'$ and $c_2'$ such that for all small enough $r$,
$$  c_1' r^{5/4} \le E (e^{-\ell_r}1_{{\mathcal H}})  \le c_2' r^{5/4},$$
where ${\mathcal H}$ is the event that $\gamma [0,\tau_r]$ did not disconnect the origin from infinity.

\medbreak
\noindent
{\bf Relation to crossing probabilities.}
We wish to understand the asymptotic behavior of the probabilities  $\pi_4 (rn, n)$ as $n \to \infty$. Note that because the number of arms is even and the colors are alternating, any point on the boundary of the annulus will be ``between two arms of opposite color''. By symmetry, if we fix a boundary point for instance $x=r_2$, and we look at the probability $\hat \pi_4 (r_1, r_2)$ that there exists 4 arms (open, closed, open, closed) ordered in this way when starting clockwise from the boundary point, then
$$ \hat \pi_4 (r_1, r_2) \le \pi_4 (r_1, r_2) \le 2 \hat \pi_4 (r_1, r_2).$$
Hence, we can estimate $\hat \pi_4$ instead of $\pi_4$.

The event corresponding to $\hat \pi_4 (r_1, r_2)$  can be explored as follows: Start a radial exploration from this particular boundary point $r_2$, and stop it when it reaches
$\partial_{r_1}'$. The exploration is not allowed (i.e. if it does, we stop) to disconnect the inner boundary from the outer boundary as it would prohibit the existence of arms. Then in the remaining domain (i.e. in the simply connected  component ${\mathcal U}$ of the complement of the path in the annulus that has points of $\partial_{r_1}'$ and of $\partial_{r_2}$ on its boundary, there must still exist two arms. This can then be explored by a chordal exploration process in this simply connected domain.

The convergence to SLE(6) results that we have derived (together with the fact that discrete explorations disconnection/hitting times converge to those of the corresponding SLE(6)) show that when $n \to \infty$, the probability $\hat \pi_4 (rn, n)$
converges as $n \to \infty$ to some function $F(r)$ that can be viewed as the probability that:
\begin {enumerate}
\item
A radial SLE started from $1$ and stopped at its first hitting time $\tau_r$ of  $ r\U$, does not disconnect the origin from $\partial \U$.
\item
We then define the connected component ${\mathcal U}$  as before in the discrete case. In this domain, we define a chordal SLE(6) that goes from one of the two boundary points corresponding to  $\gamma( \tau_r)$ to the other.  We want this SLE to successively hit $\partial \U$ and $r \partial \U$ before disconnecting the two boundaries from one another.
\end {enumerate}


Conditionally on $\gamma [0, \tau_r]$, the probability of this last event is just the asymptotic probability of the existence of two arms joining the two boundary arcs corresponding to the boundaries of  ${\mathcal U}$. This is a conformally invariant quantity, and therefore a deterministic function $m(\cdot)$ of the extremal distance between these two boundary arcs in ${\mathcal U}$. Thanks to the two-arm estimate of the first exercise sheet, we get that $ m (\ell)  e^{\ell}$ remains bounded and bounded away from $0$ when $\ell$ is large.

Putting the pieces together, we get that for some constants $c_1''$ and $c_2''$, and all small $r$,
$$
c_1'' r^{5/4}  \ge \lim_{n \to \infty} \hat \pi_4 (rn, n) \ge
c_2'' r^{5/4}.
$$

\medbreak
\noindent
{\bf Separation lemma.}
We would like to adapt the one-arm proof to this four-arm case. As we have already mentioned, there is the problem that $\pi_4$ does not correspond to increasing events. We therefore can not apply the Harris inequality to get a lower submultiplicative bound anymore.
We have to replace this by the idea of defining well-separated configurations, that we then can extend and glue together.

\medbreak

Let us fix $r_1 < r_2$. For any small $\delta>0$, we denote by $\pi^\delta_4 = \pi^\delta_4 (r_1, r_2)$ the probability of the event
that:
\begin {itemize}
\item Just as for $\hat \pi_4 (n_1, n_2)$, there exist four disjoint crossings of ${\mathcal A} (r_1, r_2)$ of alternating color, ordered in a prescribed way when starting from a prescribed point $x_0 \in \partial_{r_1}' \cup \partial_{r_2}$.
\item For any point $x$ on $\partial_{r_1}'$ (respectively on $\partial_{r_2}$), there are no three disjoint paths from $x + \delta \partial_{r_1}$ to
$\partial_{2r_1}$ (resp. from $x+ \delta \partial_{r_2}$ to $\partial_{r_2/2}$), when one restricts the percolation to the annulus ${\mathcal A}(r_1, r_2)$.
\end {itemize}

Note that for symmetry reasons, this probability does not depend on the choice of $x_0$.

One reason to define this notion is the following lemma:

\begin {lemma}
\label {19}
There exists a constant $c=c(\delta)>0$ such that for all large enough $r_1< r_2/4 <r_3/32$,
$$ \pi_4^\delta (r_1 , r_3) \ge c(\delta) \pi_4^\delta (r_1, r_2) \pi_4^\delta (2r_2, r_3).$$
\end {lemma}

The proof is again an application of the ``progressive exploration'' ideas. Here is a quick sketch -- the important part is the picture...
We explore the event corresponding to $\pi^{\delta}_4 (r_1, r_2)$ progressively: We start with an exploration process starting from any given point on the inner circle $\partial_{r_1}'$. Then, the exploration process reaches $\partial_{r_2}$ and we stop it. There, we have discovered two crossings (one open, one closed) neighboring this hitting point, provided that the path did not disconnect the origin. This means in particular that no other crossing will be allowed to land the $\delta r_2$-neighborhood of this landing point. Then, we start again from the inner circle a second exploration process to find the closest open path on the other side of the closed path. Once we are done, we start afresh to discover a fourth crossing. In this way, we have explored some regions, as indicated in the figure, and the end-points $A_1$, $A_3$ and $A_3$ of our explorations are all at distance $\delta r_2$ from each other. Note the colors of the boundaries of the unexplored regions that are close to the circle of radius $r_2$.

Finally, we look at all the undiscovered sites in ${\mathcal A}(r_1, 2r_1)$ to check that there the part of the conditions defining $\pi_4^\delta$ near to the inner circle is satisfied. Note that we do not check whether the event about no three arms in the neighborhood of the outer circle is fully satisfied (we just collect some partial information due to the explored zone): We therefore test if an event $\bar \Pi$ of probability larger than $\pi^\delta_4 (r_1, r_2)$ holds.

\begin{figure}
\centerline{\includegraphics*[height=2in]{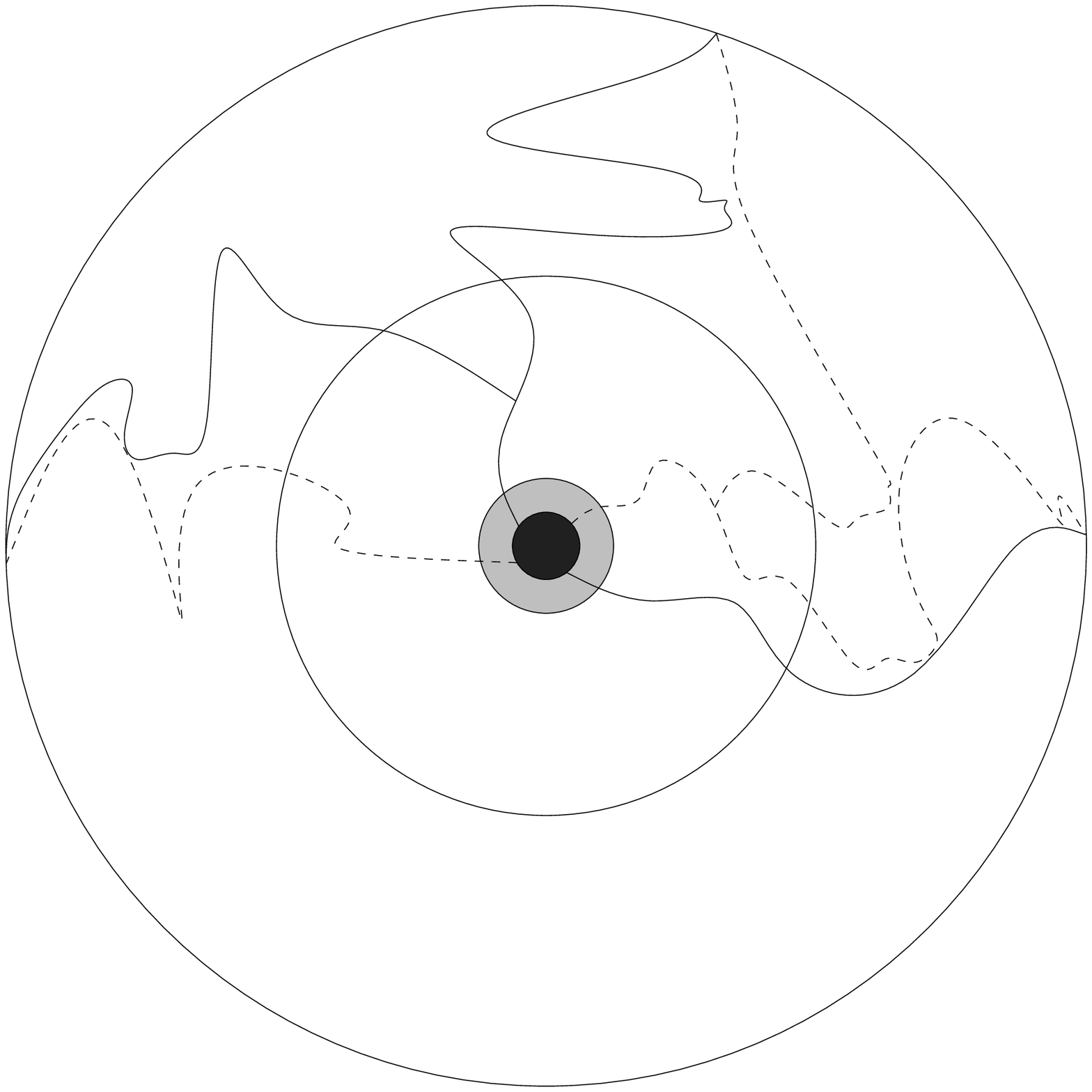}}
\caption {The event $\overline \Pi$}
\end {figure}

Now, we do the same exploration procedure in the larger annulus ${\mathcal A}(2r_2, r_3)$, except that we do it ``inwards'': We start from the outer circle, and explore inwards, and then we check the no-three-arm condition in the annulus ${\mathcal A} ({r_3/2}, r_3)$ and we leave some uncertainty about the validity
of the ``no-three-arm condition'' in the neighborhood of $\partial_{2r_2}'$.

We then leave it as a {homework exercise} (hint: see the sketchy pictures) to check that if these two events (of probability larger than $\pi^\delta_4 (r_1, r_2)$ and $\pi^\delta_4 (2r_2, r_3)$ respectively) are satisfied, then we can use the Russo-Seymour-Welsh estimates in the unexplored region and in the annulus between $\partial_{r_2}'$ and $\partial_{2r_2}$  to create, with positive probability $c(\delta)$ a configuration corresponding to $\pi_4^\delta (r_1, r_3)$.

\begin{figure}
\centerline{\includegraphics*[height=2.7in]{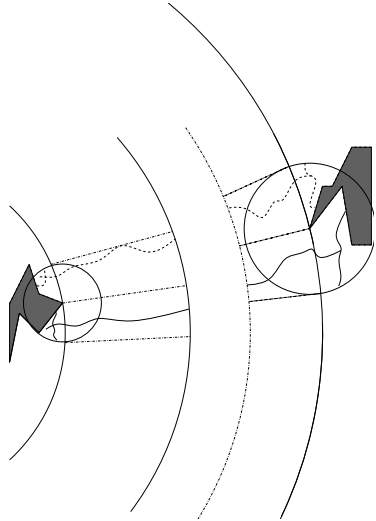}}
\caption {Extending using Russo-Seymour-Welsh near the landing points}
\end {figure}

\begin{figure}
\centerline{\includegraphics*[height=2.7in]{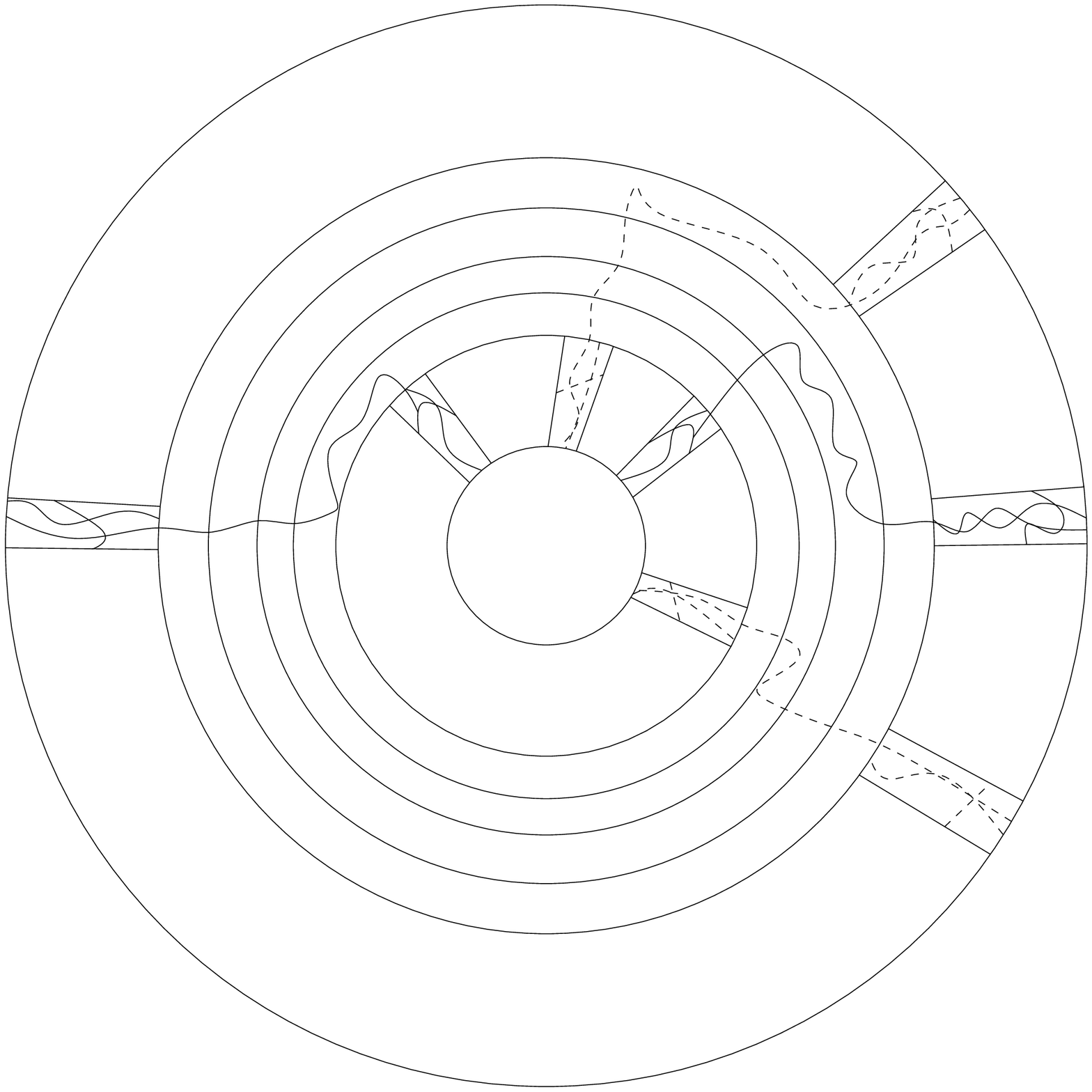}}
\caption {And then gluing using Russo-Seymour-Welsh}
\end {figure}

\medbreak

In order for this lemma to be useful, we have to see that $\pi^\delta_4$ and $\hat \pi_4$ are comparable. This is the goal of the next lemma:

\begin {lemma}
\label {20}
If one chooses $\delta$ small enough, then for all large enough $R$, and for all large enough $n \ge N(R,\delta)$,
$$\pi^\delta_4 (n, nR) \ge \frac 1 2 \hat \pi_4 (n , nR).$$
\end {lemma}

The proof is based on the one hand on the fact that that for $r_1 < r_2 /8$,
$$
\hat \pi_4 (r_1, r_2) - \pi^\delta_4 (r_1, r_2) \le \hat \pi_4 (2r_1, r_2 /2) \times a(\delta),
$$
where the probability $a(\delta)$ that the three-arms conditions fail to be satisfied near $\partial_{r_2}$ or near $\partial_{r_1}'$ goes to $0$ with $\delta$ uniformly with respect to $r_1$ and $r_2$ (because of the a priori bound on three-arms event in a half-plane). On the other hand, we see that
for large $R$,
$$ \lim_{n \to \infty} \hat \pi_4 (n,nR) \ge c_2'' R^{-5/4} \ge \frac {c_2''}{c_1'' 4^{5/4}}  (R/4)^{-5/4} \ge c_3'' \lim_{n \to \infty} \hat \pi_4 (2n, nR/2).$$
In particular, we get that for large $R$, and any large enough $n$,
$$ \hat \pi_4 (n,nR) \ge c_3''  \hat \pi_4 (2n, nR/2).$$
Now, we just choose $\delta$ in such a way that $a(\delta) \le c_3''/2$.
Then, for large $n$, we get that
$$
\pi^\delta_4 ( n, nR) \ge \hat \pi_4 (n , nR) - a(\delta) \hat \pi_4 (2n, nR/2) \ge \hat \pi_4 (n,nR) - \frac 1 2 \hat \pi_4(n,nR).
$$

\medbreak
\noindent
{\bf End of the proof.}
We are now ready to put the pieces of the puzzle together.
Let us fix any small positive $\epsilon$. Choose also a small $\delta$ so that the previous lemma holds. Then, we choose $R$ sufficiently large so that:
\begin {itemize}
\item
Lemma \ref {20} holds
\item
$
R^{-5/48 - \epsilon } \le \lim_{n\to \infty} \hat \pi_4 ( n, nR)   \le \lim_{n \to \infty}
\hat \pi_4 (2n, nR) \le R^{-5/48 + \epsilon}$.
\item
For the $2R^{- \epsilon}$ is smaller than the  $ c(\delta)$ defined in Lemma \ref {19}.
\end {itemize}
Then, we choose $j_0$ so that for all $j \ge j_0$,
$$
R^{-5/48 - 2 \epsilon} \le \pi_4 (R^{j-1}, R^{j}) \le \pi_4 ( 2R^{j-1}, R^j) \le R^{-5/48 +2 \epsilon}
$$
and
$$ \pi_4^\delta (2R^{j-1}, R^j) \ge \hat \pi_4 ( 2 R^{j-1}, R^j ) /2 .$$

Finally, we note that, just as in the case of the one-arm exponent,  on the one hand, for all $j_1 \ge j_0$ and $R^{j_1} \le n \le R^{j-1 +1}$,
\begin {eqnarray*}
\hat \pi_4 (0, n) & \le & \hat \pi_4 (0, R^{j_1}) \\
&\le & \hat \pi (0, R^{j_0}) \times \prod_{j=j_0}^{j_1} \hat \pi_4 (R^{j-1}, R^{j}) \\
&\le& c(\delta, \epsilon, j_0, R) \times (R^{-5/4 + 2 \epsilon})^{j_1 + 1}
\\
&= &c(\delta, \epsilon, j_0, R) \times n^{-5/4 + 2 \epsilon}
\end {eqnarray*}
and that on the other hand
\begin {eqnarray*}
\hat \pi_4 (0, n)
&\ge& \hat \pi_4 (0, R^{j_1+1}) \\
& \ge & \pi_4^\delta (0, R^{j_0+1}) \prod_{j=j_0+1}^{j_1} (c(\delta) \pi_4^\delta (2R^{j}, R^{j+1})) \\
&\ge& \pi_4^\delta (0, R^{j_0+1}) \prod_{j=j_0+1}^{j_1} (\frac {c(\delta)}{2} \hat \pi_4 (2R^{j}, R^{j+1})) \\
&\ge & \pi_4^\delta (0, R^{j_0 +1}) \prod_{j_0+1}^{j_1} (R^{- \epsilon} \times R^{-5/48 - 2 \epsilon} ) \\
&\ge&   c(\delta, \epsilon, j_0, R) \times n^{-5/4 + 3 \epsilon}.
\end {eqnarray*}
This concludes the proof of the theorem.

\section {Other exponents and bibliographical remarks}

\begin {itemize}
\item
The proofs presented here follow the ideas of the papers \cite {LSW5} and \cite {SmW}. In the latter one, the derivation of the other arm exponents in the plane and in the half-plane is also discussed.
 For $j \ge 2$ arms not all of the same colors, the value of the exponents turns out to be $a_j= (j^2 - 1) /12$.
The case $j=2$ with exponent $1/4$ can be derived using exactly the same arguments as in this section, based on the corresponding computation in the last lecture. In fact, the method that we presented here generalizes fairly smoothly to all cases. The only problem is that one has to use the estimates of the other moments of $|g_t' (e^{ix})|$ for SLE(6), and that the rephrasing part requires additional arguments.

\item It can be noted (see e.g. \cite {ADA}) that changing the
  prescribed order of the colors around the circle does not change the
  exponents, and that in fact, one can also change the prescribed
  colors as long as the paths are not all of the same color. The
  monochromatic multi-arm exponents are different (and so far there
  have not been shown to take any specially nice value) than the
  polychromatic ones (see \cite {LSW5} for a discussion of the
  monochromatic two-arm exponent).  The reason is that each time one
  launches exploration process in the annulus and reached the other
  side and that discovered the $j$th arm, then one can switch all
  colors in the complementary domain without changing the probability
  measure. But we need the first exploration process necessarily
  discovers two different colors. This also ``explains'' why the
  one-arm exponent's value $5/48$ does not show up naturally as $a_1$
  in the previous formula.

\item
These values for arm-exponents had been predicted by physicists using various (mathematically) non-rigorous methods (Conformal field theory, Coulomb gas, simulations, quantum gravity), see e.g. \cite {dN, N2, SRG, GA, SD, ADA, Dup}.

\item
It can be shown in various ways \cite {LW, LSW1, LSW2, LSWr,W}, see also lecture notes \cite {Wln,Wln2,Wln3} that Brownian motions and percolation clusters share the same critical exponents and that in fact \cite {LSWr, W}, their outer boundaries have exactly the same law, so that the exponents that we just have computed correspond also to the Brownian intersection exponents derived in \cite {LSW1, LSW2}.

\item Apart from the exponents that we can derive directly from
  Russo-Sey\-mour-Welsh (see the first exercise sheet), it is an open
  problem to show that $ \pi_j (j, n) n^{a_j}$ remains bounded and
  bounded from below as $n \to \infty$. The SLE based method does for
  the moment only provide the weaker statements at exponents level
  only.

\end {itemize}

\chapter{Near-critical percolation}

In this lecture, we show, based on Kesten's paper \cite {Ke}, how to
derive the behavior of $\theta (p)$ for $p$ close to the critical
value $1/2$ from the existence and values of the one-arm exponent and
the four-arm exponent that we have established in the previous
lectures.  The parameter $p$ will now vary in $[1/2, 1]$ and $P_p$
will denote the probability measure for percolation with parameter
$p$.

Due to lack of time, the proofs in this lecture will be even sketchier than in the 
other ones and we will omit some details, but we believe that
it will be enough to understand what is going on. For a more detailed presentation of the arguments, we refer to \cite {No}.

Our main goal will be to discuss what remains to be done in order to
prove the following statement (recall that $\theta (p)$ denotes the probability for percolation with parameter $p$ that the origin is in the infinite cluster):

\begin {theorem}
\label {5/36}
When $u \to 0+$, $\theta (1/2 + u )  = u^{5/36 + o(1)}$.
\end {theorem}

\section {Correlation length}

Let us define the event $H(n)$ that there exists an horizontal open crossing of the $2n \times n$ parallelogram. We define
$h_p (n) = P_p (H(n))$ the probability that such a crossing exists for percolation with parameter $p$.

For each small $\epsilon>0$ and $p > 1/2$, we define
$$ L(p, \epsilon) = \inf \{n \ : \ h_p(n) \ge 1 - \epsilon \}.$$
The fact that $L(p, \epsilon )$ is finite follows readily from the exponential decay of the connectivity probabilities by closed sites (because the set of closed sites is a sample of subcritical percolation).

Note that often, the length $L$ is defined in terms of crossings of rhombi (and this does not change its value drastically), but for our purposes, it will be simpler to work with this definition (otherwise, we would for instance need to use
and prove  an alternative version of the Russo-Seymour-Welsh formula that also shows that when the probability of crossing of a $n \times n$ rhombus is very close to one, then so is that of a $2n \times n$ parallelogram).

It follows that the probability that there exists a vertical closed crossing of the $2L(p) \times L(p)$ parallelogram for $P_p$ is smaller than
$\epsilon$.
It is easy to deduce that the probability that there exists in the $3L \times 3L$ rhombus a closed path of diameter greater than $2L$ (for the measure $P_p$ and $L=L(p)$) is not larger than a constant times $\epsilon$ (because if such a path exists, then it means that one out of a finite fixed number of $2L \times L$ parallelograms is crossed by a closed path).
From this, it is not difficult  (and we leave this as a homework exercise\footnote {There are various ways to proceed: One possibility would just be to sum over all sequences of $m/5$ non-overlapping $3L \times 3L$ rhombi centered at points $x_0=0, \ldots , x_{m/5}$ with coordinates that are multiple of $L$ and with $d(x_j, x_{j+1}) \le 6 L(p)$, of the probability that in each of the rhombi there exists a closed path of length at least $2L(p)$}) to deduce that {\em if $\epsilon$ has been chosen small enough}, then there exist absolute constants $c_1$ and $c_2$ such that for all $p$, for all $m \ge 1$,
the probability (for $P_p$) that there exists a closed path of diameter at least $m L(p)$ passing in the $(L(p)\times L(p)$ box centered at the origin is bounded from above by
$c_1 \exp (-c_2 m )$.

Until the end of the lecture, we will suppose that $\epsilon$
has been chosen in such a way.  We will keep this parameter fixed and we shall just denote $L(p, \epsilon)$ by $L(p)$.

\medbreak

A straightforward consequence of this exponential decay is that with positive probability (bounded from below independently from $p$), there is no closed loop around the origin that has a diameter greater than $L(p)$ (we will call this event ${\mathcal E}(L(p))$). This is simply due to the fact that the sum over $m$ of the probability that there exists a closed path of diameter at least $\max (mL(p))$ that goes through the $L (p) \times L(p)$ box centered at $mL(p)$ on the positive axis is bounded by a constant $c_3$ that does not depend of $p$.

Hence, the FKG inequality yields that
\begin {eqnarray*}
P_p ( 0 \leftrightarrow \partial \Lambda_{L(p)} )
& \ge & \theta (p) \\
&\ge  & P_p ( 0 \leftrightarrow \partial \Lambda_{L(p)} \hbox { and  } {\mathcal E}(L(p))) \\
& \ge & P_p ( 0 \leftrightarrow \partial \Lambda_{L(p)} ) \times P_p ({\mathcal E}(L(p))) \\
& \ge &c_3 P_p ( 0 \leftrightarrow \partial \Lambda_{L(p)}) .
\end {eqnarray*}
So, in order to control the behavior of $\theta(p)$, it is in fact sufficient to estimate
$P_p ( 0 \leftrightarrow \partial \Lambda_{L(p)})$.

\section {Outline of the proof}

Throughout this lecture, our goal will be to get bounds that hold uniformly for all $p \ge 1/2$ and $n \le L(p)$ (we just say ``uniformly for $n \le L(p)$''). The definition of $L(p)$ shows that
$ h_{p_0}( L(p_0)) $ is rather close to $1 - \epsilon$.
We will use Russo's identity (i.e. we will differentiate $h_p (n)$ with respect to $p$) to express
$h_{p_0} ( L(p_0)) - h_{1/2} ( L(p_0))$  in terms of the integral for $p=1/2$ to $p_0$ of a certain quantity.

One can similarly use Russo's formula to evaluate the derivative of $P_p (0 \leftrightarrow \partial \Lambda_n)$. In fact, we shall see that the same four-arm events govern both derivatives. This will roughly lead to the observation that
$$ \frac {d}{dp} h_p(n) \asymp \frac {d}{dp} \log P_p (0 \leftrightarrow \partial \Lambda_n )$$
uniformly for $n \le L(p)$. The definition of $L(p)$ will then imply that
$$  P_{p_0} ( 0 \leftrightarrow \partial \Lambda_{L(p_0)}) \asymp P_{1/2}(0 \leftrightarrow \partial \Lambda_{L(p_0)} )
= L(p_0)^{-5/48 + o(1)}$$
as $p_0 \to 1/2+$.

In order to evaluate $L(p)$, we will use Russo's formula again, but this time for a four-arm event. It will turn out that the
the ``four-arm up to distance $L(p)$'' event has a probability that is comparable for $P_p$ and for $P_{1/2}$. From this, we will be able to evaluate
$L(p)$ more precisely and relate its behavior to the four-arm event exponent at $p=1/2$.

Before to proceed to this proof, we need to collect some non-trivial facts: ``Uniform'' priori estimates for arm-exponents and ``uniform'' arm-separation lemmas.

\section {A priori estimates}
We note that we have uniform Russo-Seymour-Welsh estimates  for $n \le L(p)$.
Hence, those who remember the first exercise sheet will note that it implies uniform estimates for the probabilities of existence
of three-arms in the half-plane, of two-arms in the half-plane and of five-arms in the plane.

More precisely, the arm-exponents estimates that we shall use and that can be derived using the very same arguments as in the exercise sheet are the following:
Uniformly for $m \le n \le L(p)$,
\begin {itemize}
\item
The probability that there exist three disjoint arms with alternating colors joining in the upper half-plane the circles $\partial_m$ to $\partial_n$ is bounded from above by a constant times $(m/n)^2$.
\item
The probability that there exist five disjoint arms (two open and three closed) with ``alternating'' colors joining the circles $\partial_m$ and $\partial_n$ is bounded from below by a constant times $(m/n)^2$. We can even ask that some specified point on the boundary is between two ``adjacent arms'' of opposite color.
\end {itemize}

From the second estimate, we infer that the probability that there exist four arms with alternating color joining the two circles
$\partial_m$ and $\partial_n$ is bounded from below by a constant times $(m/n)^{2- \beta}$ for some $\beta>0$, uniformly for $m \le n \le L(p)$. Indeed, if we start exploring from a given boundary point of the annulus the first four arms exactly as in the previous lecture, the probability that in the remaining domain there exists (for $p>1/2$)  an additional closed arm is bounded from above by the probability of the same event at $p=1/2$ and therefore by some absolute constant times $(m/n)^\beta$.

Another estimate concerns the existence of one arm in a wedge of fixed angle. Suppose for instance that:
\begin {itemize}
\item
There exists a horizontal open crossing of the parallelogram $[0, 2n] \times [0,n]$.
\item
There exists a vertical open crossing of the parallelogram $[0, 2n] \times [0, 4n]$.
\item
There exists a horizontal open crossing of the parallelogram $[0, 8n] \times [0, 4n]$.
\item
$\cdots$
\item
There exists a vertical crossing of the parallelogram $[0, 2n \times 4^i] \times [0, 4^{i+1}\times n]$.
\end {itemize}
Then, there exists an open path in the wedge of angle $\pi/3$ that goes from distance $n$ of the origin to distance $n 4^i$.
If we combine this with the FKG inequality and our definition of $L(p)$, we see that there exists a constant $K_0$ such that  as long as $n4^{i+1} \le L(p)$,
the probability of this event is at least $4^{-i K_0}$.
Analogous arguments enable to derive the same result for other angles than $\pi/3$.

\section {Arm separation}

We will need to use again the idea of ``arm-separation'' in order to get a lower submultiplicative bound for the four-arm events.
One way to proceed is close in spirit to the  method used in the previous lecture (the method that we shall use here is in fact more general). We use the same probabilities  $\pi^\delta_4 (r_1, r_2), \hat \pi_4 (r_1, r_2)$ etc., but this time, they depend also on $p$. For notational convenience, we will drop the subscript $4$ (we will only work with the four-arms case in this section)  and denote these probabilities by $\hat \pi_p (r_1, r_2)$, $\pi^\delta_p (r_1, r_2)$. For convenience, we will work with hexagons again instead of discrete disks. So, the quantities $\pi$ will denote the probabilities for existence of arms joining certain discrete hexagons.

\begin {prop}
If $\delta$ has been chosen small enough, then for some  constant $c_0 = c_0( \epsilon, \delta)$, for all $4r_1 < r_2 \le L(p) $,
$$ \hat \pi_p (r_1, r_2) \le c_0 \pi^\delta_p (r_1, r_2) .$$
\end {prop}

Exactly as in the previous lecture, consequences (if one combines this with the RSW estimates) of this proposition go as follows:

\begin {corollary}
For some universal constant $c=c(\epsilon)$, for all $4r_1 < r_2 \le L(p)$,
$$ \hat \pi_p (r_1, r_2) \le c \hat \pi_p (2r_1, r_2/2).$$
\end {corollary}
\begin {corollary}
For some universal constant $c = c( \epsilon)$, for all $16 r_1 < 4 r_2 < r_3 \le L(p)$,
$$c \times  \hat \pi_p (r_1, r_2 ) \times \hat \pi_p (4r_2, r_3 ) \le \hat \pi_p (r_1, r_3).$$
\end {corollary}

Let us prove the proposition: We note that ($2^l \le  n$ and $l-j = 2k$)
\begin{eqnarray*}
\hat \pi_p (2^j, 2^l)
& \le & \pi_p^\delta (2^j, 2^l) + a(\delta) \hat \pi_p (2^{j+1}, 2^{l-1}) \\
& \le & \pi_p^\delta (2^j, 2^l) + a(\delta) \pi_p^\delta (2^{j+1}, 2^{l-1})+  a(\delta)^2 \hat \pi_p (2^{j+2}, 2^{l-2}) \\
& \le & \sum_{i=0}^k (a (\delta)^i \pi_p^\delta (2^{j+i}, 2^{l-i})) + a(\delta)^{k}.
\end {eqnarray*}
But,  by Russo-Seymour-Welsh and the a priori estimate for an arm in a wedge, one can prove that for some constant $c(\delta) <1$ and some exponent $K$,
$$ c(\delta) \times \pi_p^\delta (2^{j+i}, 2^{l-i}) \times  (2^i)^{-K} \le
\pi_p^\delta (2^j, 2^k)$$
i.e. that with probability that decays not faster than a power of $2^{-i}$, one can extend a configuration corresponding to $\pi^\delta_p (2^{j+i}, 2^{l-i})$ into a configuration corresponding to $\pi_p^\delta (2^j, 2^l)$.

Hence, it follows immediately that
$$
\hat \pi_p (2^j, 2^l) \le  \sum_{i=0}^{k} (c^{-1} \times a (\delta)^i \times
2^{iK}  \times \hat \pi_p (2^{j}, 2^{l}))
.$$
If we choose $\delta$ small enough so that $a(\delta) < 2^{-K}/2$, we get that
$$ \hat \pi_p (2^j, 2^l) \le (1/2c) \times \hat \pi^\delta_p (2^j, 2^l),$$
 uniformly with respect to $2^j, 2^l \le L(p)$.

The proposition (i.e. for any $r_1$ and $r_2$) follows easily by comparing the values of $\pi^\delta$ and $\hat \pi$ with
that of $r_1'= 2^j\in [r_1, 2r_1]$ and $r_2'= 2^l \in [r_2, 2r_2]$.

\section {Using differential inequalities}

\noindent
{\bf Using differential inequalities for crossing probabilities.}
We fix $n \le L(p)$, and we would like to see how the crossing probability $h_p (n)$ changes when $p$ increases a little bit. Note that the crossing probability increases with $p$.
If we couple the percolation with parameter $p$ and the percolation with parameter $p + dp$ (for very small/infinitesimal $dp$), we see that with a probability proportional to $dp$, exactly one site will differ, and that the probability that two sites differ is bounded by a constant times $dp^2$. If just one site $x$ differ, the crossing event can change only if this site $x$ is ``pivotal'' i.e. if there are four arms of alternating colors starting from its neighbors: two open ones going to the left and right boundaries of the $2n \times n$ parallelogram, and two closed ones to the top and bottom boundaries. Using the previous estimates and arguments analogous to the ones that we used in the previous paragraphs, we can conclude  that:

\begin {lemma} Uniformly for $n \le L(p)$,
$$\frac d {dp} h_p (n)  \asymp n^2 \hat \pi_p (n).$$
\end {lemma}

\noindent
{\bf Proof.}
Clearly
$$\frac d {dp} h_p (n) = \sum_x P_p ( x \hbox { is pivotal}). $$
Using the previous estimates, we see that the contribution of the $O(n^2)$ points $x$ that are at distance more than $n/4$ of the boundary of the parallelogram is at least $\hat \pi_p (n)$. This shows the lower bound for $d/dp (h_p (n))$.

For the upper bound, we have to show that the contributions due to those $x$'s that are close to the edges of the parallelogram do not matter much (one has to see that those points do not find it that much easier to be pivotal and there are less of them); we shall use a priori estimates of  probabilities of three arms in a half-plane or two arms in a wedge.

\begin{figure}
\centerline{\includegraphics*[height=2.3in]{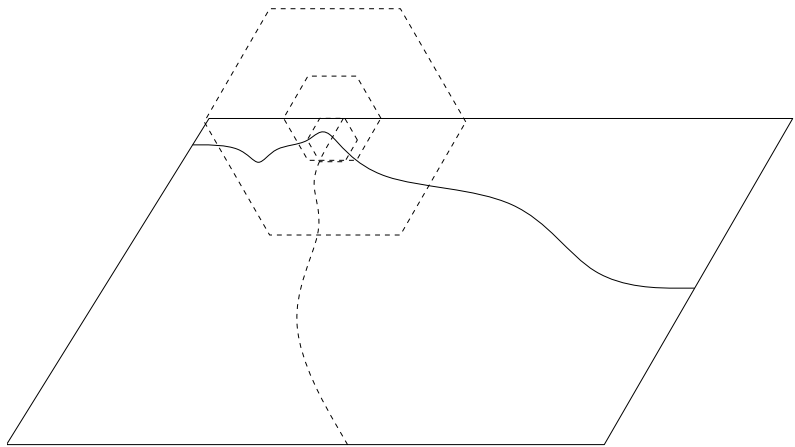}}
\caption {A pivotal point for $H(n)$}
\end {figure}

If we decompose this sum using coordinates centered at each of the four corners of the parallelogram, we see that it can be bounded by
$$ c \sum_{0 < i < j < n} \hat \pi_p (i) \pi_p^{3HP} (i, j) \pi_p^{2W} (j, n) $$
for some constant $c$, where $\pi_p^{3HP}$ and $\pi_p^{2W}$ denote respectively the three-arm probabilities in a half-plane and the two-arm probabilities in a wedge with angle $\pi /3$.
Recall the following facts (that we have just derived or that follow from a priori estimates in a half-plane):
For some $c>0$,  $\beta >0$, and all $i<j< n \le L(p)$,
\begin {itemize}
\item $\hat \pi_p (i ) \le c \hat \pi_p (n) \times (n/i)^{2- \beta}$.
\item $\pi_p^{3HP} (i,j ) \le c (i/j)^2$.
\item $\pi_p^{2W} (j, n) \le c (j/n)^{1 + \beta}$.
\end {itemize}
Hence,
$$ \frac d {dp} h_p (n) \le c' \hat \pi_p (n) \sum_{0<i<j<n} (n/i)^{2 - \beta} (i/j)^2 (j/n)^{1+ \beta}
\le c'' n^2 \hat \pi_p (n).$$
\qed

\medbreak

The following corollary will later on enable us to determine the asymptotic behavior of $L(p)$:
\begin {corollary} There exist absolute constants $c$, $c'$ and $p_1 > 1/2$ such that for all $p_0 \in (1/2,p_1)$,
$$ c \le L(p_0)^2 \int_{1/2}^{p_0}    \hat \pi_{p} ( L(p_0)) dp  \le c'. $$
\end {corollary}

\medbreak
\noindent
{\bf Proof.}
We integrate the identity of the previous lemma from $p=1/2$ to $p=p_0$ for $n = L(p_0) \le L(p)$. Note that the definition of
$L(p)$ shows that $h_{p_0} (L(p_0)) - h_{1/2} (L(p_0))$ is uniformly bounded away from $0$ and from infinity (at least for $p \in (1/2, p_1)$ for some $p_1<1$), so that
\begin {eqnarray*}
1 & \asymp & h_{p_0} ( L (p_0)) - h_{1/2} (L(p_0)) \\
& \asymp & \int_{1/2}^{p_0} \frac {d}{dp} h_p (L(p_0)) dp \\
& \asymp & \int_{1/2}^{p_0} L(p_0)^2 \hat \pi_p (L(p_0)) dp .
\end {eqnarray*}
\qed

\medbreak
\noindent
{\bf Using differential inequalities for the four arm event.}
We now use a similar argument to estimate the variation of $\hat \pi_p (n)$ with $p$, when $n$ is fixed (and not larger than $L(p)$).
Note that this time, the event $\hat \Pi_n$ corresponding to $\hat \pi_p (n)$ is not increasing anymore, but we still have an inequality:
$$  |\frac {d}{dp} \hat \pi_p (n) | \le \sum_x P_p ( x \hbox { is pivotal for } \hat \Pi_n ) .$$
Note that if $x$ is pivotal for $\hat \Pi_n$, then one has a four-arm event in each of the three ``annuli'' depicted in the following picture.

\begin{figure}
\centerline{\includegraphics*[height=2.in]{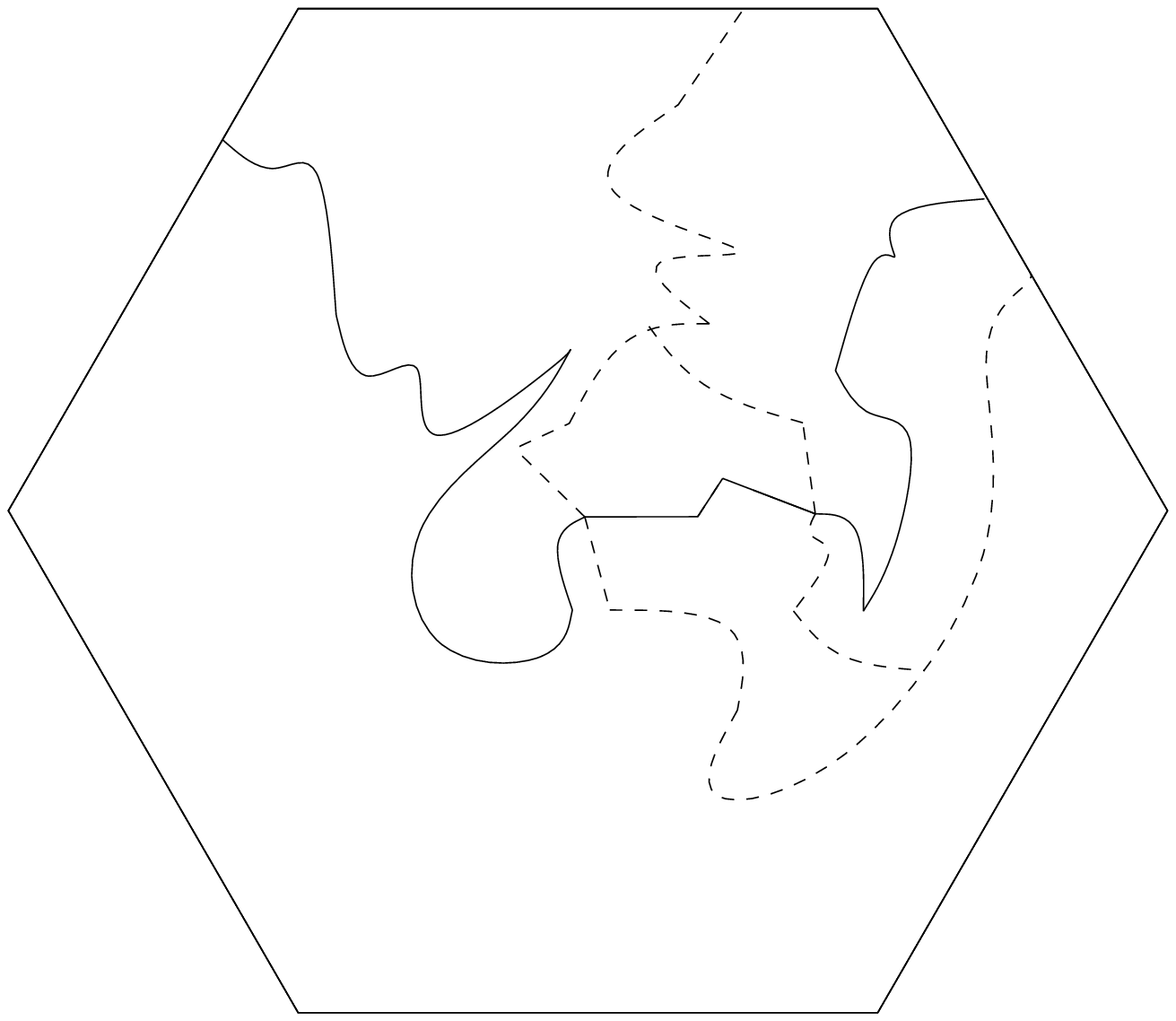}}
\centerline{\includegraphics*[height=2.in]{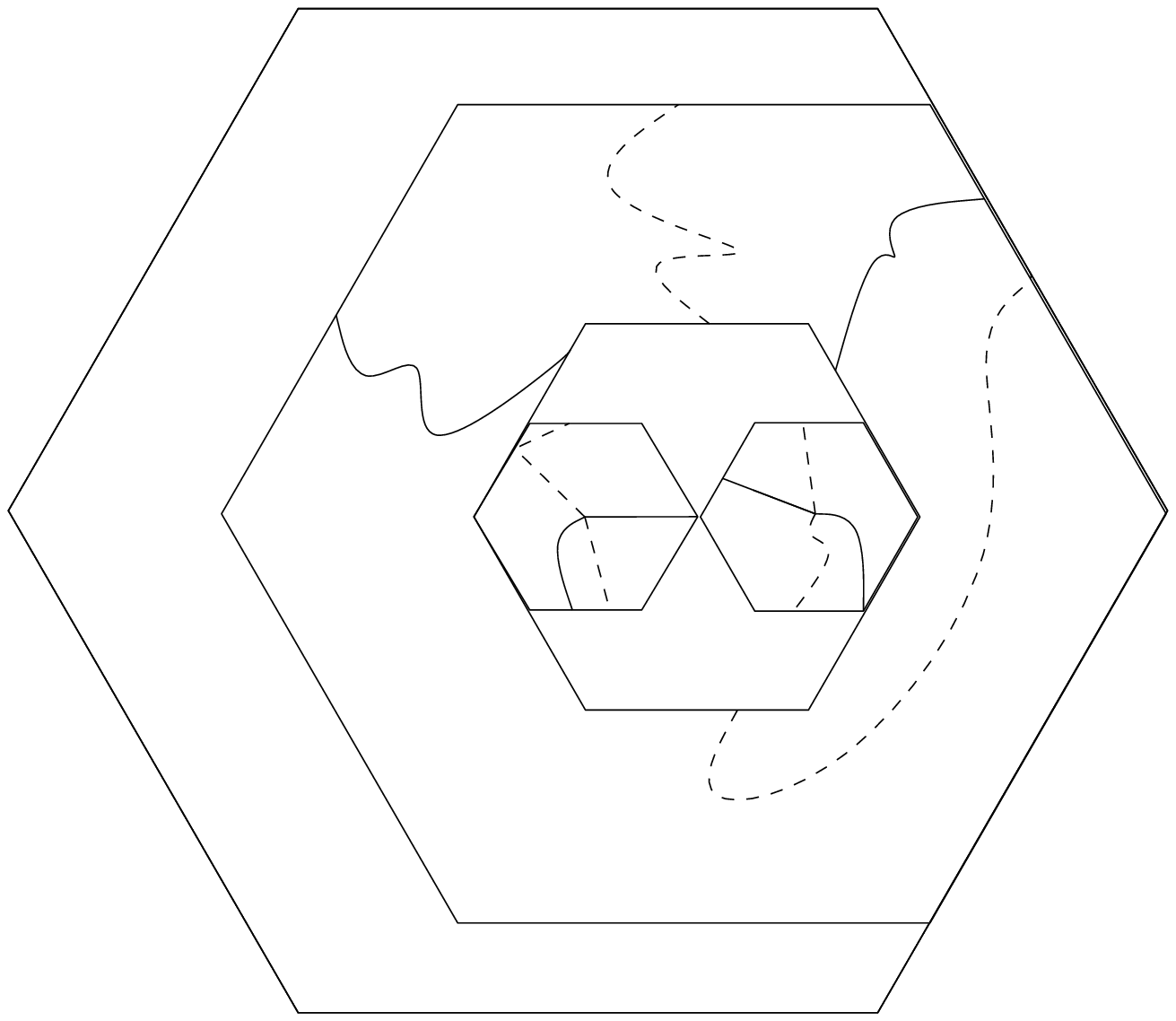}}
\caption {A pivotal point for the four arm event and the three corresponding four-arm events}
\end {figure}

In particular, using the same cut-and-pasting arguments (and the separation lemma is again instrumental) as in the previous paragraphs, we end up writing inequalities of the following type (to write this rigorously, one would need again to show that the contributions corresponding to the case where $x$ is close to the boundary of the hexagon do not matter):
\begin{eqnarray*}
|\frac {d}{dp} \hat \pi_p (n) | & \le &c  \sum_x
\hat \pi_p (\| x \| / 2)^2 \hat \pi_p ( 2\| x \|, n ) \\
& \le & c' \sum_x
\hat \pi_p (\| x \| / 2) \hat \pi_p (  n ) \\
\\
& \le &
c'' \hat \pi_p (n) \times \frac {d}{dp} h_p (n)
\end {eqnarray*}
and
$$
| \frac {d}{dp} \log \hat \pi_p(n)  | \le cst \frac {d}{dp} h_p (n).
$$
If we integrate this relation from $p=1/2$ to $p'$ for $n= L(p_0)$ with $p_0 > p'> 1/2$, we get the following lemma:
\begin {lemma}
Uniformly for $p' \in (1/2, p_0)$,
$$ \hat \pi_{p'} (L(p_0)) \asymp \hat \pi_{1/2} (L(p_0)).$$
\end {lemma}

If we furthermore combine this with the corollary derived in the previous paragraph, we get that
$$
1 \asymp
\int_{1/2}^{p_0} dp L(p_0)^2 \hat \pi_{1/2} ( L(p_0)) = (p_0 -1/2) \times L(p_0)^2 \times \hat \pi_{1/2} (L(p_0)).$$
In other words,
$$
 L(p_0)^2 \times \hat \pi_{1/2} (L(p_0)) \asymp (p_0 - 1/2)^{-1}.$$
If we plug in the value of the four-arm exponent at $p=1/2$ derived in the previous lecture, we get
the following corollary:
\begin {corollary}
\label {final3} 
When $u\to 0+$, one has
 $$L(1/2 + u) = u^{-4/3 + o(1)}.$$
\end {corollary}

\medbreak
\noindent
{\bf Using differential inequalities for the one-arm event.}
 The same argument can also be adapted (see the figure) to show that
$$\left|\frac {d}{dp} \log P_p ( 0 \leftrightarrow \partial \Lambda_n ) \right| \le c  n^2 \hat \pi_p (n)
\asymp n^2  \hat \pi_{1/2} (n).$$
\begin{figure}
\centerline{\includegraphics*[height=2.4in]{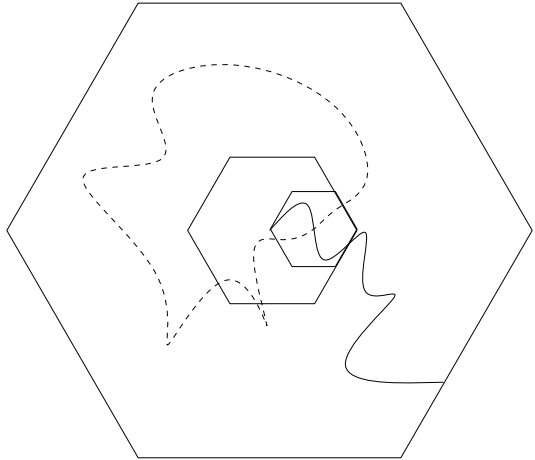}}
\caption {A pivotal point for $0 \leftrightarrow \partial \Lambda_n$, the corresponding four-arm event and one-arm event}
\end {figure}
It follows that
$$P_p ( 0 \leftrightarrow\partial \Lambda_n ) \asymp P_{1/2}(0 \leftrightarrow\partial \Lambda_n )  $$
for $n \le L(p)$.
If we plug in the value of the one-arm exponent derived in the previous lectures, we get the following lemma:

\begin {lemma}
\label {final4}
When $p \to 1/2 +$, one has
$$ P_p ( 0 \leftrightarrow \partial \Lambda_{L(p)} )
= L(p)^{-5/48 + o (1)}.$$
\end {lemma}

\bigbreak
\noindent
{\bf End of the proof of the theorem.}
Finally, combining  Corollary \ref {final3} and Lemma \ref {final4}, we get that when $p \to 1/2+$,
$$\theta (p)  \asymp P_{1/2} (0 \leftrightarrow \Lambda_{L(p)})
= L(p)^{-5/48+ o(1)} = (p-1/2)^{5/36+ o(1)}.
$$
This concludes the proof of Theorem \ref {5/36}.
Note that in order to prove this result, we have combined the input of all six lectures!

\section {Concluding remarks}
\begin {itemize}
\item
These questions on near-critical percolation are related to how percolation events depend on fluctuations or noising of the discrete realization.  For results and references in this direction, see e.g. \cite {SchSt,GPS}.
\item
It is still an open problem to generalize these Kesten's scaling relations to dependent models such as the Ising model. This would be nice since the other steps (conformal invariance, convergence to SLE) have now been rigorously established by Smirnov \cite {Sm2,Sm3,Sm4}.
\item
Another big open question is the derivation of conformal invariance of percolation for other planar lattices.
The interested reader may read \cite {Be3} for a discussion of this issue.
\end {itemize}

\eject

\end{document}